\documentclass[11pt]{article}
\usepackage{amsfonts}
\usepackage{amsmath}
\usepackage{amssymb}
\newcommand{\be}{\begin{equation}}
\newcommand{\ee}{\end{equation}}
\newcommand{\bea}{\begin{eqnarray}}
\newcommand{\eea}{\end{eqnarray}}

\newcommand{\xa}{\alpha}
\newcommand{\xb}{\beta}
\newcommand{\xg}{\gamma}
\newcommand{\xG}{\Gamma}
\newcommand{\xd}{\delta}
\newcommand{\xD}{\Delta}
\newcommand{\xe}{\varepsilon}

\newcommand{\xz}{\zeta}

\newcommand{\xm}{\mu}
\newcommand{\xn}{\nu}

\newcommand{\xs}{\sigma}

\newcommand{\xf}{\phi}

\newcommand{\xO}{\Omega}

\newcommand{\Ren}{ I \! \! R^N}

\newtheorem{theorem}{Theorem}[section]
\newtheorem{lemma}[theorem]{Lemma}
\newtheorem{prop}[theorem]{Proposition}
\newtheorem{coro}[theorem]{Corollary}
\newtheorem{defin}[theorem]{Definition}
\newtheorem{notation}[theorem]{Notation}
 
\newcommand{\qcap}{C_{\frac{2}{q},\;q'}}
\newcommand{\qsub}{\subset^q}
\newcommand{\qeq}{\sim^q}
\newcommand{\besl}{W^{\frac{2}{q},q'}(\mathbb{R}^N)\cap L^\infty(\mathbb{R}^N)}

\newtheorem{subn}{\name}

\newcommand{\bsn}[1]{\def\name{#1}
\begin{subn}}
\newcommand{\esn}{
\end{subn}}
\newtheorem{sub}{\name}[section]
\newcommand{\bs}{
\begin{sub}}
\newcommand{\es}{
\end{sub}}
\newcommand{\bsl}[1]{
\begin{sub}\label{#1}}
\newcommand{\bth}[1]{\def\name{Theorem}
\begin{sub}\label{t:#1}}
\newcommand{\blemma}[1]{\def\name{Lemma}
\begin{sub}\label{l:#1}}
\newcommand{\bcor}[1]{\def\name{Corollary}
\begin{sub}\label{c:#1}}
\newcommand{\bdef}[1]{\def\name{Definition}
\begin{sub}\label{d:#1}}
\newcommand{\bprop}[1]{\def\name{Proposition}
\begin{sub}\label{p:#1}}

\newcommand{\BA}{
\begin{array}}
\newcommand{\EA}{
\end{array}}
\newcommand{\BAN}{\renewcommand{\arraystretch}{1.2}
\setlength{\arraycolsep}{2pt}
\begin{array}}
\newcommand{\BAV}[2]{\renewcommand{\arraystretch}{#1}
\setlength{\arraycolsep}{#2}
\begin{array}}
\newcommand{\BSA}{
\begin{subarray}}
\newcommand{\ESA}{
\end{subarray}}
\newcommand{\BAL}{
\begin{aligned}}
\newcommand{\EAL}{
\end{aligned}}
\newcommand{\BALG}{
\begin{alignat}}
\newcommand{\EALG}{
\end{alignat}}
\newcommand{\BALGN}{
\begin{alignat*}}
\newcommand{\EALGN}{
\end{alignat*}}
\newcommand{\note}[1]{\textit{#1.}\hspace{2mm}}
\newcommand{\Proof}{\note{Proof}}

\newcommand{\Remark}{\note{Remark}}

\newcommand{\abs}[1]{\left |#1\right |}
\newcommand{\norm}[1]{\left \|#1\right \|}

\def\angb<#1>{\langle #1 \rangle}

\newcommand{\opname}[1]{\mbox{\rm #1}\,}

\newcommand{\dist}{\opname{dist}}
\newcommand{\myfrac}[2]{{\displaystyle \frac{#1}{#2} }}
\newcommand{\myint}[2]{{\displaystyle \int_{#1}^{#2}}}

\newcommand {\dint}{{\displaystyle \int\!\!\int}}

\newcommand{\prt}{
\partial}

\newcommand{\ti}{\times}

\def\ga{\alpha}            
             \def\ge{\epsilon}

\def\gm{\mu}        \def\gn{\nu}         \def\gp{\pi}
            
\def\gs{\sigma}       
      
\def\gx{\xi}                \def\gz{\zeta}
     \def\Gd{\Delta}      

    \def\Gs{\Sigma}      
\def\Gw{\Omega}     \def\Gx{\Xi}         

\def\CS{{\mathcal S}}      
\def\CR{{\mathcal R}}   \def\CO{{\mathcal O}}

   \def\CU{{\mathcal U}}

   \def\BBH {\mathbb H}    
       
   \def\BBN {\mathbb N}    
   \def\BBR {\mathbb R}

\medskip

\title{Classification of positive solutions of heat equation \\with supercritical absorption}

\author{\bf  Konstantinos T. Gkikas\footnote{ kugkikas@gmail.com}
                                                                           \\
       \small Centro de Modelamiento Matem\`atico\\
        \small  Universidad de Chile, Santiago de Chile, Chile
                                   \\\\
         \bf Laurent V\'{e}ron \footnote{ veronl@univ-tours.fr}\\
\small Laboratoire de Math\'{e}matiques et Physique Th\'{e}orique\\
  \small       Universit\'{e} Fran\c{c}ois-Rabelais, Tours, France
         }
\date{}
\begin{document}
\maketitle
\begin{abstract}
Let $q\geq 1+\frac{2}{N}$. We prove that any positive solution of (E) $\prt_t u-\xD u+u^q=0$ in $\mathbb{R}^N\times(0,\infty)$ admits an initial trace which is a nonnegative Borel measure, outer regular with respect to the fine topology associated to the Bessel capacity $C_{\frac{2}{q},q'}$ in $\BBR^N$ ($q'=q/q-1)$) and absolutely continuous with respect to this capacity. If $\nu$ is a nonnegative Borel measure in $\BBR^N$ with the above properties we construct a positive solution $u$ of (E) with initial trace $\gn$ and we prove that this solution is the unique $\gs$-moderate solution of (E) with such an initial trace. Finally we prove that every positive solution of (E) is $\gs$-moderate.
\end{abstract}
\tableofcontents\medskip

  \vspace{1mm}

\medskip

  \noindent {\small {\bf Aknowledgements} This research has been made possible thanks to the support of the R\'egion Centre which offered a post-doctoral position for the first author during the year 2011-2012.}

\newpage

\section{Introduction}
Let $q>1$, $Q_T=\BBR^N\times(0,T)$ with $T>0$ and $Q=\BBR^N\times(0,\infty)$. It is proved by Marcus and V\'eron \cite{MV-CPDE} that for any positive function $u\in C^{2,1}(Q_T)$ solution of
\begin{equation}\label{E1}
\prt_t  u-\Delta u+u^q=0
\end{equation}
there exists  a unique couple $(\CS,\gm)$ where $\CS$ is a closed subset  of $\BBR^N$ and $\gm$ a positive Radon measure  on $\CR:=\BBR^N\setminus\CS$  such that
\begin{equation}\label{E2}
\lim_{t\to 0}\myint{\CO}{}u(x,t) dx=\infty
\end{equation}
for all open set $\CO$ of $\BBR^N$ such that $\CS\cap\CO\neq\emptyset$, and
\begin{equation}\label{E3}
\lim_{t\to 0}\myint{\BBR^N}{}u(x,t)\gz(x) dx=\myint{\BBR^N}{}\gz(x) d\gm(x)\qquad\forall \gz\in C^\infty_0(\CR).
\end{equation}
To this couple $(\CS,\gm)$ it is associated a unique outer Borel measure $\gn$ called {\it the initial trace} of $u$ and denoted by  $tr(u)$. The set $\CS$ is {\it the singular set} of $\gn$ and the measure $\gm$ is {\it the regular set} of $\gn$.  Conversely, to any outer Borel measure $\gn$ we can associate its singular part $\CS(\gn)$ which is a closed subset of $\BBR^N$ and its regular part $\gm_\gn$ which is a positive Radon measure on $\CR(\gn)$. We denote $\gn\approx (\CS,\gm)$. When $1<q<q_c:=\frac{N+2}{N}$ Marcus and V\'eron \cite{MV-CPDE} proved that the trace operator $tr$ defines a one to one correspondence between the set $\CU_+(Q_T)$ of positive solutions of (\ref{E1}) in $Q_T$ and the set $\mathfrak B^{reg}(\BBR^N)$ of positive outer Borel measures in $\BBR^N$. This no longer the case if $q\geq q_c$ since not any closed subset of $\BBR^N$ (resp. any positive Radon measure) is eligible for being the singular set (resp. the regular part) of the the initial trace of some positive solution of (\ref{E1}). It is proved in \cite{BP-parab} that the initial value problem
\begin{equation}\label{E4}\BA {lll}
\prt_t  u-\Delta u+|u|^{q-1}u=0\qquad&\text{ in }Q\\
\phantom{\prt_t  u+u^q--}
u(.,0)=\gm&\text{ in }\BBR^N
\EA\end{equation}
where $\gm$ is a positive bounded Radon measure admits a solution if and only if $\gm$ satisfies
\begin{equation}\label{E5}
C_{\frac{2}{q},q'}(E)=0\Longrightarrow \gm(E)=0\quad\forall E\subset\BBR^N, E\text{ Borel},
\end{equation}
where $C_{\frac{2}{q},q'}$ stands for the Bessel capacity in $\BBR^N$ ($q'=q/(q-1)$). It is shown in \cite{MV-CPDE} that this result holds even if $\gm$ is unbounded; this solution is unique and denoted $u_\gm$. If $G$ is a Borel subset of $\BBR^N$ we denote by $\mathfrak M_q(G)$ the set of Borel measures $\gm$ in $G$ with the property that
\begin{equation}\label{E5'}
C_{\frac{2}{q},q'}(E)=0\Longrightarrow \gm(E)=0\quad\forall E\subset G, E\text{ Borel},
\end{equation}
 In the same article it is proved that a necessary and sufficient condition in order $\gn\approx (\CS,\gm)$ to be the initial trace of a positive solution of (\ref{E1}) is
\begin{equation}\label{E6}
\gm\in \mathfrak M_q(\CR)
\end{equation}
and
\begin{equation}\label{E7}
\CS=\prt_\gm\CS\bigcup \CS^*
\end{equation}
where
\begin{equation}\label{E8}
\prt_\gm\CS=\{z\in\CS:\gm(B_r(z)\bigcap\CS)=\infty\,,\;\forall r>0\}
\end{equation}
and
\begin{equation}\label{E9}
\CS^*=\{z\in\CS:C_{\frac{2}{q},q'}((B_r(z)\bigcap\CS)>0\,,\;\forall r>0\}.
\end{equation}
The meaning of (\ref{E7}) is that the singular set is created either by the local unboundedness of the Radon measure or because the singular set is localy non-removable. Furthermore the solution which is constructed is the maximal solution with initial trace $(\CS,\gm)$. 

A striking result due to Le Gall \cite{LeG3} shows that if $q=2$ and $N\geq 2$, a positive
 solution of (\ref{E1}) is not uniquely determinef by its initial trace $\gn\approx (\CS,\gm)$ if $\CS\neq\emptyset$. The results is actually extended to any $q\geq q_c$ in \cite{MV-CPDE}. The main point in this counter-example relies on the construction of a positive solution $u$ of (\ref{E1}) with a singular set $\CS=\BBR^N$,  with a blow-up set at $t=0$ which the union of a countable of 
 closed balls $\overline B_{\ge_n}(a_n)$ where $\{a_n\}$ is a dense set in $\BBR^N$ and the $\ge_n$ are chosen small enough 
so that $u(0,1)\leq \ga$ for some $\ga>0$ fixed. If $U_{\overline B_{\ge_n}(a_n)}$ denotes the solution with initial trace $(\overline B_{\ge_n}(a_n),0)$, then $U_{\overline B_{\ge_n}(a_n)}(0,1)\leq C(\ge_n)$ with $\lim_{\ge\to 0}C(\ge).$ 
This is a consequence of the supercriticality assumption and the estimates in \cite{MV-CVPDE}. The solution $u$ is constructed between a sub-solution and a super-solution
  \begin{equation}\label{E90}
\sup_n\{U_{\overline B_{\ge_n}(a_n)}\}\leq u\leq\sum_{n=0}^\infty U_{\overline B_{\ge_n}(a_n)}, 
\end{equation}
the right-hand side being chosen so that $\sum_{n=0}^\infty C(\ge_n)\leq \ga$.  Denoting $E=\cup_n\overline B_{\ge_n}(a_n)$, then $\abs E<\infty$ and $u$ satisfies
   \begin{equation}\label{E91}
\lim_{t\to 0}u(x,t)=0\qquad\forall x\in\BBR^N\setminus E\text{ where }\abs E<\infty,
\end{equation}
and
   \begin{equation}\label{E92}
\lim_{t\to 0}t^{\frac{1}{q-1}}u(x,t)=c_q=(q-1)^{\frac{1}{1-q}}\quad\text{ uniformly for } x\in K\subset\bigcup_n B_{\ge_n}(a_n),\; K \text{ compact}.
\end{equation}
Thus (\ref{E2}) holds for any nonempty open set $\CO\subset\BBR^N$. This counter-example points out that the trace process associated to averaging a positive solution $u$ of (\ref{E1}) on open sets and letting $t\to 0$ is not sharp enough to distinguish among solutions; this process is now called the {\it rough trace}. This is why the introduction of a finer averaging appears to be needed. This finer averaging method is constructed by using the  {\it fine topology} associated to the capacity $C_{\frac{2}{q},q'}$. It will lead us to the notion of  precise trace.

 \medskip

 A similar approach has been carried out if one considers the boundary trace problem for the positive solutions of the elliptic equation
 \begin{equation}\label{El1}
-\Gd u+|u|^{q-1}u=0\qquad\text {in }\Gw
\end{equation}
where $\Gw$ is a bounded $C^2$ domain in $\BBR^N$ ($N\geq 2$) and $q>1$. The boundary trace is defined in a somewhat similar way as the initial trace, by considering the limit in the weak sense of measures, of the restriction of $u$ to the set
$\Gs_\ge:=\{x\in\Gw:\dist (x,\Gw^c=\ge)\}$, when $\ge\to 0$. The boundary trace $tr_{\prt\Gw}(u)$ is a uniquely determined outer regular Borel measure on $\prt\Gw$, with singular part $\CS$, a closed subset of $\prt\Gw$ and regular part $\gm$, a positive Radon measure on $\CR=\prt\Gw\setminus \CS$. This equation possesses a critical exponent  $q_e=(N+1)/(N-1)$. The main contributions which lead to a complete picture of the boundary trace problem over a period of twenty years are due to Gmira and V\'eron \cite{GV-DJM}, Le Gall \cite{LeG1}, \cite{LeG2}, Dynkin and Kuznetsov \cite {Dyn},\cite{Dyn1}, \cite{DK1}
\cite{DK2}, \cite{DK3},\cite{Kuz}, Marcus and V\'eron \cite{MV-ARMA},\cite{MV-JMPA1},\cite{MV-JMPA2},\cite{MV-JEMS},\cite{MV-CONT}, \cite{MV-CVPDE}, \cite{M-JAM}, and Mselati \cite{Mse}. These contributions can be summarized as follows:\smallskip

\noindent (i) If $1<q<q_e$ the boundary trace operator establishes a one to one correspondence between the set $\CU_+(\Gw)$ of positive solutions of (\ref{El1}) and the set of positive outer regular Borel measures on $\prt\Gw$.\smallskip

\noindent (ii) If $q\geq q_e$ the boundary value problem

 \begin{equation}\label{El2}\BA {ll}
-\Gd u+|u|^{q-1}u=0\qquad&\text { in }\Gw\\
\phantom{-\Gd u+|u|^{q-1}}
u=\gm&\text { in }\prt\Gw
\EA\end{equation}
where $\gm$ is a positive Radon measure on $\prt\Gw$ admits a solution (always unique) if and only if
\begin{equation}\label{El3}
C_{\frac{2}{q},q'}(E)=0\Longrightarrow \gm(E)=0\quad\forall E\subset\prt\Gw, E\text{ Borel},
\end{equation}
where $C_{\frac{2}{q},q'}$ is the Bessel capacity in $\BBR^{N-1}$.\smallskip

\noindent (iii) If $q\geq q_e$, a outer regular Borel measure $\gn\approx (\CS,\gm)$ on $\prt\Gw$ is the boundary trace of a positive solution of (\ref{El1}) if and only if
$$
C_{\frac{2}{q},q'}(E)=0\Longrightarrow \gm(E)=0\quad\forall E\subset\CS, E\text{ Borel},
$$
and (\ref{E7}) holds with (\ref{E8}) and (\ref{E9}) where the capacity is relative to dimension $N$-$1$.\smallskip

\noindent (iv) If $q\geq q_e$ a solution is not uniquely determined by its boundary trace whenever $\CS\neq\emptyset$.\smallskip

However in \cite{MV-CONT} Marcus and V\'eron have defined a notion of {\it precise trace} for the case $q\geq q_e$ with the following properties,\smallskip

\noindent (v) If we denote by $\frak T_q$ the fine topology of $\prt\Gw$ associated with the  $C_{\frac{2}{q},q'}$-capacity, there exists a $\frak T_q$-closed subset $\CS_q$ of $\prt\Gw$ such that for every $z\in \CS_q$
\begin{equation}\label{El4}
\lim_{\ge\to 0} \myint{\Gx}{}u(\ge,\gs)dS=\infty
\end{equation}
for every $\frak T_q$-open neighborhood $\Gx$ of $z$ where $(r,\gs)\in [0,\ge_0]\ti\prt\Gw$ are the flow coordinates near
$\prt\Gw$, and for every $z\in \CR_q:=\prt\Gw\setminus \CS_q$, there exists a $\frak T_q$-open neighborhood $\Gx$ of $z$ such that
\begin{equation}\label{El4'}
\limsup_{\ge\to 0} \myint{\Gx}{}u(\ge,\gs)dS<\infty.
\end{equation} \smallskip

\noindent (vi) There exists a nonnegative Borel measure $\gm$ on $\CR_q$, outer regular for the $\frak T_q$-topology, such that
\begin{equation}\label{El5}
\lim_{\ge\to 0}u^{\Gx}_\ge=u_{\gm\chi_{_\Gx}}\qquad\text{locally uniformly in }\Gw,
\end{equation}
where $u^{\Gx}_\ge$ is the solution of
\begin{equation}\label{El6}
\BA {ll}
-\Gd v+|v|^{q-1}v=0\qquad&\text{in }\Gw_\ge:=\{x\in\Gw:\dist(x,\prt\Gw)>\ge\}\\
\phantom{-\Gd v+|v|^{q-1}}
v=u(\ge,.)\chi_{_\Gx}\qquad&\text{in }\Gs_\ge=\prt\Gw_\ge.
\EA
\end{equation}
The couple $(\CS_q,\gm)$ is uniquely determined and it is called {the precise boundary trace } of $u$. It can also be represented by a Borel measure with the $\frak T_q$-outer regularity. It is denoted by $tr^q_{\prt\Gw}(u)$.

Concerning uniqueness Dynkin and Kuznetsov introduced in \cite{DK3} the notion of $\sigma$-moderate solutions, which are elements $u$ of $\CU_+(\Gw)$ with the property that there exists an increasing sequence $\{\gm_n\}$ of nonnegative Radon measures on $\prt\Gw$ such that $u_{\gm_n}\to u$ when $n\to\infty$. In \cite{MV-CONT} Marcus and V\'eron proved that a $\sigma$-moderate positive solution of (\ref{El1}) is uniquely determined by its precise boundary trace. This precise trace is essentially the same, up to a set of zero $C_{\frac{2}{c},q'}$-capacity, as the {\it fine trace} that Dynkin and Kuznetsov introduced in \cite{DK3} using probabilistic tools such as the Brownian motion; however their construction is only valid in the range $(1,q]$ of values of $q$.  Finally, in
\cite{M-JAM}, Marcus proved that any positive solution is $\sigma$-moderate. Notice that this result was already obtained by Mselati \cite {Mse} in the case $q=2$ and then by Dynkin \cite{Dyn1} for $q_e\leq q\leq 2$  by using a combination of analytic and probabilistic techniques.\smallskip

In this article we define a notion of {\it precise initial trace} for positive solutions of (\ref{E1}) associated to the $\frak T_q$-topology, which denotes the $C_{\frac{2}{q},q'}$ fine topology of  $\BBR^N$.  We denote by $\BBH[.]$  the heat potential in $Q$ expressed by
\begin{equation}\label{Ep1}
\BBH[\gx](x,t)=\myfrac{1}{(4\gp t)^{\frac{N}{2}}}\myint{\BBR^N}{}e^{-\frac{\abs{x-y}^2}{4t}}\gx(y)dy,
\end{equation}
for all $\gx\in L^1(\BBR^N)$.
We define the {\it singular set } of $u\in\CU_+(Q_T)$
as the set of $z\in\BBR^N$ such that for any $\frak T_q$-open neighborhood $\CO\subset\BBR^N$ of $z$, there holds
\begin{equation}\label{Ep2}
\BA {ll}
\dint_{\!\!Q_T}\BBH[\chi_{_\CO}]u^qdxdt=\infty.
\EA
\end{equation}
The singular set, denoted by $\CS_q=\CS_q(u)$, is $\frak T_q$-closed. The regular set is $\CR_q:=\BBR^N\setminus\CS_q$; it is $\frak T_q$-open. If $z\in \CS_q$ and $\CO\subset\BBR^N$ is a $\frak T_q$-open neighborhood of $z$ such that
\begin{equation}\label{Ep3}
\BA {ll}
\dint_{\!\!Q_T}\BBH[\chi_{_\CO}]u^qdxdt<\infty,
\EA
\end{equation}
then for any $\eta\in L^\infty\cap W^{\frac{2}{q},q'}(\BBR^N)$ with $\frak T_q$-support contained in $\CO$ there exists
\begin{equation}\label{Ep4}
\BA {ll}
\lim_{t\to 0}\myint{\BBR^N}{}u(x,t)(\eta(x))^{2q'}dx:=\ell_\CO (\eta).
\EA
\end{equation}

As a consequence there exists a positive Borel measure $\gm$ on $\CR_q$, outer regular for the $\frak T_q$-topology, such that for  $\frak T_q$-open subset $\Gx\subset \CR_q$ there holds
\begin{equation}\label{Ep5}
\BA {ll}
\lim_{\ge\to 0}u_{\ge,\chi_{\Gx}}(.,t)=u_{\chi_{_\Gx}\gm}
\EA
\end{equation}
where $u_{\ge,\chi_{\Gx}}$ is the solution of
\begin{equation}\label{Ep6}
\BA {ll}
\prt_t  v-\Gd v+\abs v^{q-1}v=0\qquad&\text{ in } Q^\ge:=\BBR^N\ti (\ge,\infty)\\
\phantom{va+\abs v^{q-1}v}
v(.,\ge)=\chi_{_\Gx}&\text{ in } \BBR^N.
\EA
\end{equation}
The set $(\CS_q,\gm)$  is called the {\it precise initial  trace of } $u$ and denoted by $\mathrm{tr}^c(u)$. To this set we can associate a Borel measure $\gn$ on $\BBR^N$. It is absolutely continuous with respect to the $C_{\frac{2}{q},q'}$-capacity in the following sense
\begin{equation}\label{Ep6'}
\BA {ll}
\forall Q\subset\BBR^N, \frak T_q\,\text {- open }, \forall A\subset\BBR^N, \,A \text { Borel },
C_{\frac{2}{q},q'}(A)=0\Longrightarrow \gm(Q\setminus A)=\gm (Q).
\EA
\end{equation}
It is also outer regular with respect to the $\frak T_q$-topology in the sense that for every Borel set $E\subset\BBR^N$
\begin{equation}\label{Ep6''}
\BA {ll}
\gm(E)=\inf\{\gm(Q):Q\supset E,\,Q\;\frak T_q\,\text {- open }\}=\sup\{\gm(K):K\subset E,\,K\,\text { compact }\}.
\EA
\end{equation}
A measure with the above properties is called {\it $\frak T_q$-perfect}.
Similarly to Dynkin, we say that a positive solution $u$ of (\ref{E1}) is $\sigma$-moderate if the exists an increasing sequence $\{\gm_n\}$ of nonnegative Radon measures in $\BBR^N$ such that $u_{\gm_n}\to u$ when $n\to\infty$. It is proved in \cite{MV-CVPDE} that if $F\subset\BBR^N$ is a closed subset, the maximal solution $U_F$ with initial trace
$(F,0)$ coincides with the maximal  $\sigma$-moderate solution $V_F$ with the same trace and which is defined by
\begin{equation}\label{Ep7}
\BA {ll}
V_F=\sup\{u_\gm:\gm\in\frak M_q(\BBR^N), \gm(F^c)=0\}.
\EA
\end{equation}
It 
is indeed $\sigma$-moderate. Following Dynkin we define an addition among the elements of $\CU_+(Q_T)$ by
\begin{equation}\label{Ep8}
\BA {ll}
\forall (u,v)\in \CU_+(Q_T)\times\CU_+(Q_T),\, u\oplus v\text { is the largest element of $\CU_+(Q_T)$ dominated by }\, u+v.
\EA
\end{equation}
The main results of this article are the following\medskip

\noindent{\bf Theorem A}. {\it If $\gn$ is a $\frak T_q$-perfect measure with singular part $\CS_q$ and regular part $\gm$ on $\CR_q$
then $u_\gm\oplus U^{\CS_q}$ is the only $\sigma$-moderate element of $\CU_+(Q)$ with precise trace $\gn$.}
\medskip

In order to extend Marcus's result we need a parabolic counterpart of Ancona's characterization of positive solutions of Schr\"odinger equation with singular potential \cite{A-ANN}. We  prove a representation theorem valid for any positive solution of
\begin{equation}\label{Ep9}
\BA {ll}
\prt_tu-\Gd u+V(x,t)u=0\qquad\text{in }\,Q,
\EA
\end{equation}
where $V$ is a Borel function which satisfies, for some $c\geq 0$,
\begin{equation}\label{Ep10}
\BA {ll}
0\leq V(x,t)\leq \myfrac{c}{t}\qquad\text{ for almost all }\,(x,t)\in Q.
\EA
\end{equation}
Let $T$ be fixed and let $\psi$  be defined by
 $$\psi(x,t)=\myint{t}{T}\myint{\BBR^N}{}\myfrac{1}{(4\gp(s-t))^{\frac{N}{2}}}e^{-\frac{\abs{x-y}^2}{4(s-t)}}V(y,s)dy ds
 \qquad\text{in }Q_T.
 $$
 \noindent{\bf Theorem B}.{\it
There exists a kernel $\Gamma$ defined in $Q_T\times Q_T$ satisfying
\begin{equation}\label{Ep11}
\BA {ll}
c_1\frac{e^{-a_1\frac{\abs{x-y}^2}{s-t}}}{(t-s)^{\frac{N}{2}}}\leq \Gamma(x,t,y,s)\leq
c_2\frac{e^{-a_2\frac{\abs{x-y}^2}{s-t}}}{(t-s)^{\frac{N}{2}}}\,\quad\forall (x,t), (y,s)\in Q_T\times Q_T\;\text{with }s\leq t.
\EA
\end{equation}
where the $a_j$ and $c_j$ are positive contants depending on $T$ and $V$, such that for any positive solution $u$ of (\ref{Ep9}), there exists a positive Radon measure $\gm$ in $\BBR^N$ such that
\begin{equation}\label{Ep12}
\BA {ll}
u(x,t)=e^{\psi(x,t)}\myint{\BBR^N}{}\Gamma(x,t,y,0)d\gm(y)\qquad\text{ for almost all }\,(x,t)\in Q_T.
\EA
\end{equation}
}\medskip

The next result, combined with Theorem A,  shows that in the case $q\geq q_c$ the precise trace operator realizes a one to one correpondence between the set of positive solutions of (\ref{E1}) and the set of $\frak T_q$-perfect Borel measures in $\BBR^N$.\medskip

\noindent{\bf Theorem C }{\it Any positive solution of (\ref{E1}) is $\sigma$-moderate.}\medskip

Several proofs in this work are transposition to the parabolic framework of the constructions performed in \cite{MV-CONT} and \cite{M-JAM}. However, for the sake of completeness and due to the technicalities involved, we kept many of them, sometimes under an abriged form.


\section{The $\frak T_q$-fine topology}

We assume that $q\geq1+\frac{2}{N}$ and set $q'=\frac{q}{q-1}$. We recall that a set $E\subset\BBR^N$ is $(\frac{2}{q},q')$-{\it thin} at a point $a$ if
\begin{equation}\label{thin}
\myint{0}{1}\left(\myfrac{C_{\frac{2}{q},q'}(E\cap B_s(a))}{s^{N-\frac{2}{q-1}}}\right)^{q-1}\myfrac{ds}{s}<\infty.
\end{equation} 
If the value of the above integral is infinite, the set $E$ is called $(\frac{2}{q},q')$-{\it thick} at $a$. A set $U$ is  a $(\frac{2}{q},q')$-fine neighborhood of one of its point $a$ if $U^c$ is thin at $a$. It is $(\frac{2}{q},q')$-{\it finely open}, if $U^c$ is thin at any point $a\in U$. It is $(\frac{2}{q},q')$-{\it finely closed} if it complement is $(\frac{2}{q},q')$- finely open. For simplicity we will denote by $\frak T_q$ the $(\frac{2}{q},q')$-{\it fine topology} associated to these notions (see \cite[Chap 6]{AH} for a thorough discussion of these notions).
We say that a set $E\subset\Ren$ is $\frak T_q$-open (resp $\frak T_q$-closed) if it is open (resp. closed) in the $\frak T_q$-topology.
\begin{notation}
Let $A,\;B\subset\Ren.$\\
a) $A$ is $\frak T_q$-essentially contained in $B$, denoted $A\qsub B,$ if
$$\qcap(A\setminus B)=0.$$
b) The sets $A,\;B$ are $\frak T_q$-equivalent, denoted $A\qeq B,$ if
$$\qcap(A\xD B)=0.$$
c)The $\frak T_q$-closure of a set $A$ is denoted by $\widetilde{A}.$ The $\frak T_q$-interior of $A$ is denoted by $A^\diamond.$\\
d) Given $\xe>0,$ $A^\xe$ denotes the $\xe-$neighbourhood of $A$ for the standard Euclidean distance in $\BBR^N$\\
e) The set of $\frak T_q$-thick  points of $A$ is denoted by $b_q(A).$ The set of  $\frak T_q$-thin points of $A$ is denoted by $e_q(A).$
 \end{notation}
$$A \mathrm{\;is}\;\frak T_q\,\text{-open}\Leftrightarrow A\subset e_q(A^c),\qquad B \text{\;is}\;\frak T_q\,\text{closed}\Leftrightarrow b_q(B)\subset B.$$
Consequently,
$$\widetilde{A}=A\bigcup b_q(A),\qquad A^\diamond=A\cap e_q(A^c).$$
The capacity $\qcap$ possesses the Kellogg property (see  \cite[Cor. 6.3.17]{AH}), namely,
\be
\qcap(A\setminus b_q(A))=0.\label{kellog}
\ee
Therefore
$$A\qsub b_q(A)\qeq\tilde{A},$$
but, in general, $b_q(A)$ does not contain $A.$
\begin{prop}
(i) If $Q$ is a $\frak T_q$-open, then $e_q(Q^c)$ is the largest $\frak T_q$-open set that is $\frak T_q$-equivalent to $Q.$\\
(ii) If $F$ is a $\frak T_q$-closed then $b_q(F)$ is the smallest $\frak T_q$-closed set that is $\frak T_q$-equivalent to $F.$
\end{prop}
The proof is \cite[Prop. 2.1]{MV-CONT}. We collect below several facts concerning the $\frak T_q$-topology that are used throughout the paper.
\begin{prop}
Let $q\geq1+\frac{2}{N}.$\smallskip

\noindent i) Every $\frak T_q$-closed set is $\frak T_q$-quasi closed (\cite[Prop 6.4.13]{AH}).\smallskip

\noindent
ii)If E is $\frak T_q$-quasi closed then $E\qeq\widetilde{E}$ (\cite[Prop 6.4.12]{AH}).\smallskip

\noindent
iii)A set $E$ is $\frak T_q$-quasi closed if and only if there exists a sequence $\{E_m\}$ of closed subsets of $E$ such that $\qcap(E\setminus E_m)\rightarrow0$ (\cite[Prop. 6.4.9]{AH}).\smallskip

\noindent
iv) There exists a positive constant $c$ such that, for every set $E,$
$$\qcap(\widetilde{E})\leq c\qcap(E),$$
(\cite[Prop 6.4.11]{AH}).\smallskip

\noindent
v) If $E$ is $\frak T_q$-quasi closed and $F\qeq E$ then $F$ is $\frak T_q$-quasi closed.\smallskip

\noindent
vi) If $\{E_i\}$ is an increasing sequence of arbitrary Borel sets then
$$\qcap(\bigcup E_i)=\lim_{i\rightarrow\infty}\qcap(E_i).$$
vii) If $\{K_i\}$ is a decreasing sequence of compact sets then
$$\qcap(\bigcap K_i)=\lim_{i\rightarrow\infty}\qcap(K_i).$$
viii) Every Suslin set and, in particular, every Borel set $E$ satisfies
\bea
\nonumber
\qcap(E)&=&\inf\{\qcap(G):\;E\subset G,\;G\;\mathrm{open}\}\\ \nonumber
&=&\sup\{\qcap(K):\;K\subset E,\;K\;\mathrm{compact}\}.
\eea\label{propertiescapacity}
\end{prop}
For the last three statements see \cite[Sec. 2.3]{AH}. Statement (v) is an easy consequence of \cite[Prop. 6.4.9]{AH}. However note that this assertion is no longer valid if "$\frak T_q$-quasi closed" is replaced by "$\frak T_q$-closed." Only the following weaker statements holds:
$$\text{If }E\text{ is }\frak T_q\,\text{-closed and $A$ is a set such that }\qcap(A)=0\text{ then }E\cup A\text { is }\frak T_q\,\text{-closed}.$$

The next corollary is an easy consequence of (iii).
\begin{coro}
A set $E$ is $\frak T_q$-quasi closed if and only if there exists a sequence $\{E_m\}$ of $\frak T_q$-quasi closed subsets of $E$ such that $\qcap(E\setminus E_m)\rightarrow0.$\label{coro1}
\end{coro}

\begin{defin}
Let $E$ be a $\frak T_q$-quasi closed set. An increasing sequence $\{E_m\}$ of closed subsets of $E$ such that $\qcap(E\setminus E_m)\rightarrow0$ is called a $\frak T_q$-stratification of $E.$\\
(i) We say that $E_m$ is a proper $\frak T_q$-stratification of $E$ if
$$\qcap(E_{m+1}\setminus E_m)\leq\frac{1}{2^{m+1}}.$$
(ii) If $V$ is a $\frak T_q$-open set such that $\qcap(E\setminus V)=0$ we say that $V$ is a $\frak T_q$-quasi neighborhood of $E.$
\end{defin}

The following separation statement is valid in any locally compact metric space.

\begin{lemma}
Let $K$ be a closed subset of an open set $A.$ Then there exists an open set $G$ such that
$$K\subset G\subset\overline{G}\subset{A}.$$\label{anoixto}
\end{lemma}
\Proof    Let $x\in K.$ We set $B_n=B_n(x);\;n\in\mathbb{N}$ and $K_n=\overline{B}_n\cap K.$ Since $K_n$ is compact, we can easily show that there exists a decreasing sequence $\{\xe_n\}$ converging to $0$ such that $K^{\xe_n}_n\subset\overline{K^{\xe_n}_n}\subset A.$ Now we have
$$\bigcup_{n=1}^\infty K^{\frac{\xe_n}{2}}_n\subset \bigcup_{n=1}^\infty\overline{K^{\frac{\xe_n}{2}}_n}
\subset\bigcup_{n=1}^\infty K^{\xe_n}_n\subset A.$$
If we prove that the set $$\bigcup_{n=1}^\infty\overline{K^{\frac{\xe_n}{2}}_n}$$ is closed then the proof follows with $G=\bigcup_{n=1}^\infty K^{\frac{\xe_n}{2}}_n.$ We will prove it by contradiction. We assume that there exists a sequence $x_n\in\bigcup_{n=1}^\infty\overline{K^{\frac{\xe_n}{2}}_n}$ such that $x_n\rightarrow x$ and $x\notin \bigcup_{n=1}^\infty\overline{K^{\frac{\xe_n}{2}}_n}.$
We have $x_1=x_{n_1}$ such that $dist(x_{n_1},K)=\inf\{|x_{n_1}-y|:\;y\in K\}\leq\frac{\xe_1}{2}.$ Also we assert that there exists $x_{n_2}$ such that $dist(x_{n_2},K)\leq\frac{\xe_2}{2}.$ Indeed, If this is not valid then $\forall n\in\mathbb{N}$ we have $\frac{\xe_2}{2}<dist(x_n,K)\leq \frac{\xe_1}{2},$ which implies $x\in K_1.$ Thus we have clearly a contradiction. Inductively, we can construct a subsequence $\{x_{n_k}\}$ such that $dist(x_{n_k},K)\leq\frac{\xe_k}{2},\;\forall k\in\mathbb{N}.$ If we send $k$ to infinite, we reach to a contradiction, since we would have $dist(x,K)=0$ and using the fact that $K$ is closed, we would obtain that $x\in K.$\hfill$\Box$\\

In the framework of the $\frak T_q$-topology, the preceding result admits the following counterpart.
\begin{lemma}
Let $E$ be a $\frak T_q$-closed set. Then:\\
(i) Let $D$ be an open set such that $\qcap(E\setminus D^c)=0.$ Then there exists an open set $O$ such that
\be
E\qsub O\subset\widetilde{O}\qsub D.\label{1}
\ee
(ii) Let $D$ be a $\frak T_q$-open set such that $E\qsub D.$ Then there exists a $\frak T_q$-open set $O$ such that (\ref{1}) holds.\label{q-open}
\end{lemma}
\Proof    (i) Since $E\cap D\qeq E$, $E\cap D$ is $\frak T_q$-quasi closed, (see the discussion of the quasi topology in \cite[sec. 6.4]{AH}). Thus there exists a proper $\frak T_q$-stratification of $E\cap D,$ say $\{E_m\}$ and $E\qeq E'=\bigcup_{i=1}^\infty E_i$. If $E'$ is closed the result follows by Lemma \ref{anoixto}. We assume that $E'$ is not closed. Thus, we can assume without loss of generality that
$$E_{m+1}\setminus E_{m}\neq\emptyset \qquad\forall m\in\mathbb{N}.$$
We set $E_m'=G,$ where $G$ is the open set of Lemma \ref{anoixto} with $K=E_m$ and $A=D.$ Now since $\qcap(E_m\setminus E_{m-1})<\frac{1}{2^{m+1}},$ there exists an open set $D_m\supset E_{m}\setminus E_{m-1}\;;m\geq2,$ such that $\qcap(D_m)<\frac{1}{2^m}.$ Also we set $D_1=E'_1.$
Also we have by Lemma (\ref{anoixto}), $$D_m\cap E_m\subset\widetilde{D_m\cap E_m}\subset\widetilde{E_m}\subset D\qquad\forall m\in\mathbb{N}.$$
Also, since $E'=E_1\cup\bigcup_{m=2}^\infty(E_m\setminus E_{m-1})$ we have that
$$E'\subset\bigcup_{m=1}^\infty D_m\cap E_m'\subset\bigcup_{m=1}^\infty \widetilde{D_m\cap E_m'}\subset D.$$
Thus, it is enough to prove that the set $\bigcup_{m=1}^\infty \widetilde{D_m\cap E_m'}$ is $\frak T_q$-quasi closed. Indeed, for each $n>1,$ we have
\bea
\nonumber
\qcap\left(\bigcup_{m=1}^\infty \widetilde{D_m\cap E_m'}\setminus\bigcup_{m=1}^n \widetilde{D_m\cap E_m'}\right)&\leq&\qcap\left(\bigcup_{m=n+1}^\infty \widetilde{D_m\cap E_m'}\right)\leq\sum_{m=n+1}^\infty\qcap(\widetilde{D_m})\\ \nonumber
&\leq& c\sum_{m=n+1}^\infty\qcap(D_m)\leq c\sum_{m=n+1}^\infty2^{-m}.
\eea
And the result follows by Corollary \ref{coro1}, since $\bigcup_{m=1}^n \widetilde{D_m\cap E_m'}$ is $\frak T_q$-quasi closed.\\
The proof of (ii)  is same as in \cite[Lemma 2.4 (ii)]{MV-CONT}.\hfill$\Box$
\begin{lemma}(I) Let $E$ be a $\frak T_q$-closed set and $\{E_m\}$ a proper $\frak T_q$-stratification for $E$. Then there exists a decreasing sequence of open sets $\{Q_j\}$ such that $\cup E_m:=E'\subset Q_j$ for every $j\in\mathbb{N}$ and\\
(i)  $\cap_{j} Q_j= E',$    $\widetilde{Q}_{j+1}\qsub Q_j,$\\
(ii) $\lim_{j\to\infty}\qcap(Q_j)=\qcap(E).$\\
(II) If $A$ is a $\frak T_q$-open set, there exists a decreasing sequence of open sets $\{A_m\}$ such that
$$A\subset\bigcap_m A_m=:A',\qquad\qcap(A_m\setminus A')\rightarrow0 \text { as }m\to\infty,\quad A\qeq A'.$$
Furthermore there exists an increasing sequence of closed sets $\{F_j\}$ such that $F_j\subset A'$ and\\
(i) $\cup F_j=A',\qquad F_j\qsub F_{j+1}^\diamond$\\
(ii) $\qcap(F_j)\rightarrow\qcap(A')$ as $j\to\infty$.
\label{6}
\end{lemma}
\Proof
Let $\{D_j\}$ be a decreasing sequence of open sets such that $D_j\supset E,\;\forall j\in\mathbb{N}$ and
$$\lim_{j\to\infty}\qcap(D_j)=\qcap(E')=\qcap(E).$$
{\it Case 1}: $E$ is closed (thus $E_m=E$ for any $m\in\mathbb{N}$).\\
By Lemma \ref{anoixto} there exists a decreasing sequence $\{\xe_{1,n}\}$ converging to $0$, such that $\xe_{1,1}<1,$ and
$$E\subset Q_1=\bigcup_{n=1}^\infty K^{\frac{\xe_{1,n}}{2}}_n\subset\overline{Q}_1\subset D_1,$$
where $K_n=B_n(x)\cap E\;,x\in E.$ Also we have proven in Lemma \ref{anoixto} that the set $\bigcup_{n=1}^\infty\overline{K^{\frac{\xe_{1,n}}{2}}_n}$ is closed.\\
Again by Lemma \ref{anoixto} there exists a decreasing sequence $\{\xe_{2,n}\}$ converging to $0$, such that $\xe_{2,n}\leq\xe_{1,n}$ for all $ n$ and
$$E\subset Q_2=\bigcup_{n=1}^\infty K^{\frac{\xe_{2,n}}{4}}_n\subset\overline{Q}_2\subset D_2.$$
We note here that
$$\overline{Q}_2\subset\bigcup_{n=1}^\infty\overline{K^{\frac{\xe_{2,n}}{4}}_n}\subset\bigcup_{n=1}^\infty K^{\frac{\xe_{1,n}}{2}}_n,$$
and since $\overline{K^{\frac{\xe_{2,n}}{4}}_n}$ is closed we have
$$Q_2\subset\overline{Q}_2\subset Q_1.$$
By induction, we construct a decreasing sequence $\{\xe_{j,n}\}$ converging to $0$ with respect to $n$, such that $\forall n\in\mathbb{N}:$ $\xe_{j,n}\leq\xe_{k,n}$ for all $ j\geq k,$
$$E\subset Q_j=\bigcup_{n=1}^\infty K^{\frac{\xe_{j,n}}{2^{j+1}}}_n\subset\overline{Q}_j\subset D_j,$$
and
$$Q_j\subset\overline{Q}_j\subset Q_{j-1}.$$
Now note that
$$E\subset Q_j\subset E^{\frac{1}{2^{j}}},$$
thus $E=\cap Q_j.$
Finally,
$$\qcap(E)\leq\lim\qcap(Q_j)\leq\lim\qcap(D_j)=\qcap(E),$$
and the result follows in this case.\\
{\it Case 2}:  $E$ is not closed.\\
There exists a proper $\frak T_q$-stratification of $E,$ say $\{E_m\}$ and $E\qeq E'=\bigcup_{i=1}^\infty E_i$. Also by the Case 1, we can assume without loss of generality that
$$E_{m+1}\setminus E_{m}\neq\emptyset \qquad\forall m\in\mathbb{N}.$$
Let us denote by $Q_j^m$ the sets denoted by $Q_j$ in the previous case if we replace $E$ by $E_m.$ Since there holds $\qcap(\widetilde{E_m\setminus E_{m-1}})\leq c\qcap(E_m\setminus E_1),$ we can choose an open set $D_m^1$ such that $\qcap(D_m^1)\leq\frac{c}{2^m}.$ In view of Lemma (\ref{q-open}) the set
$$Q_1=\bigcup_{m=1}^\infty D_m^1\cap Q_1^m$$
is open and
$$E'\subset Q_1\subset\widetilde{Q}_1\subset D_1.$$
Furthermore the set $$\bigcup_{m=1}^\infty \widetilde{D_m^1\cap Q_1^m}$$
is $\frak T_q$-quasi closed.
By Lemma \ref{q-open} there exists an open set $D_m^2$ such that
$$D_m^2\subset\widetilde{D}^2_m\subset D_m^1.$$
By induction, we construct a sequence of open sets $\{D_m^j\}$ such that
$$D_m^j\subset\widetilde{D}^j_m\subset D_m^{j-1}\qquad\qcap(D_m^j)\leq\frac{c}{2^m}.$$
Thus in view of Lemma \ref{q-open} the set
$$Q_j=\bigcup_{m=1}^\infty D_m^j\cap Q_j^m$$ is open and
the set $$\bigcup_{m=1}^\infty \widetilde{D_m^j\cap Q_j^m}$$
is $\frak T_q$-quasi closed.
For any $m$ we have
$$D_m^j\cap Q_j^m\subset\widetilde{D_m^j\cap Q_j^m}\subset\widetilde{D_m^j}\cap \widetilde{Q_j^m}\subset D_m^{j-1}\cap Q_{j-1}^m.$$
Thus
$$Q_j\subset\widetilde{Q_j}\subset\bigcup_{m=1}^\infty \widetilde{D_m^j\cap Q_j^m}\subset\bigcup_{m=1}^\infty D_m^{j-1}\cap Q_{j-1}^m\subset D_j.$$
Since the set $\bigcup_{m=1}^\infty \widetilde{D_m^j\cap Q_j^m}$ is $\frak T_q$- quasi closed
we have $$Q_j\subset\widetilde{Q_j}\subset Q_{j-1}.$$
Finally
$$E'\subset Q_j\subset {E'}^{\frac{1}{2^{j}}},$$
thus $E'=\cap Q_j.$
The result follows in this case since
$$\qcap(E)\leq\lim\qcap(Q_j)\leq\lim\qcap(D_j)=\qcap(E).$$
(II) The proof is same as in \cite[Lemma 2.6 (b)]{MV-CONT} and we omit it.\hfill$\Box$
\\

The next results are respectively proved in \cite[Lemma 2.5]{MV-CONT} and  \cite[Lemma 2.7]{MV-CONT}.

\begin{prop}Let $E$ be a bounded $\frak T_q$-open set and let $\mathcal{D}$ be a cover of $E$ consisting of $\frak T_q$-open sets. Then, for every $\xe>0$ there exists an open set $O_\xe$ such that $\qcap(O_\xe)<\xe$ and $E\setminus O_\xe$ is covered by a finite subfamily of $\mathcal{D}.$\label{cover}
\end{prop}
\begin{prop}Let $Q$ be a $\frak T_q$-open set. Then, for every $\xi\in Q,$ there exists a $\frak T_q$-open set $O_\xi$ such that
\begin{equation}\label{perioxi}
\xi\in Q_\xi\subset\widetilde{Q}_\xi\subset Q.
\end{equation}
\end{prop}

\setcounter{equation}{0}
\section{Lattice structure of $\CU_+(Q)$}
Consider the equation
\be
\prt_t  u-\xD u+|u|^{q-1}u=0,\quad\mathrm{in}\;Q_\infty=\mathbb{R}^N\times(0,T],\text{ where } q\geq1+\frac{2}{N}.\label{maineq}
\ee
A function $u\in L^q_{loc}(Q_T)$ is a subsolution (resp. supersolution) of the equation if $\prt_t  u-\xD u+ |u|^{q-1}u\leq0$ (resp. $\geq0$) holds in the sense of distributions.

If $u\in L^q_{loc}(Q_T)$ is a subsolution of the equation then by Kato's inequality $(\prt_t-\xD)|u|+|u|^q\leq 0$  in the sense of distributions. Thus $|u|$ is a subsolution of the heat equation and consequently $u\in L^\infty_{loc}(Q_T).$ If $u\in L^q_{loc}(Q_T)$ is a solution then $u\in C^{2,1}(Q_T).$
\begin{prop} Let $u$ be a non-negative function in $L^\infty_{loc}(Q_T).$\\
(i) If $u$ is a subsolution of (\ref{maineq}), there exists a minimal solution $v$ dominating $u,$\\
i.e. $u\leq v\leq U$ for any solution $U\geq u.$\\
(ii) If $u$ is a supersolution of (\ref{maineq}), there exists a maximal solution $w$ dominated by $u,$\\
i.e. $V\leq w\leq u$ for any solution $V\leq u.$\\
All the above inequalities hold almost everywhere .\label{subsolution}
\end{prop}
\Proof   (i) Let $\{J_\xe\}$ be a filter of mollifiers in $\mathbb{R}^{N+1}.$ If $u$ is extended by zero outside of $Q_T,$ then the function $u_\xe=J_\xe*u$ belong to $C^\infty(\mathbb{R}^{N+1})$, $\lim_{\xe\rightarrow0}u_\xe=\widetilde{u}=u$ a.e. in $\mathbb{R}^{N+1}$ and $u_\xe\rightarrow u$ in $L^q_{loc}(Q_T).$ We note that we can choose $\xe>0$ small enough such that the function $u_\xe$ is a subsolution in
$B_R(0)\times(s,\infty)$ where $R>0$ and $0<s.$
Let $v_\xe$ be the positive solution of
\begin{equation}\label{1'}\BA{ll}
\prt_t  v-\xD v+ |v|^{q-1}v=0,\qquad&\mathrm{in}\;B_R(0)\times(s,\infty),\\
\phantom{\prt_t  v-\xD v+ |v|^{q-1}}v=u_\xe,\qquad&\mathrm{on}\;\partial B_R(0)\times(s,\infty),\\
\phantom{,v-\xD v+ |v|}
v(.,s)=u_\xe(.,s)\qquad&\mathrm{in}\;B_R(0).
\EA\end{equation}
In view of the proof of Lemma 2.4 and Remark 2.5 in \cite{MV-CPDE} we can prove that $v_\xe\geq u_\xe. $ Since $v_\xe$ is a subsolution of the heat equation, we have $v_\xe\leq||u_\xe||_{L^\infty(B_R(0)\times(s,T])}\leq||u||_{L^\infty(B_R(0)\times(s,T])}.$ Thus there exists a decreasing sequence $\xe_j$ converging to $0$
such that $v_{\xe_j}\rightarrow v$ in $L^q(B_R(0)\times(s,T]),$ $u\leq v\leq||u||_{L^\infty(B_R(0)\times(s,T])};0<s<T<\infty$ and $v$ is a positive solution of
\begin{equation}\label{2}\BA{ll}
\prt_t  v-\xD v+ |v|^{q-1}v=0,\qquad&\mathrm{in}\;B_R(0)\times(s,T],\\
\phantom{\prt_t  v-\xD v+ |v|^{q-1}}v=u,\qquad&\mathrm{on}\;\partial B_R(0)\times(s,T],\\
\phantom{,v-\xD v+ |v|}v(.,s)=u(.,s)\qquad&\mathrm{in}\;B_R(0).
\EA\end{equation}
Let $\{R_j\}$ be an increasing sequence tending to infinity and $s_j$ be a decreasing one converging to $0$. Let $v_j$ be the positive solution of the above problem with $R=R_j$ and $s=s_j.$ Since $v_j\geq u,$ we have by the maximum principle that $v_{j+1}\geq v_j.$ Thus, by Keller-Osserman inequality and standard parabolic regularity results, there exists a subsequence, say $\{v_j\},$ such that
$v_j\rightarrow v$ locally uniformly in $Q_T.$ The results follows in this case by the construction of $v.$\\
(ii)
Since $u\in L^q(B_R(0)\times(s,T])$ there exists a solution $w$ of the problem
\begin{equation}\label{3'}\BA{ll}
\prt_t  w-\xD w +|u|^{q}=0,\qquad&\mathrm{in}\;B_R(0)\times(s,T]\\
\phantom{\prt_t  -\xD w +|u|^{q}}
w=0,\qquad&\mathrm{on}\;\partial B_R(0)\times(s,T]\\
\phantom{-\xD w/ |u|^{q}}
w(.,s)=0\qquad&\mathrm{in}\;B_R(0).
\EA\end{equation}
Hence $u+w$ is supersolution of the heat equation with boundary and initial data $u.$ Consequently, $u+w\geq z$ where $z$ is the solution of the heat equation with boundary and initial data $u.$ Also, the function $z-w$ is a subsolution, thus there exists a solution $v\leq u$ of the problem (\ref{2}) with boundary and initial data $u.$ As before, let $\{R_j\}$ be an increasing sequence tending to infinity and $s_j$ be a decreasing sequence tending to $0$. Let $v_j$ be the positive solution of the problem (\ref{2}) with $R=R_j$ and $s=s_j.$ Since $v_j\leq u,$ we have by maximum principle that $v_{j+1}\leq v_j.$ Thus by standard parabolic arguments, there exists a subsequence, say $\{v_j\},$ such that
$v_j\rightarrow v$ locally uniformly in $Q_\infty.$ Again, the construction of $v$ implies the result.\hfill$\Box$
\begin{prop}
Let $u$ and $v$ be nonnegative, locally bounded functions in $Q_T.$\\
(i) If $u$ and $v$ are subsolutions (resp. supersolutions) then $\max(u,v)$ is a subsolution (resp. $\min(u,v)$ is a supersolution).\\
(ii) If $u$ and $v$ are supersolutions then $u+v$ is a supersolution.\\
(iii) If $u$ is a subsolution and $v$ is a supersolution then $(u-v)_+$ is a subsolution.\label{maxsub}
\end{prop}
\Proof    The first two statements are immediate  consequence of the parabolic Kato's inequality. The third statement is verified in a  similar way since
$$(\frac{d}{dt}-\xD)(u-v)_+\leq sign_+(u-v)(\frac{d}{dt}-\xD)(u-v)\leq-sign_+(u-v)(u^q-v^q)\leq-(u-v)^q_+.$$
\hfill$\Box$
\begin{notation}
Let $u,\;v$ be nonnegative, locally bounded functions in $Q_T.$\\
(a) If $u$ is a subsolution, $[u]_\dag$ denotes the smallest solution dominating $u.$\\
(b) If $u$ is a supersolution, $[u]^\dag$ denotes the largest solution dominated by $u.$\\
(c) If $u,\;v$ are subsolutions then $u\vee v:=[\max(u,v)]_\dag.$\\
(d) If $u,\;v$ are supersolutions then $u\wedge v:=[\inf(u,v)]^\dag$ and $u\oplus v=[u+v]^\dag.$\\
(e) If $u$ is a subsolution and $v$ is a supersolution then $u\ominus v:=[(u-v)_+]_\dag.$\label{34}
\end{notation}
\begin{prop}
(i) Let $\{u_k\}$ be a sequence of positive, continuous subsolutions of (\ref{maineq}). Then $U:=\sup u_k$ is a subsolution. The statement remains valid if subsolution is replaced by supersolution and $\sup$ by $\inf.$\\
(ii) (\cite{Dyn}) Let $\mathcal{T}$ be a family of positive solutions of (\ref{maineq}). Suppose that, for every $u_1$ and $u_2$ belonging to $\mathcal{T}$ there exists $v\in\mathcal{T}$ such that
$$\max(u_1,u_2)\leq v,\qquad \mathrm{resp.}\;\min(u_1,u_2)\geq v.$$
Then there exists a monotone sequence $\{u_n\}$ in $\mathcal{T}$ such that
$$u_n\uparrow\sup\mathcal{T},\qquad\mathrm{resp.}\;u_n\downarrow\inf T.$$
Thus $\sup \mathcal{T}$ (resp. $\inf\mathcal{T}$) is a solution.\label{sygklish}
\end{prop}
\Proof   (i) Set $v_j=\max\left(u_1,u_2,...,u_j\right)=\max\left(\max(u_1,u_2),\max(\max(u_1,u_2),u_3),...,\max(\max(...),u_j)\right).$ By proposition \ref{maxsub} $v_j$ is a subsolution and $v_{j+1}\geq v_j.$ Thus the positive solution $[v_j]_\dag$ is increasing with respect to $j.$ Also by Keller-Osserman inequality, we have that $[v_j]_\dag\rightarrow \widetilde{v},$ where $\widetilde{v}$ is a positive solution. Thus $v_j\rightarrow v$ where $v$ is a subsolution of (\ref{maineq}). Now since $u_i\leq v$ for each $i\in\mathbb{N},$ we have that $U\leq v.$ But $v_j\leq U$ for each $j\in\mathbb{N},$ which implies $v\leq U.$ And thus $v=U.$ The proof for "$\inf$" is similar and we omit it.\\
(ii) The proof is similar as the one in \cite{Dyn}. Let $A = { (x_n,t_n)}$ be a countable dense subset of $Q_T$ and let $u_{nm}\in \mathcal{T}$ satisfy the condition
$\sup_m u_m(x_n,t_n)= w(x_n,t_n).$ Since $\mathcal{T}$ is closed with respect to $\vee$, there exists an increasing sequence
of $v_n\in\mathcal{T}$ such that $v = \lim_{n\to\infty} v_n$, coincides with $w$ on $A.$ We claim that $v=w$ everywhere. Indeed,
$v\leq u$. Suppose $u\in \mathcal{T}$. Then $u\leq w$ and therefore $u \leq v$ on $A$. Since $A$ is everywhere dense and $u,\;v$
are continuous, $u\leq v$ everywhere in $Q_\infty,$ which implies $u\geq w = \sup u.$\hfill$\Box$\\

As a consequence we have the following result which extends to equation (\ref{E1}) what Dynkin proved for (\ref{El1}) \cite[Theorem 5.1]{Dyn}.
\begin{theorem}\label{lat} The set $\CU_+(Q_T)$ is a complete lattice stable for the laws $\oplus$ and $\ominus$.
\end{theorem}
\setcounter{equation}{0}
\section{Partition of unity in Besov spaces}
\begin{lemma}
Let $U\subset\mathbb{R}^N$ be a $\frak T_q$-open set and $z\in U.$ Then there exists a function $f\in W^{\frac{2}{q},q'}(\mathbb{R}^N)$ with compact support in $U$ such that $f(z)>0.$ In particular, there exists a bounded $\frak T_q$-open set $V$ such that $\overline{V}\subset U.$
\end{lemma}
\Proof    We suppose that $z$ is not an interior point of $U$ with respect to Euclidean topology, since otherwise the result is obvious. Since $U$ is $\frak T_q$-open we have that $U^c$ is thin at $z.$ Also by the assumption on $z,$ we have that $z\in \overline{U^c}\setminus U.$ By \cite[p. 174]{AH}, we can find an open set $W\supset U^c$, $z\in \overline{W}\setminus W$ and $W$ is thin at $z$. \\
We recall that for a set $E$ with positive $C_{\frac{2}{q},q'}$-capacity, $F^{E}:=\mathcal{V}^{\xm_E}=G_{\frac{1}{q}}\ast(G_{\frac{1}{q}}\ast\gm_E)^{p-1}$ where $\xm_E$ is the capacitary measure on $E$. Then, by \cite[ Proposition 6.3.14]{AH}, there exists $r>0$ small enough such that
$$\mathcal{V}^\xm(z)<\frac{1}{2},$$
where $\xm$ is the capacitary measure of $B(z,r)\cap W$ and $\mathcal{V}^\xm$ the corresponding Besov potential (see \cite[Theorems 2.2.7, 2.5.6 ]{AH}). By \cite[Theorem 6.3.9]{AH},  $\mathcal{V}^\xm\geq1$ quasi everywhere (abr. q.a.e.) on $B(z,r)\cap W,$ and by \cite[Proposition 2.6.7]{AH} $\mathcal{V}^\xm\geq1$ everywhere on $B(z,r)\cap W.$ Thus
$$\mathcal{V}^\xm(z)<\frac{1}{2}<1\leq\mathcal{V}^\xm(x),\qquad\forall x\in B(z,r)\cap W.$$
Thus we can find $r_0>0$ small enough such that
$$\mathcal{V}^\xm(z)<\frac{1}{2}<1\leq\inf\{\mathcal{V}^\xm(x):\;x\in B(z,r_0)\setminus U\}.$$
Now let $0\leq H(t)$ be a smooth nondecreasing function such that $H(t)=t$ for $t\geq\frac{1}{4}$ and $H(t)=0$ for $t\leq0.$
Also let $\eta\in C_0^\infty(\mathbb{R}^N)$ such that $0\leq\eta\leq1,$ $\mathrm{supp}\;\eta\subset B(z,r_0)$ and $\eta(z)=1.$ Then the function
$$f(z)=\eta H(1-\mathcal{V}^\xm),$$
belongs to $W^{\frac{2}{q},{q'}}(\mathbb{R}^N).$ Since by definition $\mathcal{V}^\xm$ is lower semicontinuous,  the set $\{1-u\geq0\}$
is closed. Hence the support of $f$ is compact and
$$\mathrm{supp}f\subset\mathrm{supp}\eta\cap\{1-u\geq0\}\subset U.$$
\hfill$\Box$
\begin{lemma}
Let $U$ be a $\frak T_q$-open set and $z\in U.$ Then there exists a $\frak T_q$-open set $V,$ such that $z\in V\subset U,$ and a function $\psi\in W^{\frac{2}{q},q'}(\mathbb{R}^N)$ such that $\psi=1$ q.a.e. on $V$ and $\psi=0$ outside $U.$\label{cutoff}
\end{lemma}
\Proof    As before, we assume that $z$ is not an interior point of $U.$ Let $\mathcal{V}^\xm$ be the Besov potential of the previous lemma, with
$$\mathcal{V}^\xm(z)<\frac{1}{4},\qquad\mathcal{V}^\xm=1\qquad\mathrm{on}\;B(z,r_0)\setminus U.$$
By \cite[Proposition 6.3.10]{AH} $\mathcal{V}^\xm$ is quasi continuous, that we can find a $\frak T_q$-open set $W$ which contains $z$ such that
$$\mathcal{V}^\xm(x)\leq\frac{1}{4},\;\;\text{ q.a.e. on } W.$$
Let $\eta\in C_0^\infty(\mathbb{R}^N)$ such that $0\leq\eta\leq1,$ $\mathrm{supp}\;\eta\subset B(z,r_0)$ and $\eta(x)=1,\forall \;x\in B(z,\frac{r_0}{2}).$
Set
$$f=2\eta H\left(1-H\left(\frac{1}{2}-\mathcal{V}^\xm(x)\right)-\mathcal{V}^\xm(x)\right).$$
Then $f\in W^{\frac{2}{q},q'}(\mathbb{R}^N),$ $0\leq f\leq1$ and $f=0$ on $B(z,r_0)\setminus U.$ Also, $f=1$ on $B(z,\frac{r_0}{2})\cap W$ and
$f=0$ outside of $B(z,r_0)\cap U.$ \hfill$\Box$
\begin{lemma}
Let $\frac{2}{q}\leq1,$ $K$ be a compact set and $U$ be a $\frak T_q$-open set such that $K\subset U.$ Also, let $\{U_j\}$ be a sequence of $\frak T_q$-open subsets of $U$ covering $U$ up to a set of zero $\qcap$-capacity $Z.$ We assume that there exists a nonnegative $u\in W^{\frac{2}{q},q'}(\mathbb{R}^N)\cap L^\infty(\mathbb{R}^N)$ with $\frak T_q$-{supp}\;$ u\subset K\subset U.$ Then there exist $m(k)\in\mathbb{N}$ and nonnegative functions $u_{k,j}\in L^\infty(\mathbb{R}^N)$ with $\frak T_q$-{supp}\;$ u_{k,j}\subset U_j,$ such that
\be
\sum_{j=1}^{m(k)}u_{k,j}\leq u\label{xwris}
\ee
and
$$\lim_{k\rightarrow\infty}||u-\sum_{j=1}^{m(k)}u_{k,j}||_{W^{\frac{2}{q},q'}(\mathbb{R}^N)}=0.$$\label{partition}
\end{lemma}

\noindent\Remark  If $u$ changes sign, the conclusion of Lemma remains valid without inequality (\ref{xwris}).\\

\noindent\Proof    Without loss of generality we can assume that $U$ and the $\cup_jU_j$ are bounded. For any $j\geq0,$ there exists open sets $G_{k,j}$ such that $\qcap(G_{k,j})\leq 2^{-k-j},$ $Z\subset G_{k,0}$ and for $j\geq1,$ the sets $U_j\cup G_{k,j}$ are open. Also the sets
$$G_k=\bigcup_{j=0}^\infty G_{k,j},\qquad \bigcup_{j=1}^\infty G_k\bigcup U_j$$
are open and $\qcap(G_k)\rightarrow0$ when $k\to\infty$.

Since $ G_k$ is open, its Besov potential $F^{ G_k}$ is larger or equal to $1$ everywhere on $ G_k$ \cite[Theorems 2.5.6, 2.6.7 ]{AH}). Also we have
$$||\mathcal{V}^{\xm_k}||_{W^{\frac{2}{q},q'}(\mathbb{R}^N)}^{q'}\leq A\qcap(G_k),$$
where $A$ is a positive constant which depends only on $n,\;q.$ Now consider a smooth nondecreasing function $H$ such that $H(t)=1$ for $t\geq1$ and $H(t)=t$ for $t\leq\frac{1}{2},$ then the function $\xf_k=H(\mathcal{V}^{\xm_k})$ belongs to $ W^{\frac{2}{q},q'}(\mathbb{R}^N)$, satisfies $0\leq\xf_k\leq1,$ $\xf_k=1$ on $G_k$ and there exists a constant $A'(n,q)>0$ such that
$$
||\xf_k||_{W^{\frac{2}{q},q'}(\mathbb{R}^N)}^{q'}\leq A'\qcap(G_k).
$$
Set $\psi_k=1-\xf_k.$ By Lebesgue's dominated theorem
\begin{equation}\label{PU}
||u-\psi_ku||_{W^{\frac{2}{q},q'}(\mathbb{R}^N)}^{q'}\rightarrow0.
\end{equation}
Thus it is enough to prove that
\be
u\psi_k=\sum_{j=1}^{m(k)}u_{k,j}.\label{1''}
\ee
Fix $k\in\mathbb{N}.$ Then there exist open balls $B_{k,j,i}$, for $i,j=1,2...$, such that
$$\overline{B}_{k,j,i}\subset U_j\bigcup G_k,\qquad\mathrm{and}\;\bigcup_{j=1}^\infty G_k\bigcup U_j=\bigcup_{i,j=1}^\infty B_{k,j,i}.$$
Since $K$ is compact, there exists $m(k)\in\mathbb{N}$ such that
$$K\subset\bigcup_{i,j=1}^{m(k)} B_{k,j,i}.$$
Now consider $w_{k,j,i}\in C_0^\infty(\mathbb{R}^N)$ such that
$$\{w_{k,j,i}>0\}=B_{k,j,i}.$$
If we set
$$u_{k,j}=u\psi_k\frac{\sum_{i=1}^{m(k)}w_{k,j,i}}{\sum_{i,j=1}^{m(k)}w_{k,j,i}},$$
then  $u_{k,j}\in L^\infty(\mathbb{R}^N),$ satisfies $1$ and
$$\frak T_q\text{-}\mathrm{supp}u_{k,j}\subset (K\setminus G_k)\cap B_{k,j,i}\subset U_j.$$
\hfill$\Box$

\noindent\Remark We conjecture that the result still holds if $\frac{2}{q}>1$, but we have not been able to prove (\ref{PU}).

\setcounter{equation}{0}

\section{The regular set and its properties}
Let $q>1,$ $T>0.$ If $Q_T=\mathbb{R}^N\times(0,T),$ we recall that $\mathcal{U}_+(Q_T)$ is the set of positive solutions $u$ of
\be
\prt_t  u-\xD u+u^q=0\qquad\mathrm{in} \;\;Q_T.\label{2.14}
\ee
If a function $\xz$ is defined in $\mathbb{R}^N.$ We denote by $\frak T_q$-$\mathrm{supp}(\xz)$ the $\frak T_q$-closure of the set where $\abs\xz>0.$
Let $U$ be a Borel subset of $\mathbb{R}^N$ and $\chi_U$ be the characteristic function of $U.$ We set
$$\mathbb{H}(\chi_U)(x,t)=\frac{1}{(4\pi t)^\frac{N}{2}}\int_{\mathbb{R}^N}e^{-\frac{|x-y|^2}{4t}}\chi_Udy.$$
For any $\xi\in\mathbb{R}^N$ the following dichotomy occurs:\\\\
(i) either there exists a $\frak T_q$-open bounded neighborhood $U=U_\xi$ of $\xi$ such that
\be
\int_0^T\int_{\mathbb{R}^N}u^q\mathbb{H}[\chi_U]^{2q'}dxdt<\infty,\label{2.15}
\ee
where $q'=\frac{q}{q-1},$\\
(ii) or for any $\frak T_q$-open neighborhood $U$ of $\xi$
\be
\int_0^T\int_{\mathbb{R}^N}u^q\mathbb{H}[\chi_U]^{2q'}dxdt=\infty.\label{2.16}
\ee
\begin{defin}
 The set of $\xi\in\mathbb{R}^N$ such that (i) occurs is $\frak T_q$-open. It is denoted by $\mathcal{R}_q(u)$ and called the regular set of $u.$ Its complement $\mathcal{S}_q(u)=\mathbb{R}^N\setminus\mathcal{R}_q(u)$ is $\frak T_q$-closed and called the singular set of $u.$\label{def2.1}
\end{defin}
\begin{prop}\label{menest}
Let $\eta\in \besl$ with $\frak T_q$-support in a $\frak T_q$-open bounded set $U.$ Also let $u\in\mathcal{U}_+(Q_T)$ satisfy
$$
M_U=\int_0^T\int_{\mathbb{R}^N}u^q\mathbb{H}[\chi_U]^{2q'}dxdt<\infty.
$$
Then there exists
\be
l(\eta):=\lim_{t\rightarrow0}\int_{\mathbb{R}^N}u(x,t)\mathbb{H}[\eta]^{2q'}_+dx.\label{2.19}
\ee
Furthermore
\be
|l(\eta)|\leq C(M_U,q)\left(||\eta||^{2q'}_{W^{\frac{2}{q},q'}}+||\eta||^{2q'}_{L^\infty(\mathbb{R}^N)}\right).\label{2.20}
\ee\label{lem2.3}
\end{prop}
\Proof    Put $h=\mathbb{H}[\eta]$ and $\xf(r)=r^{2q'}_+.$ Since $|\eta|\leq||\eta||_{L^\infty}\chi_U,$ there holds
\be
\left|\int_0^T\int_{\mathbb{R}^N}u^q\xf(h)dxdt\right|\leq||\eta||_{L^\infty}^{2q'}\int_0^T\int_{\mathbb{R}^N}u^q\mathbb{H}[\chi_U]^{2q'}dxdt
:=||\eta||_{L^\infty}^{2q'}M_U<\infty.\label{2.21}
\ee
Moreover
\be
\int_s^t\int_{\mathbb{R}^N}(-u(\partial_t\xf(h)+\xD\xf(h)))+u^q\xf(h)dxd\tau=\int_{\mathbb{R}^N}u\xf(h)(.,s)dx-\int_{\mathbb{R}^N}u\xf(h)(.,t)dx.\label{2.22}
\ee
But
$$\partial_t\xf(h)+\xD\xf(h)=2q'\xf(h)h^{-2}_+(2h_+\prt_t  h+(2q'-1)|\nabla h|^2).$$
By H\"{o}lder
\bea
\nonumber
\bigg|\int_s^t\int_{\mathbb{R}^N}u(\partial_t\xf(h)+\xD\xf(h))dxd\tau\bigg|
\phantom{-------------------------}
\\ \nonumber
\leq\left(\int_s^t\int_{\mathbb{R}^N}u^q\xf(h)dxd\tau\right)^\frac{1}{q}
\left(\int_s^t\int_{\mathbb{R}^N}\xf(h)^{-\frac{q'}{q}}|(\partial_t\xf(h)+\xD\xf(h))|^{q'}dxd\tau\right)^\frac{1}{q'}\\ \nonumber
\leq 4q'\left(\int_s^t\int_{\mathbb{R}^N}u^q\xf(h)dxd\tau\right)^\frac{1}{q}
\left(\int_s^t\int_{\mathbb{R}^N}(h_+|\prt_t h|+|\nabla h|^2)^{q'}dxd\tau\right)^\frac{1}{q'}.\phantom{,,--}
\eea
By standard regularity properties of the heat kernel
$$\int_s^t\int_{\mathbb{R}^N}|\prt_t h|^{q'}dxd\tau\leq\int_0^T\int_{\mathbb{R}^N}|\prt_t h|^{q'}dxd\tau\leq||\eta||^{q'}_{W^{\frac{2}{q},q'}},$$
and by Gagliardo-Nirenberg inequality and the maximum principle
$$\int_s^t\int_{\mathbb{R}^N}|\nabla h|^{2q'}dxd\tau\leq\int_0^T\int_{\mathbb{R}^N}|\nabla h|^{2q'}dxd\tau\leq C||\eta||^{q'}_{L^\infty}||\xD h||^{q'}_{L^{q'}}=C||\eta||^{q'}_{L^\infty}||\prt_t h||^{q'}_{L^{q'}}.$$
Therefore,
\be
\bigg|\int_s^t\int_{\mathbb{R}^N}u(\partial_t\xf(h)+\xD\xf(h))dxd\tau\bigg|\leq C\left(\int_s^t\int_{\mathbb{R}^N}u^q\xf(h)dxd\tau\right)^\frac{1}{q}
||\eta||_{L^\infty}||\eta||_{W^{\frac{2}{q},q'}}.\label{2.23}
\ee
This implies that the left-hand side of (\ref{2.22}) tends to 0 when $s,t\rightarrow0$, thus there exists
$$l(\eta):=\lim_{s\rightarrow0}\int_{\mathbb{R}^N}u\xf(h)(x,s)dx.$$
From (\ref{2.22}) it follows 
\be
\int_0^T\int_{\mathbb{R}^N}(-u(\partial_t\xf(h)+\xD\xf(h)))+u^q\xf(h)dxd\tau+\int_{\mathbb{R}^N}u\xf(h)(.,T)dx=l(\eta).\label{2.24}
\ee
Since $|u\xf(h)(.,T)|\leq C(T)||\eta||^{2q'}_{L^\infty},$ we derive
\be
|l(\eta)|\leq C_1||\eta||^{2q'}_{L^\infty}+C||\eta||^{q'}_{L^\infty}||\eta||^{q'}_{W^{\frac{2}{q},q'}}\leq C\left(||\eta||_{L^\infty}+||\eta||_{W^{\frac{2}{q},q'}}\right)^{2q'}.\label{2.25}
\ee
\begin{prop}
Let the assumptions of Lemma \ref{lem2.3} be satisfied. Then
\be
\lim_{t\rightarrow0}\int_Uu(x,t)\eta^{2q'}_+(x)dx=l(\eta).\label{2.26}
\ee
\label{lem2.4}
\end{prop}
\Proof    Using (\ref{2.21}) with $h$ replaced by $h_s(x,t):=\mathbb{H}[\eta](x,t-s),$ we get
\be
\int_s^T\int_{\mathbb{R}^N}(-u(\partial_t\xf(h_s)+\xD\xf(h_s)))+u^q\xf(h_s)dxd\tau+\int_{\mathbb{R}^N}u\xf(h_s)(.,T)dx=\int_{\mathbb{R}^N}u\xf(h_s)(.,s)dx.
\label{2.27}
\ee
When $s\rightarrow0$
$$\int_{\mathbb{R}^N}u\xf(h_s)(.,T)dx\rightarrow\int_{\mathbb{R}^N}u\xf(h)(.,T)dx,$$
and
$$\int_s^T\int_{\mathbb{R}^N}u^q\xf(h_s)dxd\tau\rightarrow\int_0^T\int_{\mathbb{R}^N}u^q\xf(h)dxd\tau,$$
by the dominated convergence theorem. Furthermore,
\bea
\nonumber
\bigg| \int_0^{T-s}\int_{\mathbb{R}^N}(u(x,t+s)-u(x,t))(\partial_t\xf(h)+\xD\xf(h))dxdt\bigg|\phantom{-------------}\\ 
\nonumber
\leq C\left(\int_0^{T-s}\int_{\mathbb{R}^N}|u(x,t+s)-u(x,t)|^qh_+^{2q'}dxdt\right)^\frac{1}{q}||\eta||^{q'}_{L^\infty}||\eta||^{q'}_{W^{\frac{2}{q},q'}},
\eea
which tends to zero with $s$. Finally,
$$\lim_{s\to 0}\int_{T-s}^T\int_{\mathbb{R}^N}u^q\xf(h)dxd\tau=0.$$
Subtracting (\ref{2.22}) to (\ref{2.27}), we derive
$$\lim_{s\rightarrow0}\int_{\mathbb{R}^N}u(.,s)(\xf(h)(.,s)-\xf(\eta))dx=0,$$
which implies the claim.\hfill$\Box$\\

The next statement obtained by contradiction  with the use of Lemma \ref{lem2.3} and Lemma \ref{lem2.4} will be very useful in the sequel.
\begin{prop}
Assume that $U$ is a bounded $\frak T_q$-open set and
\be
\lim_{t\rightarrow0}\int_Uu(x,t)\eta^{2q'}(x)dx=\infty,\label{2.28}
\ee
for some $0\leq\eta\in\besl$ with $\frak T_q$-support in $U,$ then
\be
\int_0^T\int_{\mathbb{R}^N}u^q\mathbb{H}[\eta]^{2q'}dxdt=\infty\label{2.29}.
\ee\label{lem2.5}
\end{prop}

\begin{prop}
Let $\xi\in\mathcal{S}_q(u).$ Then for any $\frak T_q$-open set $G$ which contains $\xi,$ there holds
\be
\lim_{t\rightarrow0}\int_Gu(x,t)dx=\infty.\label{2.30}
\ee
\end{prop}
\Proof    If $\xi\in\mathcal{S}_q(u)$ and if $G$ is $\frak T_q$-open and contains $\xi,$ then by Lemma \ref{cutoff} there exist $\eta\in\besl$ and a  $\frak T_q$-open set $D\subset G$ such that $\eta=1$ on $D,$  $\eta=0$ outside of $G$ and $0\leq\eta\leq1.$ Thus
$$\infty=\int_0^T\int_{\mathbb{R}^N}u^q\mathbb{H}[\chi_D]^{2q'}dxdt\leq\int_0^T\int_{\mathbb{R}^N}u^q\mathbb{H}.[\eta]^{2q'}dxdt,$$
Therefore
$$
\lim_{t\rightarrow0}\int_{\mathbb{R}^N}u\mathbb{H}[\eta]^{2q'}dx=\infty,
$$
which implies
$$
\lim_{t\rightarrow0}\int_{\mathbb{R}^N}u\eta^{2q'}dx=\infty,
$$
and the result follows by the properties of $\eta.$\hfill$\Box$

\subsection {Moderate solutions}\label{moderatesec}
We first recall some classical results concerning initial value problem with initial measure data. A solution $u$ of (\ref{maineq}) is called {\it moderate} if $u\in L^q(K)$ for any compact $K\subset\overline{Q}_\infty$. Then there exists a unique Radon measure $\mu$ such that
\be
\lim_{t\rightarrow0}\int_{\mathbb{R}^N}u(x,t)\xz(x)dx=\int_{\mathbb{R}^N}\xz(x)d\xm\qquad\forall\xz\in C_0^\infty(\mathbb{R}^N).\label{5.18}
\ee
Equivalently
$$-\int\int_{Q_\infty} u(\xf_t+\xD \xf)dxdt+\int\int_{Q_\infty}|u|^{q-1}u\xf dxdt=\int_{\mathbb{R}^N}\xf(x,0)d\xm,$$
for all $\xf\in C^{1,1;1}(\overline{Q}_\infty),$ with compact support.\smallskip

The above measure has the property that it vanishes on Borel sets with $\qcap$-capacity zero. There exists an sequence $\{\xm_n\}\subset W^{-\frac{2}{q},q}(\mathbb{R}^N)$ of Radon measures such that $\xm_n\rightharpoonup\xm$ in the weak* topology. If we assume that $u$ is a positive moderate solution, or equivalently that the initial measure $\xm$ is positive, then the previous sequence can be constructed as being increasing and particularly $\{\xm_n\}\subset W^{-\frac{2}{q},q}(\mathbb{R}^N)\cap\mathfrak{M}^b_+(\mathbb{R}^N)$, where $\mathfrak{M}^b_+(\mathbb{R}^N)$ is the set of all positive bounded Radon measures in $\mathbb{R}^N.$

If $\xn\in W^{-\frac{2}{q},q}(\mathbb{R}^N)\cap\mathfrak{M}^b_+(\mathbb{R}^N),$ then we have for some constant $C>0$ independent on $\xn$(see Lemma 3.2-\cite{MV-CVPDE})
\be
C^{-1}||\xn||_{ W^{-\frac{2}{q},q}(\mathbb{R}^N)}\leq||\BBH[\xn]||_{L^q(Q_T)}\leq C||\xn||_{ W^{-\frac{2}{q},q}(\mathbb{R}^N)},\label{21'}
\ee
where we recall that $\BBH[\xn]$ denotes the heat potential of $\xn$ in $Q$.
\begin{lemma}
Let $u$ be a moderate positive solution with initial data $\xm.$ Then for any $T>0$ and bounded $\frak T_q$-open set we have
$$\int_0^T\int_{\mathbb{R}^N}u^q(t,x)\mathbb{H}^{2q'}[\chi_{O}]dxdt<\infty.$$
\end{lemma}
\Proof
Let $0\leq\eta\in C^\infty_0(\mathbb{R}^N)$ and $\eta=1$ on $O$ and $s<T.$ We define here $h=\mathbb{H}[\eta](x,t),$ $h_s=\mathbb{H}[\eta](x,t-s)$ and $\xf(r)=|r|^{2q'}.$ Then we have
$$
\int_s^T\int_{\mathbb{R}^N}u(x,t)\left(\partial_t\xf(h_s)+\xD\xf(h_s)\right)+|u|^q\xf(h_s)dxdt+\int_{\mathbb{R}^N}u\xf(h_s)(.,T)dx=\int_{\mathbb{R}^N}u(x,s)\xf(\eta)dx.$$
In view of Proposition \ref{lem2.3}, (\ref{2.23}) and H\"{o}lder's inequality, there exists a constant $c=c(q,N)$ such that
$$
\int_s^T\int_{\mathbb{R}^N}|u|^q\xf(h_s)dxdt+\int_{\mathbb{R}^N}u\xf(h_s)(.,T)dx\leq c\left(\int_{\mathbb{R}^N}u(x,s)\xf(\eta)dx+||\eta||_{L^\infty}^{2q'}||\eta||_{W^{\frac{2}{q},q'}}^{2q'}\right).$$
Using Fatou's lemma and the fact that, for any bounded Borel set $E$
$$\limsup_{s\rightarrow0}\int_{E}u(x,s)dx<\infty,$$
we conclude the proof.\hfill$\Box$
\begin{theorem}
Let u be a positive moderate solution with $\xm$ as initial data, then\\
(i) $\xm$ is regular relative to the  $\frak T_q$-topology.\\
(ii) For each quasi continuous function $\xf\in L^\infty(\mathbb{R}^N)$ with bounded $\frak T_q$-support in $\mathbb{R}^N,$ we have
$$
\lim_{t\rightarrow0}\int_{\mathbb{R}^N}u(x,t)\xf(x)dx=\int_{\mathbb{R}^N}\xf(x)d\xm.
$$\label{moderate}
\end{theorem}
\Proof    The proof is similar to the one given \cite{MV-CONT}.\\
(i) Every Radon measure on $\mathbb{R}^N$ is regular in the usual Euclidean topology, i.e.
\bea
\nonumber
\xm(E)&=&\inf\{\xm(D):\;E\subset D,\;D\;\mathrm{open}\}\\ \nonumber
      &=& \inf\{\xm(K):\;K\subset E,\;K\;\mathrm{compact}\},
\eea
for any Borel set $E.$ But if $D$ is open and contains $E$, it is $\frak T_q$-open, hence
$$\xm(E)\leq\inf\{\xm(D):\;E\subset D,\;D\,\;\frak T_q\text{-open}\}\leq\inf\{\xm(D):\;E\subset D,\;D\;\mathrm{open}\}=\xm(E),$$
and the result follows.\\
(ii) Since the measure $\xm_t=u(t,x)dx\rightharpoonup\xm$ in the weak* topology we have
$$\limsup_{t\rightarrow0}\xm_t(E)\leq\xm(E),\qquad\liminf_{t\rightarrow0}\xm_t(A)\geq\xm(A),$$
for any compact set $E,$ respectively, open set $A.$
This extends to any bounded $\frak T_q$-closed set $E$ (resp. $\frak T_q$-open set A).\\
Indeed, let $E$ be a $\frak T_q$-closed set and $\{K_n\}$ be an increasing sequence of closed sets such that $\qcap(E\setminus K_n)\rightarrow0.$
Then for any $m\in\mathbb{N}$ and any open set $E\subset O$  we have
$$\limsup_{t\rightarrow0}\xm_t(E)\leq\limsup_{t\rightarrow0}\xm_t(K_m)+\limsup_{t\rightarrow0}\xm_t(E\setminus K_m)\leq\xm(O)+\limsup_{t\rightarrow0}\xm_t(E\setminus K_m).$$
Now we assert that
$$\lim_{m\rightarrow\infty}\limsup_{t\rightarrow0}\xm_t(E\setminus K_m)=0.$$
We will prove it by contradiction. We assume that $\lim_{m\rightarrow\infty}\limsup_{t\rightarrow0}\xm_t(E\setminus K_m)=\xe>0.$\\
Let $\{t_n\}$ be a decreasing sequence tending to $0$ and $\lim_{n\rightarrow\infty}\xm_{t_n}(E\setminus K_m)=\limsup_{t\rightarrow0}\xm_t(E\setminus K_m).$ Then there exists subsequence of positive solutions $\{u_k^m\}_{k=1}^\infty$ with initial data $\xm_{t_{n_k}}\chi_{E\setminus K_m}$ such that $u_k^m\rightarrow u^m$ for any $m\in\mathbb{N}.$ Since $u$ is a moderate solution and $u_k^m\leq u$,  $u^m$ is a moderate solution too. Also by construction, the sequence $\{u^m\}$ is nonincreasing and $u_m\leq U_{E\setminus K_m}.$ By proposition \ref{sygklisi1} we have $U_{E\setminus K_m}\rightarrow0$ which implies $u_m\rightarrow0$ and $$\lim_{m\rightarrow\infty}\lim_{k\rightarrow\infty}\xm_{t_{n_k}}(E\setminus K_m)=0.$$
The proof follows in the case where $E$ is $\frak T_q$-closed. The proof is similar in the other case.\\
If $A$ is $\frak T_q$-open and
$$\xm(A)=\xm(\widetilde{A}),$$
then
$$\lim_{t\to 0}\xm_t(A)=\xm(A).$$

Without loss of generality we may assume that $\xf\geq0$ (since otherwise we set $\xf=\xf^+-\xf^-$) and $\xf\leq1.$ Given $k\in\mathbb{N}$ and  $m=0,...,2^k-1$ choose a number $a_{m,k}$ in the interval $(m2^{-k},(m+1)2^{-k})$ such that $\xm(\xf^{-1}(\{a_{m,k}\}))=0.$ Put
$$A_{m,k}=\xf^{-1}((a_{m,k},(a_{m+1,k}]),\;\;m=1,...,2^k-1,\qquad A_{0,k}=\xf^{-1}((a_{0,k},(a_{1,k}]),$$
then we note that since $\xf$ has compact support the above sets are bounded and
\be
\lim_{t\to 0}\xm_t(A_{m,k})=\xm(A_{m,k}),\forall m\geq0,\;k\in\mathbb{N}.\label{5.20}
\ee
Consider the step function
$\xf_k=\sum_{\xm=0}^{2^k-1}m2^{-k}\chi_{A_{m,k}},$
then $\xf_k\uparrow\xf$ uniformly, and by (\ref{5.20}),
$$\lim_{t\rightarrow0}\int_{\mathbb{R}^N}u(x,t)\xf_kdx=\int_{\mathbb{R}^N}\xf_kd\xm,\qquad\forall\xz\in C_0^\infty(\mathbb{R}^N).
$$
This completes the proof of (ii).\hfill$\Box$


\subsection{Vanishing properties}

\begin{defin}
A continuous function $u\in\mathcal{U}_+(Q_T)$ vanishes on a $\frak T_q$-open subset $G\subset\mathbb{R}^N,$ if for any $\eta\in\besl$ with $\frak T_q$-$\mathrm{supp}(\eta)\qsub G,$ there holds
\be
\lim_{t\rightarrow0}\int_{G}u(x,t)\eta^{2q'}_+(x)dxdt=0.\label{2.31}
\ee
When this is case we write $u\approx_G0.$ We denote by $\mathcal{U}_G(Q_T)$ the set of $u\in\mathcal{U}_+(Q_T)$ which vanish on $G.$\label{def2.8}
\end{defin}
We have the following simple result.
\begin{prop}
Let $A$ be a $\frak T_q$-open subset of $\mathbb{R}^N$ and $u_1,\;u_2\in\mathcal{U}_+(Q_T).$\\
If $u_2\approx_A0$ and $u_1\leq u_2$ then $u_1\approx_A0.$\label{anisothta}
\end{prop}
\begin{prop}
Let $G, G'$ be $\frak T_q$-open sets such that $G\qeq G'$. If $u\in\mathcal{U}_G(Q_T)$ then $u\in\mathcal{U}_{G'}(Q_T)$
\end{prop}
\Proof If $\eta\in\besl$ with $\frak T_q$-$\mathrm{supp}(\xz)\qsub G,$ then $\frak T_q$-$\mathrm{supp}(\xz)\qsub G'.$ Since $|G\setminus G'|=|G'\setminus G|=0$ the result follows.\hfill$\Box$\medskip

If $G$ is an open subset, this notion coincides with the usual definition of vanishing, since we can take test function $\eta \in C_0^\infty(G).$ In that case $u\in C(Q_T\cup\{G\times\{0\}\}).$
\begin{lemma}
Assume $u\in\mathcal{U}_G(Q_T).$ Then for any $\eta\in\besl$ with $\frak T_q$-$\mathrm{supp}(\eta)\qsub G,$ there holds
\be
\int_0^T\int_{\mathbb{R}^N}u^q\mathbb{H}[\eta]^{2q'}_+dxdt+\int_{\mathbb{R}^N}u(x,T)\mathbb{H}[\eta]_+^{2q'}(x,T)dx\leq C_1||\eta||^{q'}_{L^\infty}||\eta||^{q'}_{W^{\frac{2}{q},q'}}.\label{2.32}
\ee\label{lem2.9}
\end{lemma}
\Proof  If $u\in\mathcal{U}_G(Q_T)$ and $\eta\in\besl$ with $\frak T_q$-$\mathrm{supp}(\eta)\qsub G,$ there holds, with $h=\mathbb{H}[\eta]$ and $\xf(r)=r^{2q'}_+.$
\be
\int_0^T\int_{\mathbb{R}^N}(-u(\partial_t\xf(h)+\xD\xf(h)))+u^q\xf(h)dxd\tau+\int_{\mathbb{R}^N}u\xf(h)(.,T)dx=0.\label{2.33}
\ee
Therefore (\ref{2.32}) follows from (\ref{2.23}).\hfill$\Box$
\begin{lemma}
Let  $G\subset\BBR^N$ be a $\frak T_q$-open set. Then there exists a nondecreasing sequence $\{u_n\}\subset\mathcal{U}_G(Q_T)$ which converges to $\sup\mathcal{U}_G(Q_T).$ Furthermore $\sup\mathcal{U}_G(Q_T)\in\mathcal{U}_G(Q_T).$
\end{lemma}
\Proof    If $u_1$ and $u_2$ belongs to $\mathcal{U}_G(Q_T),$ then $u_1+u_2$ is a supersolution and it satisfies (\ref{2.31}). Therefore $u_1\vee u_2$ is a solution which is smaller than $u_1+u_2,$ thus $u_1\vee u_2\in \mathcal{U}_G(Q_T).$ By Proposition \ref{sygklish} there exists a increasing sequence $\{u_n\}\subset\mathcal{U}_G(Q_T)$ which converges to $u:=\sup\mathcal{U}_G(Q_T).$ By (\ref{2.33}),
\be
\int_0^T\int_{\mathbb{R}^N}(-u_n(\partial_t\xf(h)+\xD\xf(h)))+u^q_n\xf(h)dxd\tau+\int_{\mathbb{R}^N}u_n\xf(h)(.,T)dx=0.\label{2.35}
\ee
Now, $u_n^q\xf(h)\uparrow u^q\xf(h)$ in $L^1(Q_T)$ and $u_n\xf(h)(.,T)\uparrow u\xf(h)(.,T)$ in $L^1(\mathbb{R}^N).$ If $E$ is any Borel subset of $Q_T,$ there holds by H\"{o}lder's inequality, as in (\ref{2.23})
\be
\bigg|\int_0^T\int_{E}u_n(\partial_t\xf(h)+\xD\xf(h))dxd\tau\bigg|\leq C\left(\int_0^T\int_{E}u^q_n\xf(h)dxd\tau\right)^\frac{1}{q}
||\eta||_{L^\infty}||\eta||_{W^{\frac{2}{q},q'}}.\label{2.36}
\ee
The right-hand side tends to zero when $|E|\rightarrow0,$ thus by Vitali's convergence theorem, we derive
\be
\int_0^T\int_{\mathbb{R}^N}(-u(\partial_t\xf(h)+\xD\xf(h)))+u^q\xf(h)dxd\tau+\int_{\mathbb{R}^N}u\xf(h)(.,T)dx=0,\label{2.37}
\ee
from (\ref{2.35}). Thus $u\in\mathcal{U}_G(Q_T).$\hfill$\Box$

\begin{defin}
(a) Let $u\in\mathcal{U}_+(Q_T)$ and let $A$ denote the union of all $\frak T_q$-open sets on which $u$ vanishes. Then $A^c$ is called the fine initial support of $u,$ to be denoted by $\frak T_q$-supp$\,(u)$.\\
(b) Let $F$ be a Borel subset of $\mathbb{R}^N.$ We denote by $U_F$ the maximal element of $\mathcal{U}_{{\widetilde{F}}^{\,c}}(Q_T).$\label{4}
\end{defin}

\subsection{Maximal solutions}
\begin{defin}
Let $\mathfrak{M}^b_+(\mathbb{R}^N)$ be the set of all positive bounded Radon measures in $\mathbb{R}^N.$ Also let $u_\xm\in \mathcal{U}_+(Q_T)$ be the moderate solution with initial data $\xm.$\\
For any Borel set $E\subset\mathbb{R}^N$ of positive $\qcap$-capacity put
$$\mathcal{V}_{mod}(E)=\{u_\xm:\;\xm\in W^{-\frac{2}{q},q'}(\mathbb{R}^N)\cap\mathfrak{M}^b_+(\mathbb{R}^N),\;\xm(E^c)=0\}.$$
$$V_E=\sup V_{mod}(E).$$
\end{defin}

The following result due to Marcus and V\'eron \cite{MV-CVPDE} shows that the maximal solution which vanishes on an open set is indeed $\gs$-moderate. This is obtained by proving a capacitary quasi-representation of the solution via a Wiener type test.
\begin{prop}\label{est} Let $F$ be a closed subset of $\mathbb{R}^N$ and $q\geq 1+\frac{2}{N}$. Then there exist two positive constants $C_1,\;C_2 > 0,$ depending only
on $N$ and $q$ such that
\begin{equation}\label{9}\BA{lll}\displaystyle
C_1t^{-\frac{1}{q-1}}\sum_{k=0}^\infty\left(k+1\right)^{\frac{N}{2}}e^{-\frac{k}{4}}
\qcap\left(\frac{F\cap F_k(x,t)}{\sqrt{(k+1)t}}\right)\leq V_F(x,t)\leq U_F(x,t)\\\phantom{-------}
\displaystyle\leq C_2t^{-\frac{1}{q-1}}\sum_{k=0}^\infty\left(k+1\right)^{\frac{N}{2}}e^{-\frac{k}{4}}
\qcap\left(\frac{F\cap F_k(x,t)}{\sqrt{(k+1)t}}\right)\quad\forall (x,t)\in Q,
\EA
\end{equation}
where $F_k(x,t)=\{y\in\mathbb{R}^N:\;\sqrt{kt}\leq|x-y|\leq\sqrt{(k+1)t}\}.$ As a consequence
$U_F=V_F.$
\end{prop}

\noindent\Remark We recall that the main argument for proving uniqueness is the fact that
\begin{equation}\label{equiv}
U_F\leq \frac{C_2}{C_1}V_F\qquad\text{ in }Q.
\end{equation}
This argument introduced in \cite{MV-ARMA} for elliptic equations has been extended to parabolic equations in \cite{MV-CPDE}, \cite{MV-CVPDE}.
\begin{defin}
Let $F$ be a Borel subset of $\mathbb{R}^N.$ We denote by $U_F$ the maximal element of $\mathcal{U}_{\widetilde{F}^c}(Q_T).$
\end{defin}
\begin{prop}
If $\{A_n\}$ is a sequence of Borel sets such that $\qcap(A_n)\rightarrow0,$ then $U_{A_n}\rightarrow0.$\label{sygklisi1}
\end{prop}
\Proof    Let $O_n$ be an open set such that $A_n\subset O_n$  and $\qcap(O_n)\leq\qcap(A_n)+\frac{1}{n}.$ Now since $O_n$ is open, $\qcap$ is an outer measure, by (\ref{kellog}) and (iv)-Proposition \ref{propertiescapacity}, we have
$$\qcap(\overline{O}_n)=\qcap\left((\overline{O}_n\cap b_q(O_n))\bigcup(\overline{O}_n\cap e_q(\overline{O}_n))\right)\leq \qcap(\widetilde{O}_n)\leq c\qcap(O_n).$$
Thus $\qcap(\overline{O}_n)\rightarrow0.$ The result follows by
$$U_{A_n}\leq U_{\overline{O}_n}$$
and
by (\ref{9}).\hfill$\Box$
\begin{coro}
Let $E$ be a Borel set such that $\qcap(E)=0.$ If $u\in\mathcal{U}_{\widetilde{E}^c}(Q_T)$ then $u=0.$ In particular $U_E\equiv0.$\label{coro0cap}
\end{coro}
\begin{prop}
Let $E,\;F$ be Borel sets.\\
(i) If $E,\;F$ are $\frak T_q$-closed, then $U_E\wedge U_F=U_{E\cap F}.$\\
(ii) If $E,\;F$ are $\frak T_q$-closed, then
\bea
\nonumber
U_E<U_F&\Leftrightarrow& [E\qsub F\;\mathrm{and}\;\qcap(F\setminus E)>0],\\
U_E=U_F&\Leftrightarrow& \;E\qeq F.
\eea
(iii) If $F_n$ is a decreasing sequence of $\frak T_q$-closed sets, then
$$\lim_{n\to\infty} U_{F_n}=U_F\;\;\mathrm{where}\;\;F=\cap F_n.$$
(iv) Let $A$ be a $\frak T_q$-open set and $u\in\mathcal{U}_+(Q_T).$ Suppose that $u$ vanishes $\frak T_q$-locally in $A,$ i.e. for every point $\xs\in A$ there exists a $\frak T_q$-open set $A_\xs$ such that
$$\xs\in A\subset A,\;\;\;u\approx_{A_\xs}0.$$
Then $u$ vanishes on $A.$ In particular any $u\in\mathcal{U}_+(Q_T)$ vanishes on the complement of $\frak T_q$-supp$\,(u)$.\label{7}
\end{prop}
\Proof    The proof is similar to the one in \cite{MV-CONT} dealing with elliptic equations.\\
(i) $U_E\wedge U_F$ is the largest solution under $\inf(U_E, U_F)$ and therefore, by definition, it is the largest solution which vanishes outside $E\cap F.$\\
(ii) By (\ref{9}) $U_E$ and $U_F$ satisfies the same capacitary quasi-representation up to universal constants. By the Remark after Proposition \ref{est} ,
$$
E\qeq F\Rightarrow \myfrac{C_1}{C_2}U_E\leq U_F\leq  \myfrac{C_2}{C_1}U_E\Rightarrow U_E=U_F.$$
The proof of
$$E\qsub F\Rightarrow U_E\leq U_F.
$$
follows from Proposition \ref{est} and the fact that $U_E=V_E$ and $U_F=V_F$ and $V_E\leq V_F$. In addition,
$$\qcap(F\setminus E)>0\Rightarrow U_E\neq U_F.$$
Indeed, if $K$ is a compact subset of $F\setminus E$ of positive capacity, then $U_K>0$ and $U_K\leq U_F$ but $U_K\nleq U_E.$ Therefore $U_E=U_F$ implies $E\qeq F$ and $U_E\leq U_F$ implies $E\qsub F.$\\
(iii) If $V:=\lim_{n\to\infty} U_{F_n}$ then $U_F\leq V.$ But $\frak T_q$-supp$\,(V)\subset F_n$ for each $n\in\mathbb{N}$ and consequently $V\leq U_F.$\\
(iv) First assume that $A$ is a countable union of $\frak T_q$-open sets $\{A_n\}$ such that $u\approx_{A_n}0$ for each $n.$ Then $u$ vanishes on $\cup_{i=1}^kA_k$ for each $k.$ Therefore we can assume that the sequence $A_k$ is increasing. Put $F_n=A_n^c.$ Then $u\subset U_{F_n}$ and by (iii), $U_{F_n}\downarrow U_F$ where $F=A^c.$ Thus $u\leq U_F,$ i.e.,which is equivalent to $u\approx_A0.$

We turn to the general case. It is known that the $\frak T_q$-fine topology possesses the {\it quasi-Lindel\"{o}f property }(see \cite[Sec. 6.5.11]{AH}) as any topology associated to a Bessel capacity $C_{\ga,p}$. Therefore $A$ is covered, up to a set of capacity zero, by a countable subcover of $\{A_\xs:\;\xs\in A\}.$ Therefore the previous argument implies that $u\approx_A0.$\hfill$\Box$
\begin{prop}
(i) Let $E$ be a $\frak T_q$-closed set. Then
\bea
\nonumber
U_E&=&\inf\{U_D:\; E\subset D,\;D\;\mathrm{open}\}\\
&=&\sup\{U_K:\;K\subset E,\;K\;\mathrm{closed}\}.\label{5}
\eea
(ii) If $E,\;F$ are two Borel sets then
$$U_E=U_{F\cap E}\oplus U_{E\setminus F}.$$
(iii) Let $E,\;F_n,\;n=1,2,... $ be Borel sets and let $u$ be a positive solution of (\ref{maineq}). If either $\qcap(E\bigtriangleup F_n)\rightarrow0$ or $\widetilde{F}_n\downarrow \widetilde{E}$ then
$$U_{F_n}\rightarrow U_E.$$\label{12}
\end{prop}
\Proof    (i) Let $\{Q_j\}$ be the decreasing sequence of open sets of Lemma \ref{6}-(I) such that $\cap Q_j=\cap\widetilde{Q}_j=E'\qeq E.$ Thus by Proposition \ref{7} (iii) we have that $U_{Q_j}\rightarrow U_E,$ this implies the first equality in (i).

Let $\{F_n\}$ be a nondecreasing sequence of closed subsets of $E$ such that $\qcap(E\setminus F_n)\rightarrow0.$ Let $D_1$ and $D_2$ be open sets such that $F_n\subset D_1$ and $E\setminus F_n\subset D_2.$ Also set $D_3=(\widetilde{D}_1\cup \widetilde{D}_2)^c.$ Let $u_\xb^{(i)}$ be the positive solution of
\bea
\nonumber
\prt_t  u-\xD u+u^q &=&0\quad\quad\quad\quad\;\mathrm{in }\;\;\;\mathbb{R}^N\times(\xb,T]\\
u(.,\xb)&=& \chi_{\widetilde{D}_i}U_E(.)\;\;\;\;\mathrm{on }\;\;\mathbb{R}^N,
\eea
where $0<\xb< T.$
For any $(x,t)\in\mathbb{R}^N\times(\xb,T]$ we have
$$U_E\leq u_\xb^{(1)} +u_\xb^{(2)}+ u_\xb^{(3)}.$$
Letting $\xb\rightarrow0$ (taking an subsequence if it is necessary) we have $ u_\xb^{(i)}\rightarrow u^{(i)}$ and
$$U_E\leq u^{(1)}+ u^{(2)} +u^{(3)}\qquad\mathrm{in} \;\;Q_T,$$
But $u^{(i)}\leq U_{D_i}$ thus
$$U_E\leq U_{D_1}+U_{D_2}+u^{(3)}.$$
Now $u^{(3)}\leq U_{D_3}$ and $u^{(3)}\leq U_E$ thus by Proposition \ref{12}-(i) $u^{(3)}\leq U_{D_3\cap E}.$ But $D_1\cup D_2$ is an open set and thus $\qcap(D_3\cap E)=0,$ which implies by Corollary \ref{coro0cap}  that $u^{(3)}=0.$ Finally we have that
$$U_E\leq U_{D_1}+U_{D_2}.$$
Since $D_i$ is arbitrary, we have by the first assertion of this Proposition
\be
U_E\leq U_{F_n}+ U_{E\setminus F_n}.\label{1*}
\ee
But $\qcap(E\setminus F_n)\rightarrow0,$ thus by Proposition \ref{sygklisi1},
we have
$$U_E\leq\lim_{n\rightarrow\infty}U_{F_n}\Rightarrow U_E=\lim_{n\rightarrow\infty}U_{F_n},$$
since $U_{F_n}\leq U_E$ for any $n\in\mathbb{N}.$\\
(ii) By similar argument as in the proof of (\ref{1*}) we can prove that
$$U_E\leq U_{F\cap E}+U_{E\setminus F}\Rightarrow U_E\leq U_{F\cap E}\oplus U_{E\setminus F}.$$
On the other hand both $U_{F\cap E}$ and $U_{E\setminus F}$ vanish outside of $\widetilde{E}.$ Consequently $U_{F\cap E}\oplus U_{E\setminus F}$ vanishes outside $\widetilde{E}$ so that
$$U_E\geq U_{F\cap E}\oplus U_{E\setminus F},$$
and the result follows in this statement.\\
(iii) The previous statement implies,
\be
U_E\leq U_{F_n\cap E}+U_{E\setminus F_n},\qquad U_{F_n}\leq U_{F_n\cap E} +U_{F_n\setminus E}.\label{11}
\ee
If $\qcap(E\bigtriangleup F_n)\rightarrow0$ then Proposition \ref{sygklisi1} implies $U_{E\bigtriangleup F_n}\rightarrow0.$ And the result follows in this case by (\ref{11}).\\
If $\widetilde{F}_n\downarrow \widetilde{E}$ the result follows in this case by Proposition \ref{7}(iii).\hfill$\Box$\\

This implies the following extension of Proposition \ref{est} to merely $\frak T_q$-closed sets.
\begin{prop}
If $E$ is a $\frak T_q$-closed set, then $V_E$ and $U_E$ satisfy the capacitary estimates (\ref{est}). Thus
$U_E=V_E$
and the maximal solution $U_E$ is $\xs$-moderate.\label{31}
\end{prop}
\emph{Remark.} Actually the estimates hold for any Borel set $E.$ Indeed by definition, $U_E=U_{\widetilde{E}}$ and
$$\qcap\left(\frac{E\cap F_n(x,t)}{\sqrt{(n+1)t}}\right)\sim\qcap\left(\frac{\widetilde{E}\cap F_n(x,t)}{\sqrt{(n+1)t}}\right).$$
\Proof    The proof is same as in \cite{MV-CONT}.\\
Let $\{E_k\}$ be a $\frak T_q$-stratification of $E.$ If $u\in\mathcal{V}_{mod}$ and $\xm=\mathrm{tr}u$ then $u_\xm=\sup u_{\xm_k}$ where $\xm_k=\xm\chi_{E_k}.$ Hence $V_E=\sup V_{E_k}.$ By proposition \ref{9}, $U_{E_k}=V_{E_k}.$ These facts and Proposition \ref{12}(c) implies $U_E=V_E.$ Since $U_{E_k}$ satisfies the capacitary estimates (\ref{est}) and
$$\qcap\left(\frac{E_k\cap F_n(x,t)}{\sqrt{(n+1)t}}\right)\rightarrow\qcap\left(\frac{\widetilde{E}\cap F_n(x,t)}{\sqrt{(n+1)t}}\right)\quad\text{as }n\to\infty.$$
it follows that $U_E$ satisfies the corresponding capacitary estimates.\hfill$\Box$
\subsection{Localization}
\begin{defin}
Let $A$ be a Borel subset of $\mathbb{R}^N,$ we denote by $[u]_A$ the supremum of the $v\in \mathcal{U}_+(Q_T)$ which are dominated by $u$ and vanishes on $\widetilde{A}^c.$
\end{defin}
We note here that $[u]_A=u\wedge U_A$
\begin{lemma}
If $G\subset\mathbb{R}^N$ is a $\frak T_q$-open set and $u\in\mathcal{U}_G(Q_T),$ then
$$u=\sup\{v\in\mathcal{U}_G(Q_T):\;v\leq u,\;v\; \mathrm{vanishes\;on\;an\;open\;neighborhood\;of\;}G\}.$$\label{3}
\end{lemma}
\Proof
 Set $A=G^c$ and let $\{A_n\}$ be a sequence of closed subsets of $A$, such that $\qcap(A\setminus A_n)\rightarrow0.$ By Proposition \ref{12} we have
 $$U_A\leq U_{A_n}+U_{A\setminus A_n},$$
 thus
 $$u=u\wedge U_A\leq u\wedge U_{A_n}+u\wedge(U_{A\setminus A_n}).$$
By Proposition \ref{sygklisi1}, we have
$$U_{A\setminus A_n}\rightarrow0.$$
Thus $$u=\lim_{n\rightarrow\infty}u\wedge U_{A_n},$$
and the result follows.\hfill$\Box$\medskip

The next result points out the set-regularity of the correspondence $E\mapsto [u]_E$.

\begin{prop}
Let $u\in\mathcal{U}_+(Q_T).$\\
(i) If $E$ is $\frak T_q$-closed then,
\bea
[u]_E&=&\inf\{[u]_D:\;E\subset D,\;D\mathrm{\;open}\}.\label{10}\\
     &=&\sup\{[u]_F:\;F\subset E,\;F\mathrm{\;closed}\}. \label{16}
\eea
(ii) If $E,\;F$ are two Borel sets then
\be
[u]_E\leq[u]_{F\cap E}+[u]_{E\setminus F},\label{13}
\ee
and
\be
[[u]_E]_F=[[u]_F]_E=[u]_{F\cap E}.\label{14}
\ee
(iii) Let $E,\;F_n,\;n=1,2,... $ be Borel sets and let $u$ be a positive solution of (\ref{maineq}). If either $\qcap(E\bigtriangleup F_n)\rightarrow0$ or $\widetilde{F}_n\downarrow \widetilde{E}$ then
$$[u]_{F_n}\rightarrow [u]_E.$$\label{17}
\end{prop}
\Proof    The proof uses a similar argument as in \cite{MV-CONT}.\\
(i) Let $\mathcal{D}=\{D\}$ be the family of sets in (\ref{10}). By (\ref{5}) (with respect to the family $\mathcal{D}$)
\be
\inf(u,U_E)=\inf(u,\inf_{D\in\mathcal{D}}U_D)=\inf_{D\in\mathcal{D}}\inf(u,U_D)\geq\inf_{D\in\mathcal{D}}[u]_D.\label{15}
\ee
Obviously
$$[u]_{D_1}\wedge[u]_{D_2}\geq[u]_{D_1\cap D_2},$$
thus we can apply Proposition \ref{sygklish} and obtain that the function $v:=\inf_{D\in\mathcal{D}}[u]_D$ is a solution of (\ref{maineq}). Hence (\ref{15}) implies $[u]_E\geq v.$ The opposite inequality is obvious.

For the equality (\ref{16}), Firstly, we note that the set  $\{v\in\mathcal{U}_+(Q_T):\;u\leq u,\;\frak T_q\text{-supp}\,(v)\qsub E\}$ is closed under $\vee.$ Thus, by Proposition \ref{sygklish}, there exists an increasing sequence $\{v_n\}$ such that $v_n\approx_{E^c}=0$ and $\lim_{n\to\infty} v_n=[u]_E.$ Since $v_n$ is an increasing sequence by Proposition \ref{3} we can construct an increasing sequence $\{w_n\}$ such that each $w_n$ vanishes on an open neighborhood $B_n$ of $E,$ $B_n\subset B_{n+1}$ and $\lim_{n\to\infty} w_n=[u]_E.$ Now set $K_n=B_n^c,$ then
$$w_n\leq[u]_{K_n}\leq[u]_E.$$
Letting $n$ tend to infinity, we obtain the desired result.\\
(ii) Let $v\in\mathcal{U}_+(Q_T),\;v\leq u$ and $\frak T_q$-supp$\,(v)\subset E.$
Let $D$ and $D'$ be open sets such that $\widetilde{E\cap F}\subset D$ and $\widetilde{E\setminus F}\subset D'.$ By Lemma 2.8-\cite{MV-CPDE}, there exists a unique solution $v^1_j$, where $\frac{1}{[T]}<j\in\mathbb{N}$, of the problem
\bea
\nonumber
\prt_t  u-\xD u+|u|^{q-1}u&=&0,\qquad\qquad\qquad\,\,\;\mathrm{in}\;\mathbb{R}^N\times(\frac{1}{j},T]\\ \nonumber
u(.,\frac{1}{j})&=& \chi_D(.)v(.,\frac{1}{j})\qquad\,\,\mathrm{in}\;\mathbb{R}^N.
\eea
Also we consider $v^2_j$ and $v^3_j$ the unique solutions of the above problem with initial data $\chi_{D'}(x)v(x,\frac{1}{j})$ and $\chi_{(D_1\cup D_2)^c}.$
In view of the proof of Proposition \ref{12} we can prove that $v\leq v^1_j+v^2_j+v^3_j.$ By standard arguments there exists a subsequence, say $\{v^{i}_j\}$, $i=1,2,3$, such that $v^i_j\rightarrow v^i$ and $v\leq v^1+v^2+v^3.$ Since $v$ vanishes outside of $E,$ it vanishes outside of $(D_1\cup D_2),$ consequently $v(x,\frac{1}{j})\chi_{\chi_{(D_1\cup D_2)^c}}\rightarrow0,$ as $j\to\infty$, which implies $v^3_j\rightarrow0.$ Thus we have
$$v\leq v^1+v^2\leq[u]_D+[u]_{D'}.$$
By (\ref{10}) we have
$$v\leq[u]_{F\cap E}+[u]_{E\setminus F},$$
since $v\in\{w\in\mathcal{U}_+(Q_T):\;w\leq u,\;\frak T_q\text{-supp}\,(w)\qsub E\}$ is arbitrary the result follows in the case where $E$ is closed. In the general case the result follows by (\ref{16}).

Put $A=\widetilde{E}$ and $B=\widetilde{F}.$ It follows directly from the definition that
$$[[u]_A]_B\leq\inf(u,U_A,U_B).$$
The largest solution dominated by $u$ and vanishing on $A^c\cup B^c$ is $[u]_{A\cap B}.$ Thus
$$[[u]_A]_B\leq[u]_{A\cap B}.$$
On the other hand $$[u]_{A\cap B}=[[u]_{A\cap B}]_B\leq[[u]_A]_B,$$
this proves (\ref{14}).\smallskip

\noindent(iii) By (\ref{13})
$$[u]_E\leq[u]_{F_n\cap E}+[u]_{E\setminus F_n},\qquad[u]_{F_n}\leq[u]_{F_n\cap E}+[u]_{F_n\setminus E}.$$
If $\qcap(E\bigtriangleup F_n)\rightarrow0,$ then by Proposition (\ref{sygklisi1})(c) we have that $U_{E\bigtriangleup F_n}\rightarrow0.$ Since $[u]_{E\setminus F_n},\;[u]_{F_n\setminus E}\leq U_{E\bigtriangleup F_n},$ the result follows by the above inequalities, if we let $n$ go to infinite.

If $\widetilde{F}_n\downarrow \widetilde{E}.$
By Proposition (\ref{sygklisi1})(c) we have $U_{E_n}\rightarrow U_E,$ thus
$$[u]_E\leq\lim_{n\to\infty}[u]_{F_n}=\lim_{n\to\infty} u\wedge U_{F_n}\leq\lim_{n\to\infty}\inf(u,U_{F_n})\leq\inf(u,U_E).$$
And since $[u]_E$ is the largest solution under $\inf(u,U_E)$ and the function $v=\lim_{n\to\infty}[u]_{F_n}$ is a solution of (\ref{maineq}), we have that $U_E\leq v,$ and the proof of (\ref{13}) is complete.\hfill$\Box$
\begin{defin}
Let $\xm$ be a positive Radon measure on $\mathbb{R}^N$ which vanishes on compact sets of $\qcap$-capacity zero.\\
(a) The $\frak T_q$-support of $\xm$, denoted $\frak T_q$-$\mathrm{supp}(\xm)$, is the intersection of all $\frak T_q$-closed sets $F$ such that $\xm(F^c)=0.$\\
(b) We say that $\xm$ is concentrated on a Borel set $E$ if $\xm(E^c)=0.$
\end{defin}
\begin{prop}
If $\xm$ is a measure as in the previous definition then,
$$\frak T_q\text{-supp}\,(\xm)\qeq \frak T_q\text{-supp}\,(u_\xm).$$\label{33}
\end{prop}
\Proof
Put $F=\mathrm{supp}^qu_\xm.$ By Proposition \ref{7}(iv) $u_\xm$ vanishes on $F^c$ and by Proposition \ref{3}(c) there exists an increasing sequence of positive solutions $u_n$ such than each function $u_n$ vanishes outside a closed subset $F,$ say $F_n,$ and $u_n\uparrow u_\xm.$ If $S_n:=\frak T_q\text{-supp}\,(u_n)$ then $S_n\subset F_n$ and $\{S_n\}$ increases. Thus $\{\overline{S}_n\}$ is an increasing sequence of closed subsets of $F$ and, setting $\xm_n=\xm\chi_{\overline{S}_n},$ we find $u_n\leq u_{\xm_n}\leq u_\xm$ so that $u_{\xm_n}\uparrow u_\xm.$ This, in turn, implies
$$\xm_n\uparrow\xm,\qquad \frak T_q\text{-supp}(\xm)\qsub\widetilde{\bigcup_{n=1}^\infty\overline{S}_n}\subset F.$$
If $D$ is an open set and $\xm(D)=0$ it is clear that $u_\xm$ vanishes on $D.$ Therefore $u_{\xm_n}$ vanishes outside of $\overline{S}_n,$ thus outside $\frak T_q$-$\mathrm{supp}\,(\xm).$ Consequently $u_\xm$ vanishes outside $\frak T_q$-$\mathrm{supp}(\xm),$ i.e., $F\qsub \frak T_q$-supp$\,(\xm).$\medskip

\noindent{\it Second proof.}
The result follows by Proposition \ref{moderate} and Definition \ref{def2.8}\hfill$\Box$
\begin{defin}
Let $u$ be a positive solution and $A$ a Borel set. Put
$$[u]^A:=\sup\{[u]_F:\;F\qsub A,\;F\; q\mathrm{-closed}\}.$$\label{32}
\end{defin}
\begin{defin}
Let $\xb>0,$ $u\in C(Q_T).$ For any Borel set $A$ we denote by $u^A_\xb$ the positive solution of
\bea
\nonumber
\prt_t  v-\xD v+|v|^{q-1}v&=&0\,\qquad\qquad\qquad\;\;\mathrm{in}\;\mathbb{R}^N\times(\xb,\infty)\\ \nonumber
v(.,\xb)&=& \chi_A(.)u(.,\xb)\qquad\mathrm{in}\;\mathbb{R}^N.
\eea
\end{defin}
\begin{prop}
Let $u$ be a positive solution of (\ref{maineq}) and put $E=\frak T_q\text{-supp}\,(u).$\\
(i) If $D$ is a $\frak T_q$-open set such that $E\qsub D$, then
\be
[u]^D=\lim_{\xb\rightarrow0} u^D_\xb=[u]_D=u.\label{4.9}
\ee
(ii) If A is a $\frak T_q$-open set, then
\be
u\approx_{A}0\Leftrightarrow u^Q=\lim_{\xb\rightarrow0}u_\xb^Q=0,\qquad\forall Q\;\frak T_q\mathrm{-open}:\;\widetilde{Q}\qsub A.\label{4.10}
\ee
(iii) Finally,
\be
u\approx_{A}0\Leftrightarrow[u]^A=0.\label{4.11}
\ee\label{lem4.5}
\end{prop}
\Proof    The proof is similar as in the one as in \cite{MV-CONT}\\
{\it Case 1}: $E$ is closed. Since $u$ vanishes in $E^c,$ it yields $u\in C(Q_\infty\cup E^c)$ and $u=0$ on $E^c.$ If, in addition, $D$ is an open neighborhood of $E$, we have $$\lim_{t\rightarrow0}\int_{E^c}\xf(x)u(x,t)dx=0,\qquad\forall\xf\in C_0(E^c).$$
Thus,
$$\lim_{\xb\to 0} u_\xb^{D^c}=0.$$
Since
$$u_\xb^D\leq u\leq u_\xb^D+u_\xb^{D^c},\qquad\forall t\geq\xb,$$
it follows
\be
u=\lim_{\xb\to 0} u_\xb^D.\label{4.13}
\ee
If we assume that $D$ is $\frak T_q$-open and $E\qsub D$ then, for every $\xe>0,$ there exists an open set $O_\xe$ such that $D\subset O_\xe,\;E\subset O_\xe$ and $\qcap(O_\xe')<\xe$ where $O_\xe'=O_\xe\setminus D.$ Therefore
$$u_\xb^{O_\xe}(x,t)-u_\xb^D(x,t)\leq U_{O_\xe'}(x,t-\xb),\qquad\forall t\geq\xb.$$
We note here that $\lim_{\xe\rightarrow0} U_{O_\xe'}(x,t-\xb)=0$ holds uniformly with respect to $\xb.$ Since $\lim_{\xb\rightarrow0} u_\xb^{O_\xe}(x,t)=u$
it follows that $u=\lim_{\xb\to 0} u_\xb^D.$ The same arguments shows that $\lim_{\xb\to 0} u_\xb^{D^c}=0.$ Thus we have
$$u=\lim_{\xb\to 0} u_\xb^D\leq[u]_D\leq u.$$
Hence $u=[u]_D.$ By Lemma \ref{q-open}, there exists a $\frak T_q$-open set $Q$ such that $E\qsub Q\subset\widetilde{Q}\qsub D,$ then $u=[u]_Q\leq[u]^D.$ Hence $u=[u]^D.$

In addition, there holds $E\qsub A^c\qsub\widetilde{Q}^c.$ Thus the direction "$\Rightarrow$" in  (\ref{4.10}) follows by the previous argument if we replace $D$ by $\widetilde{Q}^c.$ For the opposite direction, by Proposition \ref{perioxi}, for any $\xi\in A,$ there exists a $\frak T_q$-open set $O_\xi$ such that $\widetilde{O}_\xi\qsub A.$ Using (i) we infer $u=\lim_{\xb\to 0} u^{\widetilde{O}_\xi^c}_\xb.$ Finally, since $u^{\widetilde{O}_\xi^c}_\xb\approx_{Q_\xi}0$ for all $\xb>0,$ it implies $u\approx_{O_\xi}0$ by Proposition \ref{sygklisi1}(i), and the result follows in this case by Proposition \ref{7}(iv).\smallskip

\noindent{\it Case 2}. Assume $E$ is $\frak T_q$-closed. Let $\{E_n\}$ be a $\frak T_q$-stratification of $E$ such that $\qcap(E\setminus E_n)\rightarrow0.$ If $D$ is a $\frak T_q$-open such that $E\qsub D$ then, by the first case we have,
\be
\lim_{\xb\rightarrow 0}([u]_{E_n})^D_\xb=[u]_{E_n}.\label{4.14}
\ee
By (\ref{13}) and the definition of $u_\xb^D$, and since $[u]_E=u$,
\be
u_\xb^D=\left([u]_E\right)_\xb^D\leq\left([u]_{E\cap E_n}\right)_\xb^D+\left([u]_{E\setminus E_n}\right)_\xb^D=\left([u]_{E_n}\right)_\xb^D
+\left([u]_{E\setminus E_n}\right)_\xb^D.\label{4.15}
\ee
Let $\{\xb_k\}$ be a decreasing sequence converging to $0$ such that the following limits exist
$$w:=\lim_{k\rightarrow\infty}u_{\xb_k}^D,\qquad w_n=\lim_{k\rightarrow\infty}\left([u]_{E\setminus E_n}\right)_{\xb_k}^D,\qquad n=1,2,...\;.$$
Then by (\ref{4.14}) and (\ref{4.15}),
$$[u]_{E_n}\leq w\leq[u]_{E_n}+w_n\leq[u]_{E_n}+U_{E\setminus E_n}.$$
Further, by (\ref{16}) and Proposition \ref{12}(c)
$$[u]_{E_n}\rightarrow[u]_E=u,\qquad U_{E\setminus E_n}\rightarrow0.$$
Hence $w=u.$ This implies (\ref{4.13}), which in turn implies (\ref{4.9}).\\
To verify (\ref{4.10}) in the direction $\Rightarrow$ we apply (\ref{4.15}) with $D$ replaced by $Q.$
We obtain
$$\left([u]_E\right)_\xb^{Q}\leq\left([u]_{E_n}\right)_\xb^{Q}
+\left([u]_{E\setminus E_n}\right)_\xb^{Q}.$$
By the first case we have
$$\lim_{\xb\rightarrow0}\left([u]_{E_n}\right)_\xb^{Q}=0.$$
 There exists a decreasing sequence converging to $0$, still denoted by  $\{\xb_k\}$,  such that the following limits exist
$$\lim_{k\rightarrow\infty}u_{\xb_k}^Q,\qquad \lim_{k\rightarrow\infty}\left([u]_{E\setminus E_n}\right)_{\xb_k}^Q,\qquad n=1,2,...\;.$$
Then
$$\lim_{k\rightarrow\infty}u_{\xb_k}^Q\leq\lim_{k\rightarrow\infty}\left([u]_{E\setminus E_n}\right)_{\xb_k}^Q\leq U_{E\setminus E_n},$$
since $U_{E\setminus E_n}\rightarrow0$ we obtain (\ref{4.10}) in the direction $\Rightarrow.$ The assertion in the opposite direction is proved as in Case 1. This complete the proofs of (i) and (ii).

Finally we prove (iii). First assume that $u\approx_A0.$ If $F$ is a $\frak T_q$-closed set such that $F\qsub A,$ then by Lemma \ref{q-open} there exists
a $\frak T_q$-open set $\frak T_q$ such that $F\qsub Q\subset\widetilde{Q}\qsub A.$ Therefore, applying (\ref{4.9}) to $v:=[u]_F$ and using (\ref{4.10}) we obtain
$$v=\lim_{\xb\to 0} v_\xb^Q\leq\lim_{\xb\to 0} u_\xb^Q=0.$$
By definition of $[u]^A,$ this implies $[u]^A=0.$

If $[u]^A=0,$ then for any $\frak T_q$-open set $Q\subset\widetilde{Q}\qsub A$ there holds $[u]_Q=0.$. Now since $\frak T_q\text{-supp}\,(u_\xb^Q)\qsub \widetilde{Q}$ we have for some subsequence $\xb_k\downarrow0,$ $\lim_{k\rightarrow\infty}u^Q_{\xb_k}\leq[u]_Q=0.$ Thus $u\approx_Q0$  by (\ref{4.10}). Applying once again Proposition \ref{perioxi} and Proposition \ref{7}(iv) we conclude $u\approx_A0.$\hfill$\Box$
\begin{defin}
Let $u,\;v\in\mathcal{U}_+(Q_T)$ and let $A$ be a $\frak T_q$-open set. We say that $u=v$ on $A$ if $u\ominus v$ and $v\ominus u$ vanishes on $A.$ This relation is denoted by $u\approx_A v.$
\end{defin}
\begin{prop}
Let $u,\;v\in\mathcal{U}_+(Q_T)$ and let $A$ be a $\frak T_q$-open set. Then,\\
(i)
\be
u\approx_A v\Leftrightarrow\lim_{\xb\rightarrow0}|u-v|^Q_\xb=0.\label{4.16}
\ee
for every $\frak T_q$-open set $Q$ such that $\widetilde{Q}\qsub A.$\\
(ii)
\be
u\approx_A v\Leftrightarrow[u]_F=[v]_F,\label{4.17}
\ee
for every $\frak T_q$-closed set $F$ such that $F\qsub A.$
\end{prop}
\Proof    The proof is similar, but in a parabolic framework, to the elliptic one in \cite{MV-CONT}.\\
By definition $u\approx_Av$ is equivalent to $u\ominus v\approx_A=0$ and $v\ominus u\approx_A=0.$ Hence, by (\ref{4.10}) we have
$w_\xb=(u\ominus v)^Q_\xb\rightarrow_{\xb\rightarrow0}0.$ Set $f_\xb=\left((u-v)_+\right)^Q_\xb$ and consider the problem
\bea
\nonumber
\prt_t  w-\xD w +|w|^{q}&=&0,\qquad\mathrm{in}\;B_j(0)\times(\xb,\infty)\\ \nonumber
w&=&0,\qquad\mathrm{on}\;\partial B_j(0)\times(\xb,\infty)\\ \nonumber
w(.,\xb)&=&\xm,\qquad\mathrm{in}\;B_j(0).
\eea
Let $w_j$ and $f_j$ be solutions of the above problem, with initial data $\chi_{Q}(u\ominus v)(x,\xb)$ and $\chi_{Q}(u-v)_+(x,\xb).$ By \cite[Lemma 2.7]{MV-CPDE}, the sequences $\{w_j\}$ and $\{f_j\}$ are increasing. Also, we recall that $u\ominus v$ is the smallest solution which dominates the subsolution $(u-v)_+,$ thus $w_j\geq v_j,\;\forall j\in\mathbb{N}.$ Furthermore, in view of \cite[Lemma 2.8]{MV-CPDE}, there holds $\lim_{j\to\infty} w_j=w_\xb$ and $\lim_{j\to\infty} f_j=f_\xb.$ Thus $w_\xb\geq f_\xb$, and letting $\xb\to 0$ we derive
$$\left((u-v)_+\right)^Q_\xb\rightarrow0.$$
By the same argument we have
$$\left((v-u)_+\right)^Q_\xb\rightarrow0,$$
this implies (\ref{4.16}) in the direction $\Rightarrow$.

For the opposite direction, we consider the problem
\bea
\nonumber
\prt_t  w-\xD w +|w|^{q}&=&0,\qquad\mathrm{in}\;B_j(0)\times(\xb,\infty)\\ \nonumber
w&=&h,\qquad\mathrm{on}\;\partial B_j(0)\times(\xb,\infty)\\ \nonumber
w(.,\xb)&=&\xm,\qquad\mathrm{in}\;B_j(0).
\eea
Let $Q\subset\widetilde{Q}\qsub A $ be a $\frak T_q$-open set and $w_j$ be the solution of the above problem, with $h=\chi_{Q}(|u-v|)$ and $\xm=\chi_{Q}|u-v|dx.$ Also, let $f_j$ be the solution of the above problem with $h=\chi_{Q^c}|u-v|$ and $\xm=\chi_{Q^c}|u-v|dx,$ then
$$|u-v|\leq w_j+f_j.$$
In view of \cite[Lemma 2.8]{MV-CPDE}, there exist subsequences, say $\{w_j\}$ and $\{f_j\}$, satisfying $\lim_{j\to\infty} w_j=w$ and $\lim_{j\to\infty} f_j=f,$ such that $(w,f)$ solves the problem
\bea
\nonumber
\prt_t  v-\xD v+|v|^{q-1}v&=&0,\qquad\qquad\mathrm{in}\;\mathbb{R}^N\times(\xb,\infty)\\ \nonumber
v(.,\xb)&=& \xm\qquad\qquad\;\,\mathrm{in}\;\mathbb{R}^N,
\eea
with initial data $\xm=\chi_{Q}|u-v|dx$ and $\xm=\chi_{Q^c}|u-v|dx$ respectively. By uniqueness of the problem (see \cite[Lemma 2.8]{MV-CPDE}), we have
$w=|u-v|^Q_\xb$ and $f=|u-v|^{Q^c}_\xb.$ Let $\xb_k$ be a decreasing sequence such that the following limit exists
$$\lim_{k\rightarrow\infty}|u-v|^{Q^c}_{\xb_{k}}.$$
Since $\lim_{\xb\to 0}|u-v|^Q_\xb=0,$  we have
$$|u-v|\leq\lim_{k\rightarrow\infty}|u-v|^{Q^c}_{\xb_{k}}.$$
Now since $|u-v|^{Q^c}_{\xb_{k}}\approx_{Q}0,$ by Proposition \ref{sygklisi1}(i) we have $\lim_{k\rightarrow\infty}|u-v|^{Q^c}_{\xb_{k}}\approx_Q=0.$
Using the fact  that $u\ominus v$ is the smallest solution which dominates the subsolution $(u-v)_+$, there holds $\max\{u\ominus v,\;v\ominus u\}\leq\lim_{k\rightarrow\infty}|u-v|^{Q^c}_{\xb_{k}}$ and the result follows in this case by Propositions \ref{3} and \ref{7}(iv).\\
(ii) We assume that $u\approx_A v.$\\
For any two positive solutions $u,\;v$ we have
\be
u+(v-u)_+\leq v +(u-v)_+\leq v+u\ominus v\label{4.19}
\ee
If $F$ is a $\frak T_q$-closed set and $Q$ a $\frak T_q$-open set such that $F\qsub Q,$ we claim that
\be
[u]_F\leq [v]_Q+[u\ominus v]_Q. \label{4.20}
\ee
To verify this inequality, we observe first that (see (\ref{13}))
$$u=[u]_{\mathbb{R}^N}\leq[u]_{Q}+[u]_{Q^c},$$
thus by (\ref{4.19})
$$
[u]_F\leq[u]_{\mathbb{R}^N}\leq v+u\ominus v\leq[v]_Q+[v]_{Q^c}+[u\ominus v]_Q+[u\ominus v]_{Q^c}.
$$
The subsolution $w:=\left([u]_F-([v]_Q+[u\ominus v]_Q)\right)_+$ is dominated by the supersolution $[u\ominus v]_{Q^c}+[v]_{Q^c}.$ By definition we have
$$w\leq[w]_\dag\leq[u\ominus v]_{Q^c}\oplus[v]_{Q^c}\leq[u\ominus v]_{Q^c}+[v]_{Q^c}.$$
Thus $[w]_\dag\approx_Q0.$ But $w\leq[u]_F$ which implies $[w]_\dag\leq[u]_F,$ that is $\frak T_q\text{-supp}\,([w]_\dag)\qsub F\qsub Q.$ Taking into account that $[w]_\dag\approx_Q0$ we have that $w=[w]_\dag=0$ and the proof of (\ref{4.20}) is completed.\\
If we choose a $\frak T_q$-open set $Q$ such that $F\qsub Q\subset\widetilde{Q}\qsub A$ (see Lemma \ref{q-open}), and using the fact that $u\ominus v\approx_A=0\Rightarrow[u\ominus v]_F=0$ (see (\ref{4.11})) and (\ref{4.20}), we infer
$$[u]_F\leq[v]_Q.$$
Now by Lemma \ref{6}(I), we can construct a decreasing sequence $\{Q_j\}$ of open sets such that $\cap Q_j\qeq F,$ thus by Proposition \ref{17}(iii) we have
$$[u]_F\leq\lim_{n\to\infty}[v]_{Q_n}=[v]_F.$$
Similarly, $[v]_F\leq[u]_F$ and hence the equality holds.

Next we assume that $[u]_F=[v]_F$ for any $\frak T_q$-closed set $F\qsub A.$ If $Q$ is a $\frak T_q$-open set such that $F\qsub Q\subset\widetilde{Q}\qsub A$ (see Lemma \ref{q-open}), we have
$$u\ominus v\leq\left([u]_Q\oplus[u]_{Q^c}\right)\ominus[v]_Q,$$
where in the last inequality we have used the fact that
$$u=[u]_{\mathbb{R}^N}\leq[u]_{Q}+[u]_{Q^c}\Rightarrow\;u\leq[u]_{Q}\oplus[u]_{Q^c}\leq[u]_{Q}+[u]_{Q^c}.$$
Since $\left([u]_Q\oplus[u]_{Q^c}\right)\ominus([v]_Q)$ is the smallest solution dominating $\left(([u]_Q\oplus[u]_{Q^c})-[v]_Q\right)_+,$ we have
$$\left(([u]_Q\oplus[u]_{Q^c})-[v]_Q\right)_+\leq\left(([u]_Q+[u]_{Q^c})-[v]_Q\right)_+=[u]_Q+[u]_{Q^c}-[v]_Q=[u]_{Q^c},$$
since by assumption we have $[u]_Q=[v]_Q.$ Thus we have
$$[u\ominus v]_F\leq u\ominus v\leq[u]_{Q^c},$$
This means $\frak T_q\text{-supp}\,([u\ominus v]_F)\qsub F$ and $[u\ominus v]_F\approx_{Q}0,$ which in turn implies $[u\ominus v]_F=0,$ and by \ref{4.11} $u\ominus v\approx_A=0.$
Similarly, $v\ominus u\approx_A0.$\hfill$\Box$
\begin{coro}
If $A$ is a $\frak T_q$-open set, the relation $\approx_A$ is an equivalence relation in $\mathcal{U}_+(Q_T).$
\end{coro}
\Proof    This is an immediate consequence of (\ref{4.16}).\hfill$\Box$
\setcounter{equation}{0}

\section{The precise initial trace}
\subsection{The regular initial set}
\begin{lemma}
Let $u\in\mathcal{U}_+(Q_T)$ and $Q$ be a $\frak T_q$-open set. Then for any $\eta\in W^{\frac{2}{q},q'}(\mathbb{R}^N)\cap L^\infty(\mathbb{R}^N)$ with $\frak T_q$-support in $\widetilde{Q}^c,$ we have
$$\int_0^T\int_{\mathbb{R}^N}(u\wedge U_Q)^q(t,x)\mathbb{H}^{2q'}[\eta]_+dxdt<\infty.$$\label{maximal1}
\end{lemma}
\Proof    By  Proposition \ref{anisothta} and the properties of $U_Q$, there holds
$$
\lim_{t\rightarrow0}\int_{Q}u\wedge U_Q(x,t)\eta_+^{2q'}(x)dx=0,
$$
and the result follows by the estimates in Lemma \ref{lem2.9}.\hfill$\Box$
\begin{prop}
Let $u\in\mathcal{U}_+(Q_T)$ and $Q$ be a $\frak T_q$-open set. We assume that $u\wedge U_Q$ is a moderate solution with initial data $\xm.$ Then for any $\xi\in Q$ there exists a $\frak T_q$-open set $O_\xi\subset Q$ such that
$$\int_0^T\int_{\mathbb{R}^N}u^q(t,x)\mathbb{H}^{2q'}[\chi_{O_\xi}]_+dxdt<\infty.$$
Furthermore, for any $\eta\in W^{\frac{2}{q},q'}(\mathbb{R}^N)\cap L^\infty(\mathbb{R}^N)$ with $\frak T_q$-support in $Q,$ we have
$$\lim_{t\rightarrow0}\int_{Q}u(x,t)\eta^{2q'}_+(x)dx=\int_{Q}\eta^{2q'}d\xm.$$\label{maxtrace}
\end{prop}
\Proof
Let $\eta\in W^{\frac{2}{q},q'}(\mathbb{R}^N)\cap L^\infty(\mathbb{R}^N)$ with $\frak T_q$-support in $Q.$ Since $\eta_+^{2q'}$ is a quasi continuous function we have by Lemma \ref{moderate} that
$$\lim_{t\rightarrow0}\int_{Q}u\wedge U_Q(x,t)\eta^{2q'}_+(x)dx=\int_{Q}\eta^{2q'}d\xm.$$
Using the properties of $U_{Q^c}$,
$$\lim_{t\rightarrow0}\int_{Q}u\wedge U_{Q^c}(x,t)\eta^{2q'}_+(x)dx=0.$$
Combining all above and using the fact that $u\leq u\wedge U_{Q}+u\wedge U_{Q^c}$ we get
\bea
\nonumber
\int_{Q}\eta^{2q'}d\xm=\lim_{t\rightarrow0}\int_{Q}u\wedge U_Q(x,t)\eta^{2q'}_+(x)dx
&\leq&\lim_{t\rightarrow0}\int_{Q}u(x,t)\eta^{2q'}_+(x)dx\\ \nonumber
\leq\lim_{t\rightarrow0}\int_{Q}u\wedge U_Q(x,t)\eta^{2q'}_+(x)dx&+&\lim_{t\rightarrow0}\int_{Q}u\wedge U_{Q^c}(x,t)\eta^{2q'}_+(x)dx\\ \nonumber
&=&\int_{Q}\eta^{2q'}d\xm+0.\label{4*}
\eea
In view of the proof of Lemma \ref{lem2.3} and by \ref{lem2.4} there holds
\be
\int_0^T\int_{\mathbb{R}^N}(u\wedge U_{Q})^q(t,x)\mathbb{H}^{2q'}[\eta]_+dxdt<\infty,\label{2*}
\ee
for any $\eta\in W^{\frac{2}{q},q'}(\mathbb{R}^N)\cap L^\infty(\mathbb{R}^N)$ with $\frak T_q$-support in $Q.$ By Lemma \ref{cutoff}, there exists $\eta\in W^{\frac{2}{q},q'}(\mathbb{R}^N)\cap L^\infty(\mathbb{R}^N)$ such that $0\leq\eta\leq1,$ $\eta=1$ on $O_\xi\subset Q$ and $\frak T_q$-supp$\,(\eta)\subset Q.$ Thus we have by (\ref{2*}) and the properties of $\eta$,
\be
\int_0^T\int_{\mathbb{R}^N}(u\wedge U_{Q^c})^q(t,x)\mathbb{H}^{2q'}[\chi_{O_\xi}]dxdt<\infty.\label{3*}
\ee
\hfill$\Box$
\begin{defin} (Section 10.1-\cite{AH})
Let $Q$ be a Borel set. We denote $W^{\frac{2}{q},q'}(E^c)$ the closure of the space of $C^\infty$ functions (with respect the norm $||\cdot||_{W^{\frac{2}{q},q'}}$) with compact support in $E^c$.
\end{defin}
\begin{prop}
Let $u$ be a positive solution of (\ref{maineq}) and $Q$ a bounded $\frak T_q$-open sets such that
\begin{equation}\label{MEA1}
\int_0^T\int_{\mathbb{R}^N}u^q(t,x)\mathbb{H}^{2q'}[\chi_Q]dxdt<\infty.
\end{equation}
(i) Then, there exists an increasing sequence of $\frak T_q$-open set $\{Q_n\}$ satisfying
$Q_n\subset Q,$ $\widetilde{Q}_n\qsub Q_{n+1}$ and $Q_0:=\bigcup_{n=1}^\infty Q_n\qeq Q,$ such that the solution
$v_n=u\wedge Q_n$ is moderate,  $v_n\uparrow [u]_Q$, $\mathrm{tr}(v_n)\rightarrow \xm_Q$.\smallskip

\noindent (ii) For any $\eta\in W^{\frac{2}{q},q'}(Q)$ we have
$$
\lim_{t\rightarrow0}\int_{Q}u(x,t)\eta^{2q'}_+(x)dx=\int_{Q}\eta^{2q'}_+(x)d\xm_Q.
$$\label{measure1}
\end{prop}
\Proof   We choose a point $z\in Q.$ Then by Lemma \ref{cutoff} there exist a $\frak T_q$-open set $V,$ such that $z\in V\subset\widetilde{V}\subset Q,$ and a function $\psi\in W^{\frac{2}{q},q'}(\mathbb{R}^N)$ such that $\psi=1$ q.a.e. on $V$ and $\psi=0$ outside $Q.$ By Lemma \ref{perioxi}, there exists a $\frak T_q$-open set $z\in O_z\subset\widetilde{O}_z\subset V. $

We assert that the function
\be
v_z=u\wedge U_{O_z}\label{reg**}
\ee
is a moderate solution.

Indeed, let $B_R(0)$ be a ball with radius $R$ large enough such that $Q\subset\subset B_R(0).$ Also, let $0\leq\eta\leq1$ be a smooth function with compact support in $B_{2R}(0)$ and $\eta=1$ on $B_R(0).$ Then the function $\zeta=(1-\psi)\eta\in W^{\frac{2}{q},q'}(\mathbb{R}^N)\cap L^\infty(\mathbb{R}^N)$ with compact support in $B_{2R}(0)\setminus \widetilde{V}.$ Now
\bea
\nonumber
\int_0^T\int_{\mathbb{R}^N}v_z^q(t,x)\mathbb{H}^{2q'}[\chi_{B_R(0)}]dxdt&\leq&\int_0^T\int_{\mathbb{R}^N}v_z^q(t,x)\mathbb{H}^{2q'}[\psi]dxdt
+\int_0^T\int_{\mathbb{R}^N}v_z^q(t,x)\mathbb{H}^{2q'}[1-\psi]dxdt\\ \nonumber
&\leq&\int_0^T\int_{\mathbb{R}^N}v_z^q(t,x)\mathbb{H}^{2q'}[\psi]dxdt+\int_0^T\int_{\mathbb{R}^N}v_z^q(t,x)\mathbb{H}^{2q'}[\zeta]dxdt<\infty,
\eea
where the first integral in the last inequality is finite by assumption and the second integral is finite by Lemma \ref{maximal1}. Thus since $B_R(0)$ is arbitrary, the function $u\wedge O_z$ is a moderate solution.

By the {\it quasi-Lindel\"off property} there exists a non decreasing sequence of $\frak T_q$- open set $\{O_n\}$ such that $Q\qeq\cup O_n$ and (by the above arguments) the solution $u\wedge U_{O_n}$ is moderate for any $n\in\mathbb{N}.$ Now, by Lemma \ref{6} (II)(i)-(ii), for any $n\in\mathbb{N},$ there exists an increasing sequence $\{A_{n,j}\}$ of $\frak T_q$-open sets such that $\widetilde{A}_{n,j}\qsub A_{n,j+1}\qsub E_n$ and $\bigcup_{j=1}^\infty A_{n,j}\qeq E_n.$ Put
$$Q_n=\bigcup_{k+j=n}A_{k,j}.$$
Then
$$\widetilde{Q}_n\subset\bigcup_{k+j=n}\widetilde{A}_{k,j}\qsub\bigcup_{k+j=n}\widetilde{A}_{k,j+1}=Q_{n+1}.$$
Hence, $$Q_0:=\bigcup Q_n\qeq Q.$$

Now, we will prove that $v_n=u\wedge U_{Q_n}\rightarrow u\wedge U_{Q}.$ By Proposition \ref{17}(ii) we have
$$ u\wedge U_{Q}\leq u\wedge U_{Q_n}+u\wedge U_{Q\setminus Q_n}.$$
Since $Q\setminus Q_n\downarrow F$ with $\qcap(F)=0,$ we have by Proposition \ref{17}(iii) that
$$u\wedge U_{Q\setminus Q_n}\rightarrow0.$$
The opposite inequality is obvious and the result follows in this assertion.
By Lemma \ref{17}(ii) $v_n=[v_{n+k}]_{Q_n},\;\forall k\in\mathbb{N}.$ Therefore
\be
\xm_n(Q_n)=\xm_{n+k}(Q_n)=\xm_Q(Q_n).\label{defmQ}
\ee
(ii) First we assume that the function $\eta\in W^{\frac{2}{q},q'}(Q)$ has compact support in $Q.$ Then by Lemma \ref{partition} there exists a function $\eta_k$ such that $\frak T_q$-supp $\,(\eta_k)\subset Q_k$, and
\be
||\eta-\eta_k||_{ W^{\frac{2}{q},q'}}\leq\frac{1}{k},\label{est1}
\ee
and $|\eta_k|\leq|\eta|.$
By Lebesgue's dominated theorem, we can assume that $\eta_k$ satisfies
$$\int_0^T\int_{\mathbb{R}^N}u^q(t,x)(\mathbb{H}[\eta-\eta_k])^{2q'}dxdt<\frac{1}{k}$$
Also in view of Proposition \ref{menest} and (\ref{2.22})-(\ref{2.26}),
$$
\lim_{t\rightarrow0}\int_{Q}u(x,t)\eta^{2q'}(x)dx\leq C||\eta||_{L^\infty(\mathbb{R}^N)}^{q'}||\eta||_{ W^{\frac{2}{q},q'}}^{q'}+\int_0^T\int_{\mathbb{R}^N}u^q(t,x)(\mathbb{H}[\eta])^{2q'}dxdt,
$$
But by (\ref{defmQ}) and Lemma \ref{maxtrace} we have
\bea
\nonumber
\left(\int_{Q}\eta_k^{2q'}(x)d\xm_Q\right)^{\frac{1}{2q'}}&=&\lim_{t\rightarrow0}\left(\int_{Q}u(x,t)\eta_k^{2q'}(x)dx\right)^{\frac{1}{2q'}}
\\ \nonumber
&\leq&\lim_{t\rightarrow0}\left(\int_{Q}u(x,t)\eta^{2q'}(x)dx\right)^{\frac{1}{2q'}}\\ \nonumber
&\leq&\lim_{t\rightarrow0}\left(\int_{Q}u(x,t)\left(\eta-\eta_k\right)^{2q'}(x)dx\right)^{\frac{1}{2q'}}
+\lim_{t\rightarrow0}\left(\int_{Q}u(x,t)\eta_k^{2q'}(x)dx\right)^{\frac{1}{2q'}}\\ \nonumber
&\leq&\left(\int_{Q}\eta_k^{2q'}(x)d\xm_Q\right)^{\frac{1}{2q'}}+C||\eta-\eta_k||_{L^\infty(\mathbb{R}^N)}^{\frac{1}{2}}||\eta-\eta_k||_{ W^{\frac{2}{q},q'}}^{\frac{1}{2}}\\ \nonumber
&+&\left(\int_0^T\int_{\mathbb{R}^N}u^q(t,x)(\mathbb{H}[\eta-\eta_k])^{2q'}dxdt\right)^\frac{1}{2q'}\\ \nonumber
&\leq&\left(\int_{Q}\eta_k^{2q'}(x)d\xm_Q\right)^{\frac{1}{2q'}}+C\frac{1}{\sqrt{k}}||\eta||_{L^\infty(\mathbb{R}^N)}^{\frac{1}{2}}
+\left(\frac{1}{k}\right)^\frac{1}{2q'}.
\eea
The result follows in this case by letting $k\to\infty$.\\
For the general case, by theorem 10.1.1 in \cite{AH}, there exists a function $\eta_k$ with compact support in $Q$ such that
\be
||\eta-\eta_k||_{ W^{\frac{2}{q},q'}}\leq\frac{1}{k},\label{est2}
\ee
and $|\eta_k|\leq|\eta|.$ The result follows as above.\hfill$\Box$\\

\noindent\Remark By Lemma \ref{maxtrace} and (\ref{reg**}), we have that the definition of the regular points in the elliptic case (see \cite{MV-CONT}) coincides with our definition of the regular points.
%
\begin{lemma}
Let $Q$ be a $\frak T_q$-open set and $u\in \CU_+(Q_T)$ satisfy (\ref{MEA1}). Then\\
i)
\be
[u]_Q=\sup\{[u]_F:\;F\qsub Q,\;F\;\,\frak T_q\mathrm{-closed}\}.\label{5.25*i}
\ee
ii) For every $\frak T_q$-open set $O\subset\widetilde{O}\qsub Q$ such that $[u]_O$ is a moderate solution we have
\be
\xm_Q\chi_{\widetilde{O}}=\mathrm{tr}'[u]_O)=\mathrm{tr}([u_Q]_O).\label{5.25*}
\ee
Finally, $\xm_Q$ is $\frak T_q$-locally finite on $Q$ and $\xs$-finite on $Q':=\cup Q_n.$\\
iii) If $\{w_n\}\subset\mathcal{U}_+(Q_T)$ is a nondecreasing sequence of moderate solutions such that $\frak T_q\text{-supp}(w_n)\qsub Q$ and $w_n\uparrow [u]_Q$, then $\mathrm{tr}(w_n)=\xn_n\uparrow\xm_Q.$\label{measure1*}\\
\end{lemma}
\Proof    i) Let $u^*$ denote the right-hand side of (\ref{5.25*i}). By Proposition \ref{sygklish} there exists a nondecreasing sequence $\{[u]_{F_n}\}$ such that $[u]_{F_n}\uparrow u^*.$ We consider the function $[u]_{Q_n}$ of Proposition  \ref{measure1}. Then by Proposition \ref{17} we have
$$[u]_{F_n}\leq [u]_{F_n\cap Q_m}+[u]_{F_n\setminus Q_m}.$$
Now we note that $F_n\setminus Q_m$ is a $\frak T_q$-closed set and $\cap_{m=1}^\infty F_n\setminus Q_m=A$ with $\qcap (A)=0.$ Thus by Proposition \ref{7} we have that $\lim_{m\rightarrow\infty}U_{F_n\setminus Q_m}=0$ which implies $\lim_{m\rightarrow\infty}[u]_{F_n\setminus Q_m}=0.$ Thus $[u]_{F_n}\leq\lim[u]_{Q_m}=u_Q.$ Letting $n\rightarrow\infty$ we have $u^*\leq u_Q.$ By definition of $u^*$ we have that $u_Q\leq u^*,$ thus $u^*=u_Q.$\\\\
ii) Put $\xm_O=\mathrm{tr}([u]_O).$ If $F$ is a $\frak T_q$-closed set such that $F\qsub O$ , by Proposition \ref{17}-(ii) we have
\be
\mathrm{tr}([u]_F)=\mathrm{tr}([[u]_O]_F)=\xm_O\chi_F.\label{5.35*}
\ee
In particular the compatibility condition holds: if $O'\subset\widetilde{O}'\qsub Q$ is $\frak T_q$-open set such that $[u]_{O'}$ is a moderate solution
\be
\xm_{O\cap O'}=\xm_O\chi_{\widetilde{O}\cap\widetilde{O}'}=\xm_{O'}\chi_{\widetilde{O}\cap\widetilde{O}'}.\label{5.36*}
\ee
With the notation of (\ref{defmQ}), $[v_{n+k}]_{Q_k}=v_k$ and hence $\xm_{n+k}\chi_{\widetilde{Q}_k}=\xm_k$ for every $k\in\mathbb{N}.$

Since $[u]_F$ is moderate, we have by (\ref{5.36*})
\be
[v_n]_F=[u]_{F\cap\widetilde{Q}_n}\uparrow[u]_F.\label{5.37*}
\ee
In addition, $[u_Q]_F\geq\lim_{n\to\infty}[v_n]_F=[u]_F$, jointly with $u_Q\leq u$, leads to,
\be
[u]_F=[u_Q]_F.\label{5.38*}
\ee

By (\ref{5.35*}) and (\ref{5.37*}), if $F$ is a $\frak T_q$-closed subset of $\mathcal{R}_q(u)$ and $[u]_F$ is moderate,
\be
\mathrm{tr}([u]_F)=\lim_{n\to\infty}\mathrm{tr}([v_n]_F)=\lim_{n\to\infty}\xm_n\chi_F=\xm_{\mathcal{R}_q}\chi_F,\label{5.39*}
\ee
which implies (\ref{5.25*}).

Since $Q':=\cup Q_n\qeq Q,$ $\xm_Q$ is $\xs$-finite on $Q'.$ The assertion that $\xm_Q$ is $\frak T_q$-locally finite on $Q$ is a consequence of the fact that every point in $Q$ is contained in a $\frak T_q$-open set $O\qsub\widetilde{O}\subset Q$ such that $[u]_{O}$ is a moderate solution (see (\ref{reg**})).\\\\
iii) If $w$ is a moderate solution and $w\leq u_Q$ and $\frak T_q\text{-supp}\,(w)\qsub Q$, then $\tau:=\mathrm{tr}(w)\leq \xm_Q.$ Indeed
$$[w]_{Q_n}\leq[u]_{Q_n}=v_n,\;[w]_{Q_n}\uparrow w\Rightarrow \mathrm{tr}([w]_{Q_n})\uparrow\tau\leq\lim_{n\to\infty}\mathrm{tr}(v_n)=\xm_Q.$$

Now, let $\{w_n\}$ be an increasing sequence of moderate solutions such that $F_n:=\frak T_q\text{-supp}\,(w_n)\qsub Q$ and $w_n\uparrow u_Q.$ We must show that, if $\xn_n:=\mathrm{tr}(w_n),$ then
\be
\xn:=\lim_{n\to\infty} \xn_n=\xm_Q.\label{5.40*}
\ee
By the previous argument $\xn\leq\xm_Q.$ The opposite inequality is obtained as follows. Let $D$ be a $\frak T_q$-open  set such that $[u]_D$ is moderate. Also, let $K$ be a compact subset of $D$ such that $\qcap(K)>0.$
$$w_n\leq[w_n]_D+[w_n]_{D^c}\Rightarrow u_Q=\lim_{n\to\infty} w_n\leq\lim_{n\to\infty}[w_n]_D+U_{D^c}.$$
The sequence $\{[w_n]_D\}$ is dominated by the moderate solution $[u_Q]_D.$ In addition $\mathrm{tr}([w_n]_D)=\xn_n\chi_{\widetilde{D}}\uparrow\xn\chi_{\widetilde{D}}.$ Hence, $\xn\chi_{\widetilde{D}}$ is a Radon measure which vanishes on sets with $\qcap$-capacity zero. Also, $[w_n]_D\uparrow u_{\xn\chi_{\widetilde{D}}},$ where $u_{\xn\chi_{\widetilde{D}}}$ is a moderate solution with initial trace $\xn\chi_{\widetilde{D}}.$ Consequently
$$u_Q=\lim_{n\to\infty} w_n\leq u_{\xn\chi_{\widetilde{D}}}+U_{D^c}.$$
This in turn implies
$$\left([u_Q]_K-u_{\xn\chi_{\widetilde{D}}}\right)_+\leq\inf(U_{D^c},U_K),$$
the function on the left being a subsolution and the one on the right a supersolution. Therefore
$$\left([u_Q]_K-u_{\xn\chi_{\widetilde{D}}}\right)_+\leq[[U]_{D^c}]_K=0.$$
Thus, $[u_Q]_K\leq u_{\xn\chi_{\widetilde{D}}}$ and hence $\xm_Q\chi_K\leq\xn\chi_{\widetilde{D}}.$ Further, if $O$ is a $\frak T_q$-open set such that $\widetilde{O}\qsub D$ then, in view of the fact that
$$\sup\{\xm_Q\chi_K:\;K\subset O,\;K\;\mathrm{compact}\}=\xm_Q\chi_O,$$
we obtain,
\be
\xm_{Q}\chi_O\leq\xn\chi_{\widetilde{D}}.\label{5.41*}
\ee
Applying this inequality to the sets $Q_m,\;Q_{m+1}$ we finally obtain
$$\xm_{Q}\chi_{Q_{m}}\leq\xn\chi_{\widetilde{Q}_{m+1}}\leq\xn\chi_{{Q}_{m+2}}.$$
Letting $m\rightarrow\infty$ we conclude that $\xm_{\mathcal{R}_q}\leq\xn.$ This completes the proof of (\ref{5.40*}).\hfill$\Box$
\subsection{$\frak T_q$-perfect measures}
\begin{defin}
Let $\xm$ be a positive Borel measure on $\mathbb{R}^N.$\\
(i) We say that $\xm$ is \textbf{essentially absolutely continuous relative to $\qcap$} if the following condition holds:

If $Q$ is a $\frak T_q$-open set and $A$ is a Borel set such that $\qcap(A)=0$ then
$$\xm(Q\setminus A)=\xm(Q).$$
This relation be denoted by $$\xm\prec\prec_f\qcap.$$
(ii) $\xm$ is \textbf{regular relative to $\frak T_q$-topology} if, for every Borel set $E,$
\bea
\nonumber
\xm(E)&=&\inf\{\xm(D):\;E\subset D,\;D\;\frak T_q\mathrm{-open}\}\\
      &=& \inf\{\xm(K):\;K\subset E,\;K\;\mathrm{compact}\}.\label{5.13}
\eea
$\xm$ is \textbf{outer regular relative to $\frak T_q$-topology} if the first equality in (\ref{5.13}) holds.\\
(iii) A positive Borel measure is called \textbf{$\frak T_q$-perfect} if it is essentially absolutely continuous relative to $\qcap$ and outer regular relative to $\frak T_q$-topology. The space of $\frak T_q$-perfect Borel measures is denoted by $\mathfrak{M}_q(\mathbb{R}^N).$
\end{defin}
\begin{prop}
If $\xm\in\mathfrak{M}_q(\mathbb{R}^N)$ and $A$ is a non-empty Borel set such that $\qcap(A)=0$, then
\be
\xm=\Bigg \{
\begin{array}{ll}\infty\;\;\;\;\mathrm{if} \;\xm(Q\setminus A)=\infty\;\;\;\;\forall Q\;\text{$\frak T_q$-open neighborhood of A},
\\[2pt]
0\;\;\;\; \mathrm{otherwise}.
\end{array}\label{5.14}
\ee
If $\xm_0$ is an essentially absolutely continuous positive measure on $\mathbb{R}^N$ and $Q$ is $\frak T_q$-open set such that $\xm_0(Q)<\infty$ then $\xm_0|_Q$ is \textbf{absolutely continuous with respect to $\qcap$ in the strong sense}, i.e., if $\{A_n\}$ is a sequence of Borel subsets of $\mathbb{R}^N$
$$\qcap(A_n)\rightarrow0\Rightarrow\xm_0(Q\cap A_n)\rightarrow0.$$
Let $\xm_0$ is an essentially absolutely continuous positive Borel measure on $\mathbb{R}^N$ and denote
\be
\xm(E)=\inf\{\xm_0(D):\;E\subset D,\;D\;\text{$\frak T_q$-open}\},\label{5.15}
\ee
for every Borel set $E$; then
\bea
\nonumber
&(a)\;\;\;\;\;\;\xm_0\leq\xm&\qquad \xm_0(Q)=\xm(Q)\qquad\forall Q\;\;\text{$\frak T_q$-open},\\
&(b)\;\;\;\;\xm|_Q=\xm_0|_Q&\text{ for every $\frak T_q$-open set } Q\;\text{ such that }\xm_0(Q)<\infty.\label{5.16}
\eea
Finally $\xm$ is $\frak T_q$-perfect; thus $\xm$ is the smallest measure in $\mathfrak{M}_q(\mathbb{R}^N)$ which dominates $\xm_0.$\label{21}
\end{prop}
\Proof    The first assertion follows immediately from the definition $\mathfrak{M}_q(\mathbb{R}^N).$ We turn to the second assertion. If $\xm_0$ is an essentially absolutely continuous positive Borel measure on $\mathbb{R}^N,$ and $Q$ is a $\frak T_q$-open set such that $\xm_{0}(Q)<\infty$ then $\xm_0\chi_Q$ is a bounded Borel measure which vanishes on sets of $\qcap-$ capacity zero. If $\{A_n\}$ is a sequence of Borel sets such that $\qcap(A_n)\rightarrow0$ and $\xm_n=\chi_{Q\cap A_n},$ then by Lemma 2.8-\cite{MV-CPDE}, there exists a unique moderate solution $u_{\xm_n}.$ Also in view of Lemma 2.8-\cite{MV-CPDE} we can prove that the sequence $\{u_{\xm_n}\}$ is decreasing. Also by Proposition \ref{sygklisi1}, we have $u_{\xm_n}\leq U_{Q\cap A_n}\rightarrow0,$ since $\qcap(Q\cap A_n)\rightarrow0.$ Thus we have that $u_{\xm_n}\rightarrow0$ locally uniformly and $\xm_n\rightharpoonup0$ weakly with respect to $C_0(\mathbb{R}^N).$ Hence $\xm(Q\cap A_n)\rightarrow0.$ Thus $\xm_0|_Q$ is absolutely continuous with respect to $\qcap$ in the strong sense.

Assertion (\ref{5.16})(a) follows from (\ref{5.15}). It is clear that $\xm$, as defined by (\ref{5.15}), is a measure. Now if $Q$ is $\frak T_q$-open set such that $\xm_0(Q)<\infty$, then $\xm(Q)<\infty$ and both $\xm_0|_Q$ and $\xm|_Q$ are regular. Since they agree on open sets, the regularity implies (\ref{5.16}) (b).

If $A$ is a Borel set such that $\qcap(A)=0$ and $Q$ is a $\frak T_q$-open set then $Q\setminus A$ is $\frak T_q$-open and consequently
$$\xm(Q)=\xm_0(Q)=\xm_0(Q\setminus A)=\xm(Q\setminus A).$$
Thus $\xm$ is essentially absolutely continuous. By (\ref{5.16}) (a) and the definition of $\xm,$ we have that $\xm$ is outer regular with respect to $\qcap.$ Thus $\xm\in\mathfrak{M}_q(\mathbb{R}^N).$\hfill$\Box$

\subsection{The initial trace on the regular set}
\begin{prop}
Let $u\in\mathcal{U}_+(Q_T).$\\
(i) There exists an increasing sequence of $\frak T_q$-open sets $\{Q_n\}$ with the properties
$Q_n\subset \mathcal{R}_q(u),$ $\widetilde{Q}_n\qsub Q_{n+1}$ and $\mathcal{R}_{q,0}(u):=\bigcup_{n=1}^\infty Q_n\qeq \mathcal{R}(u),$ such that the solution
\be
v_n=u\wedge U_{Q_n}\;\mathrm{is\;moderate}\qquad v_n\uparrow v_{\mathcal{R}_q},\qquad \mathrm{tr}(v_n)\rightarrow \xm_{\mathcal{R}_q}. \label{5.22}
\ee
(ii)
\be
v_{\mathcal{R}_q}:=\sup\{[u]_F:\;F\qsub\mathcal{R}_q(u),\;F\;\frak T_q\mathrm{-closed}\}.\label{5.21}
\ee
Thus $v_{\mathcal{R}_q}$ is $\xs$-moderate.\\
(iii) If $[u]_F$ is moderate and $F\qsub\mathcal{R}_q(u),$  there exists a $\frak T_q$-open set $Q$ such that $F\qsub Q,$ $[u]_Q$ is moderate solution and $Q\subset\mathcal{R}_q(u)$\\
(iv) For every $\frak T_q$-open set $Q,$ such that $[u]_Q$ is a moderate solution, we have
\be
\xm_{\mathcal{R}_q}\chi_{\widetilde{Q}}=\mathrm{tr}([u]_Q)=\mathrm{tr}([v_{\mathcal{R}_q}]_Q).\label{5.25}
\ee
Finally, $\xm_{\mathcal{R}_q}$ is $\frak T_q$-locally finite on $\mathcal{R}_q(u)$ and $\xs$-finite on $\mathcal{R}_{q,0}(u):=\cup Q_n.$\\
(v) If $\{w_n\}$ is a sequence of moderate solutions such that $w_n\uparrow u_{\mathcal{R}_q}$ then,
\be
\xm_{\mathcal{R}_q}=\lim_{n\to\infty}\mathrm{tr}(w_n):=\lim_{n\to\infty}\xn_n. \label{5.26}
\ee
(vi) The regularized measure $\overline{\xm}_{\mathcal{R}_q}$ given by
\be
\overline{\xm}_{\mathcal{R}_q}(E)=\inf\{\xm_{\mathcal{R}_q}(Q):\;E\subset Q,\quad Q\;\frak T_q\text{-open},\quad E\;\text{ Borel}\},\label{5.27}
\ee
is $\frak T_q$-perfect.\\
(vii) $$u\approx_{\mathcal{R}_q(u)}v_{\mathcal{R}_q}.$$
(viii) For every $\frak T_q$-closed set $F\qsub\mathcal{R}_q(u):$
\be
[u]_F=[v_{\mathcal{R}_q}]_F.\label{5.28}
\ee
If, in addition, $\xm_{\mathcal{R}_q}(F\cap K)<\infty$ for any compact $K\subset\mathbb{R}^N,$ then $[u]_F$ is moderate and
\be
\mathrm{tr}([u]_F)=\xm_{\mathcal{R}_q}\chi_F.\label{5.29}
\ee
(ix) If $F$ is a $\frak T_q$-closed set and $\qcap(F)>0$ then
\be
\xm_{\mathcal{R}_q}(F\cap K)<\infty\qquad \mathrm{for \;any \;compact}\; K\subset\mathbb{R}^N\Leftrightarrow [u]_F\;\;\mathrm{is\;moderate}.\label{5.30}
\ee\label{24}
\end{prop}
\Proof
(i) For any $z\in \mathcal{R}_q(u)$ there exist a bounded $\frak T_q$-open set $Q\subset \mathcal{R}_q(u)$ such that
$$\int_0^T\int_{\mathbb{R}^N}u^q(t,x)\mathbb{H}^{2q'}[\chi_Q]dxdt<\infty.$$

The result follows by similar arguments as in Lemma \ref{measure1}.
Also, we recall that for any $z\in\mathcal{R}_q(u)$ there exists a $\frak T_q$-open set $O_z\subset\mathcal{R}_q(u)$ such that
\be
[u]_{O_z},\label{moderate1}
\ee
is moderate.

Also we recall that
$v_n=[v_{n+k}]_{Q_n},\;\forall k\in\mathbb{N}$ and
\be
\xm_n(Q_n)=\xm_{n+k}(Q_n)=\xm_{\mathcal{R}_q}(Q_n).\label{defmQ*}
\ee
(ii) The proof is same as the one of Lemma \ref{measure1*}-a)\\
(iii) First we assume that $F$ is bounded. By definition and (\ref{moderate1}), every point in $\mathcal{R}_q(u)$ possesses a $\frak T_q$-open neighborhood $A$ such that $[u]_A$ is moderate. Then by Proposition \ref{cover}, for any $\xe>0$ there exists a $\frak T_q$-open set $Q_\xe$ such that $\qcap(F\setminus Q_\xe)<\xe$ and $[u]_{Q_\xe}$ is moderate. Since $F$ is bounded, we can assume that so is $Q_\xe$. Let $O_\xe$ be an open set containing $F\setminus Q_\xe$ such that $\qcap(O_\xe)<2\xe.$ Put
\be
F_\xe:=F\setminus O_\xe.\label{22}
\ee
Then $F_\xe$ is a $\frak T_q$-closed set, $F_\xe\subset F,$ $\qcap(F\setminus F_\xe)<2\xe$ and $F_\xe\subset Q_\xe$.\\
\emph{Assertion 1. Let $E$ be a $\frak T_q$-closed set, $D$ a $\frak T_q$-open set such that $[u]_D$ is moderate and $E\qsub D.$ Then there exists a decreasing sequence of $\frak T_q$-open sets
$\{G_n\}$ such that
\be
E\qsub G_{n+1}\subset\widetilde{G}_{n+1}\qsub G_n\qsub D,\label{5.31}
\ee
and
\be
[u]_{G_n}\rightarrow[u]_E\qquad \mathrm{in}\; L^q(K)\; \mathrm{for\;any\;compact}\; K\subset\overline{Q}_T.\label{5.32}
\ee
}

By Lemma \ref{6} and Proposition \ref{17}-(iii), there exists a decreasing sequence of $\frak T_q$-open sets $\{G_n\}$ satisfying (\ref{5.31}) and, in addition, such that $[u]_{G_n}\downarrow[u]_E$ locally uniformly in $\mathbb{R}^N.$ Since $[u]_{G_n}\leq[u]_D$ and the later is a moderate solution we obtain (\ref{5.32}).

Now we assume that $F$ is $\frak T_q$-closed set (possibly unbounded). Let $x\in F$ and $B_n=B_n(x);\;n\in\mathbb{N}.$ Set
$$E_n=\bigcup_{m=1}^n(F\cap B_n)_{\frac{1}{2^m}},$$
where $(F\cap B_n)_{\frac{1}{2^m}}$ is the set  in (\ref{22}), if we replace $F$ by $F\cap B_n$ and $\xe$ by $\frac{1}{2^m}.$ Also we assume without loss of generality that $\{E_n\}$ is an increasing sequence. Also set
$$Q_n=\bigcup_{m=1}^nQ_{\frac{1}{m}}^n,$$
where $Q^n_{\frac{1}{m}}=(F\cap B_n)_{\frac{1}{m}}.$ Also we may assume that the sequence of set $\{Q_n\}$ is increasing.
Therefore, we have that $E_n\subset E,$ $Q_n$ is $\frak T_q$-open, $[u]_{Q_n}$ is moderate and $E_n\qsub Q_n$ and $\cup E_n=E'\qeq F,$ since
\bea
\nonumber
\qcap(F\setminus\bigcup_{n=1}^\infty E_n)\leq \sum_{k=1}^n\qcap\left((F\cap B_k)\setminus\bigcup_{j=1}^\infty E_j\right)
&+&\sum_{k=n+1}^\infty\qcap\left((F\cap B_k)\setminus E_k\right)\\ \nonumber
&\leq&\frac{1}{2^n}+\sum_{k=n+1}^\infty\frac{1}{2^k},\;\;\forall n\in\mathbb{N}.
\eea
Thus by Assertion 1, it is possible to choose a sequence of $\frak T_q$-open sets $\{V_n\}$ such that
\be
E_n\qsub V_n\subset\widetilde{V}_n\qsub Q_n,\qquad ||[u]_{V_n}-[u]_{E_n}||_{L^q(B_n(0))\times(0,T]}\leq2^{-n}.\label{5.33}
\ee

We note here that since $E_n,\; Q_n$ are bounded sets, the function $[u]_{V_n},\;[u]_{E_n}$ have compact support with respect to variable "$x$" in $\mathbb{R}^N,$ thus we can take the norm in (\ref{5.33}) in whole space $\mathbb{R}^N\times(0,T].$

Because $[u]_F$ is moderate, there exists a Radon measure $\xm_F$ where $\xm_F=\mathrm{tr}([u]_F).$ Moreover, $[u]_F=[u]_{E'}$ since  $F\qeq E'$.  Finally, we have by (\ref{14}) and the fact that $E_n\qsub F$,
$$[u]_{E_n}=[u]_{E_n\cap F}=[[u]_{E_n}]_{F}.$$
Using the above equality and the fact that $[u]_{F}$ is moderate, we have that $\mathrm{tr}([u]_{E_n})=\chi_{E_n}\xm_F.$ Now since $E_n\uparrow E'\qeq F$, it implies that
$[u]_{E_n}\uparrow[u]_{F}$ $L^q(K\times[0,T])$, for compact set $K\subset\mathbb{R}^N.$
Hence, we derive from by (\ref{5.33}) that $[u]_{V_n}\rightarrow[u]_{F}$ in $L^q(K\times[0,T])$ for each compact set $K\subset\mathbb{R}^N$.\\
Let $\{V_{n_k}\}$ be a sequence such that
\be
\left(\int_0^1\int_{B_k(0)}|[u]_{V_{n_k}}-[u]_{F}|^qdxdt\right)^{\frac{1}{q}}\leq 2^{-k}.\label{5.34}
\ee
If $K$is a compact set, there exist $k_0\in\mathbb{N}$ such that $K\subset B_{k}(0),\;\forall k\geq k_0.$ Set $W=\bigcup_{k=1}^\infty V_{n_k}$, then
$$[u]_W\leq\sum_{k=1}^\infty[u]_{V_{n_k}}.$$
Thus we have
\bea
\nonumber
\left(\int_0^T\int_{K}|[u]_{W}-[u]_{F}|^qdxdt\right)^{\frac{1}{q}}&\leq&
\sum_{k=1}^{k_0}\left(\int_0^T\int_{K}|[u]_{V_{n_k}}-[u]_{F}|^qdxdt\right)\phantom{--------------}\\ \nonumber
&\phantom{-}&\phantom{------}+\sum_{k=k_0+1}^{\infty}\left(\int_0^T\int_{B_k(0)}|[u]_{V_{n_k}}-[u]_{F}|^qdxdt\right)^{\frac{1}{q}}\\ \nonumber
&\leq&\sum_{k=1}^{k_0}\left(\int_0^T\int_{K}|[u]_{V_{n_k}}-[u]_{F}|^qdxdt\right)
+\sum_{k=k_0+1}^\infty2^{-k}\phantom{--------------}\\
&<&\infty.
\eea
We recall that $F\qsub W$ and $W$ is a $\frak T_q$-open set. Using the facts that $[u]_F$ is moderate, $K$ is an abstract compact domain and the above inequality, we obtain that $[u]_W$ is moderate. Thus by Lemma \ref{maxtrace} we have that $W\subset \mathcal{R}_q(u).$\\
(iv) Let $Q$ be a $\frak T_q$-open set and $[u]_Q$ be a moderate solution, put $\xm_Q=\mathrm{tr}([u]_Q).$ If $F$ is a $\frak T_q$-closed set such that $F\qsub Q$ then, by Proposition \ref{17}-(ii),
\be
\mathrm{tr}[u]_F=\mathrm{tr}([[u]_Q]_F)=\xm_Q\chi_F.\label{5.35}
\ee
In particular the compatibility condition holds: if $Q,\;Q'$ are $\frak T_q$-open regular sets then
\be
\xm_{Q\cap Q'}=\xm_Q\chi_{\widetilde{Q}\cap\widetilde{Q}'}=\xm_{Q'}\chi_{\widetilde{Q}\cap\widetilde{Q}'}.\label{5.36}
\ee
With the notation of (i), $[v_{n+k}]_{Q_k}=v_k$ and hence $\xm_{n+k}\chi_{\widetilde{Q}_k}=\xm_k$ for every $k\in\mathbb{N}.$

Let $F$ be an arbitrary $\frak T_q$-closed regular subset of $\mathcal{R}_q(u).$ Since $[u]_F$ is moderate, we have by (\ref{5.36})
\be
[v_n]_F=[u]_{F\cap\widetilde{Q}_n}\uparrow[u]_F.\label{5.37}
\ee
In addition, $[v_{\mathcal{R}_q}]_F\geq\lim_{n\to\infty}[v_n]_F=[u]_F$, jointly with $v_{\mathcal{R}_q}\leq u$, leads to,
\be
[u]_F=[v_{\mathcal{R}_q}]_F.\label{5.38}
\ee

By (\ref{5.35}) and (\ref{5.37}), if $F$ is a $\frak T_q$-closed subset of $\mathcal{R}_q(u)$ and $[u]_F$ is moderate,
\be
\mathrm{tr}([u]_F)=\lim_{n\to\infty}\mathrm{tr}([v_n]_F)=\lim_{n\to\infty}\xm_n\chi_F=\xm_{\mathcal{R}_q}\chi_F,\label{5.39}
\ee
which implies (\ref{5.25}).

Since $\mathcal{R}_{q,0}(u)$ has a regular decomposition, $\xm_{\mathcal{R}_q}$ is $\xs$-finite on $\mathcal{R}_{q,0}(u).$ The assertion that $\xm_{\mathcal{R}_q}$ is $\frak T_q$-locally finite on $\mathcal{R}_q(u)$ is a consequence of the fact that every point $\xi\in\mathcal{R}_q(u)$ is contained in a $\frak T_q$-open set $O_\xi\subset\mathcal{R}_q(u)$ such that $[u]_{O_\xi}$ is moderate and thus $\xm_{\mathcal{R}_q}\chi_{O_\xi}<\infty$ .\\
(v) If $w$ is a moderate solution and $w\leq v_{\mathcal{R}_q}$ and $\frak T_q\text{-supp}\,(w)\qsub{\mathcal{R}_q}(u)$ then $\tau:=\mathrm{tr}(w)\leq \xm_{\mathcal{R}_q}.$ Indeed
$$[w]_{Q_n}\leq[v_{{\mathcal{R}_q}}]=v_n,\;[w]_{Q_n}\uparrow w\Rightarrow \mathrm{tr}([w]_{Q_n})\uparrow\tau\leq\lim_{n\to\infty}\mathrm{tr}(v_n)=\xm_{\mathcal{R}_q}.$$

Now, let $\{w_n\}$ be an increasing sequence of moderate solutions such that $F_n:=\frak T_q\text{-supp}\,(w_n)\qsub\mathcal{R}_q(u)$ and $w_n\uparrow v_{\mathcal{R}_q}.$ If $\xn_n:=\mathrm{tr}(w_n),$ we have to prove that
\be
\xn:=\lim_{n\to\infty} \xn_n=\xm_{\mathcal{R}_q}.\label{5.40}
\ee
By the previous argument $\xn\leq\xm_{\mathcal{R}_q}.$ The opposite inequality is obtained as follows. Let $D$ be a $\frak T_q$-open set, $[u]_D$ be moderate and let $K$ be a compact subset of $D$ such that $\qcap(K)>0.$
$$w_n\leq[w_n]_D+[w_n]_{D^c}\Rightarrow v_{\mathcal{R}_q}=\lim_{n\to\infty} w_n\leq\lim_{n\to\infty}[w_n]_D+U_{D^c}.$$
The sequence $\{[w_n]_D\}$ is dominated by the moderate function $[v_{\mathcal{R}_q}]_D.$ In addition $\mathrm{tr}([w_n]_D)=\xn_n\chi_{\widetilde{D}}\uparrow\xn\chi_{\widetilde{D}}.$ Hence, $\xn\chi_{\widetilde{D}}$ is a Radon measure which vanishes on sets with $\qcap$-capacity zero. Also, $[w_n]_D\uparrow u_{\xn\chi_{\widetilde{D}}},$ where the function $u_{\xn\chi_{\widetilde{D}}}$ on the right is the moderate solution with initial trace $\xn\chi_{\widetilde{D}}.$ Consequently
$$v_{\mathcal{R}_q}=\lim_{n\to\infty} w_n\leq u_{\xn\chi_{\widetilde{D}}}+U_{D^c}.$$
This in turn implies
$$\left([v_{\mathcal{R}_q}]_K-u_{\xn\chi_{\widetilde{D}}}\right)_+\leq\inf(U_{D^c},U_K).$$
Note that, in the previous relation, the function on the left being a subsolution and the one on the right a supersolution, we obtain
$$\left([v_{\mathcal{R}_q}]_K-u_{\xn\chi_{\widetilde{D}}}\right)_+\leq[[U]_{D^c}]_K=0.$$
Thus, $[v_{\mathcal{R}_q}]_K\leq u_{\xn\chi_{\widetilde{D}}}$ and hence $\xm_R\chi_K\leq\xn\chi_{\widetilde{D}}.$ Further, if $Q$ is a $\frak T_q$-open set such that $\widetilde{Q}\qsub D$ then, in view of the fact that
$$\sup\{\xm_{\mathcal{R}_q}\chi_K:\;K\subset Q,\;K\;\mathrm{compact}\}=\xm_{\mathcal{R}_q}\chi_Q,$$
we obtain,
\be
\xm_{\mathcal{R}_q}\chi_Q\leq\xn\chi_{\widetilde{D}}.\label{5.41}
\ee
Applying this inequality to the sets $Q_m,\;Q_{m+1}$ we finally obtain
$$\xm_{\mathcal{R}_q}\chi_{Q_{m}}\leq\xn\chi_{\widetilde{Q}_{m+1}}\leq\xn\chi_{{Q}_{m+2}}.$$
Letting $m\rightarrow\infty$ we conclude that $\xm_{\mathcal{R}_q}\leq\xn.$ This completes the proof of (\ref{5.40}) and of assertion (v).\\
(vi) The measure $\xm_{\mathcal{R}_q}$ is essentially absolutely continuous relative to $\qcap.$ Clearly this assertion follows now from Proposition \ref{21}.\\
(vii) By (\ref{13})
$$u\leq[u]_{Q_n}+[u]_{Q^c_n}.$$
Now since $Q_n^c$ is $\frak T_q$-closed and $Q_n^c\downarrow\mathcal{R}_{q,0}^c(u),$ we have by Proposition \ref{17}-(iii) that
$$[u]_{Q^c_n}\downarrow[u]_{\mathcal{R}_{q,0}^c(u)}.$$
Hence
$$\lim_{n\to\infty}\left(u-[u]_{Q_n}\right)=u-v_{{\mathcal{R}_q}}\leq[u]_{\mathcal{R}_{q,0}^c(u)},$$
therefore $u\ominus v_{{\mathcal{R}_q}}\approx_{\mathcal{R}_{q,0}(u)}0.$ Since $v_{{\mathcal{R}_q}}\leq u$ this is equivalent to the statement $u\approx_{\mathcal{R}_{q,0}(u)}v_{{\mathcal{R}_q}}.$\\
(viii) (\ref{5.28}) follows by the previous statement. Now we assume that $\xm_{\mathcal{R}_q}(F)\chi_{K}<\infty$ for any compact $K\subset\mathbb{R}^N.$ Now set $F_n=F\cap\widetilde{Q}_n.$ By (\ref{13}).
$$[u]_F\leq[u]_{F_n}+[u]_{F\setminus F_n}=[u]_{F_n}+[u]_{F\setminus \widetilde{Q}_n }\leq[u]_{F_n}+[u]_{F\setminus Q_n }.$$
Now since $F\setminus Q_n$ is a $\frak T_q$-closed set and $\cap F\setminus Q_n=G$ with $\qcap(G)=0,$ we have by Proposition \ref{17}-(iii) that
$[u]_{F\setminus Q_n }\rightarrow[u]_{G}=0.$ Hence $[u]_F=\lim_{n\to\infty} [u]_{F_n},$ and $\mathrm{tr}([u]_{F_n})=\xm_{{\mathcal{R}_q}}\chi_{F_n}\uparrow\xm_{{\mathcal{R}_q}}\chi_{F_0}=\xm_{{\mathcal{R}_q}}\chi_{F}.$ Since $\xm_R\chi_F$ is a Radon measure essentially absolutely continuous relative to $\qcap,$ $[u]_F$ is moderate and (\ref{5.29}) holds.\\
(ix) If $\xm_{\mathcal{R}_q}(F)\chi_{K}<\infty$ for any compact $K\subset\mathbb{R}^N$ then, by (viii), $[u]_F$ is moderate. Conversely, if  $[u]_F$ is moderate, by (\ref{5.25}), $\xm_{\mathcal{R}_q}(F)\chi_{K}<\infty$ for any compact $K\subset\mathbb{R}^N.$\hfill$\Box$\\

\noindent \textbf{Example.} We give below an example which shows that there exists  $u\in\mathcal{U}_+(Q_T)$ such that ${\mathcal{R}_q}(u)= \mathbb{R}^N$ but $u$ is not a moderate solution. 
Let $\eta:[0,\infty)\rightarrow[0,\infty)$ be a smooth function such that $\eta(r)>0$ for any $r>0$ and $\lim_{r\rightarrow0^+}\eta(r)=0,$ ($\eta$ tends to $0$ very fast, for example $\eta(r)=e^{-\frac{1}{r^2}}$). Let $K$ be the close set
$$K=\{(x',x_n)\in\mathbb{R}^N: |x'|\leq\eta(x_n),\;x_n\geq0\}.$$
Then $K$ is thin at the origin $0.$

Set $f(x)=\frac{1}{\eta^{n}(x_n)}$ if $x\in K$ and $f=0$ otherwise. We define the measure
$$\xm=fdx.$$
This measure possesses the following properties:\\
1. $\xm$ is $\frak T_q$-locally finite.\\
2. $\xm(Q_n)<\infty$ where $Q_n=B_{2n}(0)\setminus\overline{B}_{\frac{1}{n}}(0)$ and $\cup Q_n\qeq \mathbb{R}^N$\\
3. $\xm(F)=0$ for any $F$ such that $\qcap(F)=0.$\\
4.  There exists a non decreasing sequence of bounded Radon measures $\xm_n$ absolutely continuous with $\qcap$ such that\\
(a) $\frak T_q$-supp$\,(\xm_n)\subset\widetilde{Q}_n,$ $\xm_n(A)=\xm_{n+k}(A)$ for any $A\subset\widetilde{Q}_n$ and any $n,k\in\mathbb{N}.$\\
(b) $\lim_{n\to\infty} \xm_n=\xm$\\
5. We can construct a solution $u\in\mathcal{U}_+(Q_T)$ with respect to this measure.\\
As we see later this solution is unique since it is $\xs$-moderate (see Proposition \ref{29}).
\begin{lemma}
Let $\xm$ satisfying the conditions 1-4 as above. Then there exists an open set ${\mathcal{R}_q}\qeq\mathbb{R}^N$ such that the measure $\xm$ is a Radon measure in ${\mathcal{R}_q}.$
\end{lemma}
\Proof
We consider the ball $B_R(0)$ with $R>1.$ From \cite[Lemma 2.5]{MV-CONT} there exists a sequence of open sets $\{O_{\frac{1}{m}}\}_{m=1}^\infty$ and $n(m)\in\mathbb{N}$ such that $\qcap(O_\frac{1}{m})<\frac{1}{m},$  and
\be
\overline{B}_R(0)\setminus O_\frac{1}{m}\subset \bigcup_{i=i}^{n(m)}Q_i.\label{1111}
\ee
Now since $O_\frac{1}{m}$ is open we have $$\qcap(\overline{O}_\frac{1}{m})=\qcap(\widetilde{O}_\frac{1}{m}\bigcup (\overline{O}_\frac{1}{m}\cap e_q(O)))\leq\qcap(\widetilde{O}_\frac{1}{m})\leq c\qcap(O_\frac{1}{m})\rightarrow0,$$
where $e_q(O)$ is the set of thin points of $O.$\\
Thus if $x\in B_R(0)\setminus \bigcap_{m=1}^\infty \overline{O}_\frac{1}{m}$ there exist $r>0$ small enough and $N\in\mathbb{N}$ such that
$$B_r(x)\subset B_R(0)\setminus \bigcap_{m=1}^N\overline{O}_\frac{1}{m}. $$
Thus by the properties of $\xm$ and (\ref{1111}) we have
$$\xm(B_r(x))<\infty.$$
We define
$${\mathcal{R}_q}:=\{x\in\mathbb{R}^N:\;\exists\;r>0\;\mathrm{such\;that}\;\xm(B_r(x))<\infty\}.$$
Then the set ${\mathcal{R}_q}$ is open and by the above argument, letting $R$ go to infinity, we have that ${\mathcal{R}_q}\qeq\mathbb{R}^N.$
Also by the definition of ${\mathcal{R}_q},$ it is easy to see that $\xm(K)<\infty$ for any compact $K\subset{\mathcal{R}_q}$ and by the properties of $\xm$ we can prove that $\xm$ is Radon measure in ${\mathcal{R}_q}.$\hfill$\Box$
\subsection{The precise initial trace}
We are now in condition to define the {\it precise initial trace.}
\begin{defin}
Let $q\geq1+\frac{2}{N}$ and $u\in\mathcal{U}_+(Q_T).$\\
\textbf{a:} The solution $v_{\mathcal{R}_q}$ defined by (\ref{5.21}) is called \textbf{regular component of $u$} and will be denoted by $u_{reg}.$\\
\textbf{b}: Let $\{v_n\}$ be an increasing sequence of moderate solutions satisfying condition (\ref{5.22}) and put $\xm_{{\mathcal{R}_q}}=\xm_{{\mathcal{R}_q}}(u):=\lim_{n\to\infty}\mathrm{tr}(v_n).$ Then, the regularized measure $\overline{\xm}_{{\mathcal{R}_q}},$ defined by (\ref{5.27}), is called the \textbf{regular initial trace of $u$}. It will be denoted by $\mathrm{tr}_{{\mathcal{R}_q}}(u).$\\
\textbf{c:} The couple $(\mathrm{tr}_{{\mathcal{R}_q}}(u),\mathcal{S}_q(u))$ is called \textbf{the precise initial trace of $u$} and will be denoted by $\mathrm{tr}^c(u).$\\
\textbf{d:} Let $\xn$ be the Borel measure on $\mathbb{R}^N$ given by
\be
\xn=\Bigg \{
\begin{array}{ll}(\mathrm{tr}_{{\mathcal{R}_q}}(u))(E)\;\;\;\;\mathrm{if} \;E\subset{\mathcal{R}_q}(u),
\\[2pt]
\xn(E)=\infty\;\;\;\; \mathrm{if} \;E\cap\mathcal{S}_q(u)\neq\emptyset,
\end{array}\label{5.42}
\ee
for every Borel set $E.$ Then $\xn$ is the measure representation of the precise trace of $u,$ to be denoted by $\mathrm{tr}(u).$\label{25}
\end{defin}
Note that, by Proposition \ref{24}-(v), the measure $\xm_{{\mathcal{R}_q}}$ is independent of the choice of the sequence $\{v_n\}.$
\begin{theorem}
Assume that $u\in\mathcal{U}_+(Q_T)$ is a $\xs$-moderate solution, i.e., there exists an increasing sequence $\{u_n\}$ of positive moderate solutions such that $u_n\uparrow u.$ Let $\xm_n=\lim_{n\to\infty}\mathrm{tr}(u_n),$ $\xm_0:=\lim_{n\to\infty}\xm_n$ and set, for any Borel set $E$, 
\be
\xm(E)=\inf\left\{\xm_0(Q):\;E\subset Q,\quad Q\;\,\frak T_q\text{-open}\right\}.\label{5.43}
\ee
Then:\\
(i) $\xm$ is the precise initial trace of $u$ and $\xm$ is $\frak T_q$-perfect. In particular $\xm$ is independent of the sequence $\{u_n\}$ which appears in its definition.\\
(ii) If $A$ is a Borel set such that $\xm(A)<\infty$, then $\xm(A)=\xm_0(A).$\\
(iii) A solution $u\in\mathcal{U}_+(Q_T)$ is $\xs$-moderate if and only if
\be
u=\sup\{v\in\mathcal{U}_+(Q_T):\;v\;\mathrm{moderate},\;\;\;v\leq u\},\label{5.44}
\ee
which is equivalent to
\be
u=\sup\{u_\tau\in\mathcal{U}_+(Q_T):\;\tau\in W^{-\frac{2}{q},q}(\mathbb{R}^N)\cap\mathfrak{M}^b_+(\mathbb{R}^N),\;\;\;\tau\leq \mathrm{tr}(u)\}.\label{5.45}
\ee
(iv) If $u,\;w$ are $\xs$-moderate solutions,
\be
\mathrm{tr}(w)\leq\mathrm{tr}(u)\Leftrightarrow w\leq u.\label{5.46}
\ee\label{27}
\end{theorem}
\Proof    The proof is an adaptation of the one in \cite{MV-CONT}.\\
(i) Let $Q$ be a $\frak T_q$-open set and $A$ a Borel set such that $\qcap(A)=0.$ Then $\xm_n(A)=0$ so that $\xm_0(A)=0.$ Thus $\xm_0$ is essentially absolutely continuous and, by Proposition \ref{21}, $\xm$ is $\frak T_q$-perfect.

Let $\{D_n\}$ be the family of $\frak T_q$-open sets as in Proposition \ref{24}-(i). Put $D_n'={\mathcal{R}_q}(u)\setminus D_n$ and observe that $D_n'\downarrow E$ where $\qcap(E)=0.$ Therefore, for fixed $n$, 
$$u_{\xm_n\chi_{D'_m}}\downarrow0\quad\text{when }m\to\infty.$$
Thus there exist a subsequence, say $\{D'_n\},$ such that
$$\left(\int_0^T\int_{B_n(0)}|u_{\xm_n\chi_{D'_n}}|^qdxdt\right)^{\frac{1}{q}}\leq2^{-n}.$$
Since,
$$\xm_n({\mathcal{R}_q}(u))=\xm_n\chi_{D_n}+\xm_n\chi_{D'_n},$$
it follows that
$$\lim_{n\to\infty}\left|u_{\xm_n\chi_{{\mathcal{R}_q}(u)}}-u_{\xm_n\chi_{D_n}}\right|\leq \lim_{n\to\infty} u_{\xm_n\chi_{D'_n}}=0.$$
Thus
$$u_n\leq u_{\xm_n\chi_{D_n}}+[u]_{\mathcal{S}_q(u)}.$$
Hence
$$u-[u]_{\mathcal{S}_q(u)}\leq w:=\lim_{n\rightarrow\infty}u_{\xm_n\chi_{{\mathcal{R}_q}(u)}}=
\lim_{n\rightarrow\infty}u_{\xm_n\chi_{D_n}}\leq u_{reg}.$$
This implies $u\ominus[u]_{\mathcal{S}_q(u)}\leq u_{reg}.$ For the opposite inequality, by Proposition \ref{24}-(iv) we get
$$[u]_{D_n}\uparrow u_{reg}.$$
But by (\ref{4.20}) and using the facts that $\widetilde{D}_n\qsub D_{n+1}\subset\widetilde{D}_{n+1}\qsub{\mathcal{R}_q}(u),$ $\qcap\left(\widetilde{D}_{n+1}\cap\mathcal{S}_q(u)\right)=0,$
$$[u]_{D_n}\leq[[u]_{\mathcal{S}_q(u)}]_{D_{n+1}}+[u\ominus[u]_{\mathcal{S}_q(u)}]_{D_{n+1}}=[u\ominus[u]_{\mathcal{S}_q(u)}]_{D_{n+1}}\leq u\ominus[u]_{\mathcal{S}_q(u)}.$$
Letting $n\rightarrow\infty$ we derive $u_{reg}\leq u\ominus[u]_{\mathcal{S}_q(u)}.$ Therefore $\lim_{n\to\infty} u_{\xm_n\chi_{D_n}}=u_{reg}.$ Thus the sequence $\{u_{\xm_n\chi_{D_n}}\}$ satisfies condition (\ref{5.22}) and consequently, by Proposition \ref{24}-(iv) and Definition \ref{25},
\be
\lim_{n\to\infty}\xm_n\chi_{D_n}=\xm_{{\mathcal{R}_q}},\quad\mathrm{tr}_{\mathcal{R}_q}(u)=\overline{\xm}_{{\mathcal{R}_q}}.\label{5.47}
\ee
Next, we claim that, if $\xi\in \mathcal{S}_q(u)$ then, for every $\frak T_q$-open bounded neighborhood $Q$ of $\xi$ $\xm_n(\widetilde{Q})\rightarrow\infty.$ Indeed
let $\eta\in\besl$ with $\frak T_q$-support in $Q.$  Put $h=\mathbb{H}[\eta]$ and $\xf(r)=r^{2q'}_+.$ Then by Proposition \ref{moderate}, Lemma \ref{lem2.4} and in view of the proof of Lemma \ref{lem2.3} we have
$$
\int_0^T\int_{\mathbb{R}^N}(-u_n(\partial_t\xf(h)+\xD\xf(h)))+u^q_n\xf(h)dxd\tau+\int_{\mathbb{R}^N}u_n\xf(h)(.,T)dx=\int_{Q}\eta^{2q'}d\xm_n.
$$
In view of Lemma \ref{lem2.3}, we can prove
$$\int_0^T\int_{\mathbb{R}^N}u^q_n\xf(h)dxd\tau\leq C(q)\left(\int_{Q}\eta^{2q'}d\xm_n+||\eta||^{2q'}_{W^{\frac{2}{q},q'}}+||\eta||_{L^\infty}\right).$$
By Lemma \ref{cutoff} there exist $\eta\in\besl$ and a $\frak T_q$-open set $D\subset Q$ such that $\eta=1$ on $D,$  $\eta=0$ outside of $Q$ and $0\leq\eta\leq1.$ Letting $n\rightarrow\infty$ we have
$$\lim_{n\to\infty}\int_0^T\int_{\mathbb{R}^N}u^q_n\mathbb{H}^{2q'}[\chi_D]dxd\tau\leq C(q)\left(\lim_{n\to\infty}\int_{Q}\eta^{2q'}d\xm_n+||\eta||^{2q'}_{W^{\frac{2}{q},q'}}+||\eta||_{L^\infty}\right),$$
the assertion follows by Lemma \ref{lem2.5}.

In conclusion, if $\xi\in\mathcal{S}_q(u)$ then $\xm_0(\widetilde{Q})=\infty$ for every $\frak T_q$-open neighborhood of $\xi.$ Consequently $\xm(\xi)=\infty.$ This fact and (\ref{5.47}) imply that $\xm$ is the precise trace of $u.$\\
(ii) If $\xm(A)<\infty$ then $A$ is contained in a $\frak T_q$-open set $D$ such that $\xm_0(D)<\infty$ and, by Proposition \ref{21}, $\xm(A)=\xm_0(A).$\\
(iii) Let $u\in\mathcal{U}_+(Q_T)$ be $\xs$-moderate and put
\be
u^*:=\sup\{v:\;v\;\mathrm{moderate},\;\;\;v\leq u\}.\label{5.48}
\ee
By its definition $u^*\leq u.$ On the other hand, since there exists an increasing sequence of moderate solutions $\{u_n\}$ converging to $u,$ it follows that $u\leq u^*.$ Thus $u=u^*.$

Conversely, if $u\in\mathcal{U}_+(Q_T)$ and $u=u^*$ then by proposition \ref{sygklish}, there exists an increasing sequence of moderate solutions $\{u_n\}$ converging to $u.$ Therefore $u$ is $\xs$-moderate.

Since $u$ is $\xs$-moderate there exist an increasing sequence of moderate solutions $\{u_n\}$ converging to $u.$ In view of the discussion at the beginning of subsection \ref{moderatesec}, for any $u_n$ there exist an increasing sequence of $\{w_m\}$ such that $w_m\uparrow u_n$ and $\mathrm{tr}(w_m)\in W^{-\frac{2}{q},q}(\mathbb{R}^N)\cap\mathfrak{M}^b_+(\mathbb{R}^N).$ Thus
$$u_n\leq\sup\{u_\tau:\;\tau\in W^{-\frac{2}{q},q}(\mathbb{R}^N)\cap\mathfrak{M}^b_+(\mathbb{R}^N),\;\;\;\tau\leq \mathrm{tr}(u)\}=:u^\ddag.$$
Letting $n\rightarrow\infty,$ we have $u\leq u^\ddag.$

On the other hand, if $u$ is $\xs$-moderate, $\tau\in W^{-\frac{2}{q},q}(\mathbb{R}^N)\cap\mathfrak{M}^b_+(\mathbb{R}^N)$ and $\tau\leq \mathrm{tr}(u)$ then (with $\xm_n$ and $u_n$ as in the statement of the Proposition), $\mathrm{tr}(u_\tau\ominus u_n)=(\tau-\xm_n)_+\downarrow0.$ Hence, $u_\tau\ominus u_n\downarrow0,$ which implies, $u_\tau\leq u.$ Therefore $u^\ddag\leq u.$ Thus (\ref{5.44}) implies (\ref{5.45}) and each of them that $u$ is $\xs$-moderate. Therefore the two statements are equivalent.\\
(iv) The assertion $\Rightarrow$ is a consequence of (\ref{5.45}). To verify the assertion $\Leftarrow$ it is sufficient to show that if $w$ is moderate, $u$ is $\xs$-moderate and $w\leq u$, then $\mathrm{tr}(w)\leq\mathrm{u}.$ Let $\{u_n\}$ be an increasing sequence of positive moderate solutions converging to $u.$ Then $u_n\vee w\leq u$ and consequently $u_n\leq u_n\vee w\uparrow u.$ Therefore $\mathrm{tr}(u_n\vee w)\uparrow\xm'\leq\mathrm{tr}(u)$ so that $\mathrm{tr}(w)\leq\mathrm{tr}(u).$ \hfill$\Box$

\begin{theorem}
Let $u\in\mathcal{U}_+(Q_T)$ and put $\xn=\mathrm{tr}(u).$\smallskip

\noindent (i) $u_{reg}$ is $\xs$-moderate and $\mathrm{tr}(u_{reg})=\mathrm{tr}_{{\mathcal{R}_q}}(u).$\smallskip

\noindent
(ii) If $v\in\mathcal{U}_+(Q_T)$
\be
v\leq u\Rightarrow\mathrm{tr}(v)\leq\mathrm{tr}(u).\label{5.49}
\ee
If $F$ is a $\frak T_q$-closed set, then
\be
\mathrm{tr}([u]_F)\leq \xn\chi_F.\label{5.50}
\ee
(iii) A singular point can be characterized in terms of the measure $\xn$ as follows:
\be
\xi\in\mathcal{S}_q(u)\Leftrightarrow \xn(Q)=\infty\qquad\forall Q:\;\;\xi\in Q,\;Q\;\frak T_q\text{-open}.\label{5.51}
\ee
(iv) If $Q$ is a $\frak T_q$-open set then:
\bea
[u]_Q\;\mathrm{is\;moderate}&\Leftrightarrow&\;\exists\; \mathrm{Borel\;set}\;A:\;\qcap(A)=0,\;\xn(K\cap\widetilde{Q}\setminus A)<\infty,\label{5.53}
\eea
for any compact $K\subset\mathbb{R}^N.$\smallskip

\noindent
(v) The singular set of $u_{reg}$ may not be empty. In fact
\be
\mathcal{S}_q(u)\setminus b_q(\mathcal{S}_q(u))\subset\mathcal{S}_q(u_{reg})\subset\mathcal{S}_q(u)\cap\widetilde{{\mathcal{R}_q}(u)},\label{5.54}
\ee
where $b_q(\mathcal{S}_q(u))$ is the set of $\qcap$-thick points of $\mathcal{S}_q(u).$\smallskip

\noindent
(vi) Put
\be
\mathcal{S}_{q,0}(u):=\{\xi\in\mathbb{R}^N:\;\xn(Q\setminus\mathcal{S}_q(u))=\infty\;\;\;\;\forall\;Q\;\;\frak T_q\text{-open neighborhood of } \xi\}.\label{5.55}
\ee
Then
\be
\mathcal{S}_q(u_{reg})\setminus b_q(\mathcal{S}_q(u))\subset\mathcal{S}_{q,0}(u)\subset\mathcal{S}_q(u_{reg})\bigcup b_q(\mathcal{S}_q(u)).\label{5.56}
\ee\label{29}
\end{theorem}

\noindent\emph{Remark.} This results extends Proposition \ref{24} which deals with the regular initial trace.

\medskip

\noindent \Proof
(i) By proposition \ref{24}-(ii) $u_{reg}$ is $\xs$-moderate. The second part of the statement follows from Definition \ref{25} and Proposition \ref{27}-(i).\\
(ii) If $v\leq u$ then ${\mathcal{R}_q}(u)\subset{\mathcal{R}_q}(v)$ and by definition $v_{reg}\leq u_{reg}.$ By Proposition \ref{27}-(iv) $\mathrm{tr}(v_{reg})\leq\mathrm{tr}(u_{reg})$ and consequently $\mathrm{tr}(v)\leq\mathrm{tr}(u).$ Inequality (\ref{5.50}) is an immediate consequence of (\ref{5.49}).\\
(iii) If $\xi\in{\mathcal{R}_q}(u)$ there exists a $\frak T_q$-open  neighborhood $Q$ of $\xi$ such that $[u]_Q$ is moderate. Hence $\xn(Q)=\mathrm{tr}_{{\mathcal{R}_q}}(u)(Q)<\infty.$ If $\xi\in\mathcal{S}_q(u),$ it follows immediately from the definition of precise trace that $\xn(Q)=\infty$ for every $\frak T_q$-open neighborhood $Q$ of $\xi.$\\
(iv) If $[u]_Q$ is moderate then $Q\subset {\mathcal{R}_q}(u).$ Proposition \ref{24}-(ix) implies (\ref{5.53}) in the direction $\Rightarrow$. On the other hand,
$$\xn(K\cap\widetilde{Q}\setminus A)<\infty,\;\forall\;\mathrm{compact}\;K\subset\mathbb{R}^N\Rightarrow\widetilde{Q}\qsub{\mathcal{R}_q}(u),$$
and $\xm_{{\mathcal{R}_q}}(K\cap\widetilde{Q})=\xm_{{\mathcal{R}_q}}(K\cap\widetilde{Q}\setminus A)<\infty.$ Hence, by Proposition \ref{24}-(ix), $[u]_Q$ is moderate.\\
(v) Since $\frak T_q\text{-supp}\,(u_{reg})\subset\widetilde{{\mathcal{R}_q}(u)}$ and ${\mathcal{R}_q}(u)\subset{\mathcal{R}_q}(u_{reg})$ we have
$$\mathcal{S}_q(u_{reg})\subset\mathcal{S}_q(u)\cap\widetilde{{\mathcal{R}_q}(u)}.$$
Next we show that $\mathcal{S}_q(u)\setminus b_q(\mathcal{S}_q(u))\subset\mathcal{S}_q(u_{reg}).$

If $\xi\in\mathcal{S}_q(u)\setminus b_q(\mathcal{S}_q(u))$ then ${\mathcal{R}_q}(u)\cup\{\xi\}$ is a $\frak T_q$-open neighborhood of $\xi.$ By (i) $u_{reg}$ is $\xs$-moderate and consequently (by Proposition \ref{27}-(i)) its trace is $\frak T_q$-perfect. Therefore, if $Q_0$ is a bounded $\frak T_q$-open neighborhood of $\xi$ and $Q=Q_0\cap\left(\{\xi\}\cup{\mathcal{R}_q}(u)\right)$ then
$$\mathrm{tr}(u_{reg})(Q)=\mathrm{tr}(u_{reg})(Q\setminus\{\xi\})=\mathrm{tr}(u)(Q\setminus\{\xi\}),$$
where in the last equality we have used the fact that $Q\setminus\{\xi\}\subset{\mathcal{R}_q}(u).$ Let $D$ be a $\frak T_q$-open set such that $\xi\in D\subset\widetilde{D}\subset Q.$ If $\mathrm{tr}(u)(\widetilde{D}\setminus\{\xi\})<\infty$ then, by (iv), $[u]_D$ is moderate and $\xi\in{\mathcal{R}_q}(u),$ contrary to our assumption. Therefore $\mathrm{tr}(u)(\widetilde{D}\setminus\{\xi\})=\infty$ so that $\mathrm{tr}(u_{reg})(Q_0\setminus\{\xi\})=\infty$ for every $\frak T_q$-open bounded neighborhood $Q_0$ of $\xi,$ which implies $\xi\in\mathcal{S}_q(u_{reg}).$ This completes the proof of (\ref{5.54}).\\
(vi) If $\xi\notin b_q(\mathcal{S}_q(u))$, there exists a $\frak T_q$-open neighborhood $D$ of $\xi$ such that $\left(D\setminus\{\xi\}\right)\cap\mathcal{S}_q(u)=\emptyset$ and consequently
\be
\mathrm{tr}(u_{reg})(D\setminus\{\xi\})=\mathrm{tr}(u_{reg})(D\setminus\mathcal{S}_q(u))=\mathrm{tr}(u)(D\setminus\mathcal{S}_q(u)).\label{5.57}
\ee
If, in addition $\xi\in\mathcal{S}_{q,0}(u)$ then
$$\mathrm{tr}(u)(D\setminus\mathcal{S}_q(u))=\mathrm{tr}(u_{reg})(D\setminus\mathcal{S}_q(u))=\infty.$$
If $Q$ is an arbitrary $\frak T_q$-open neighborhood of $\xi$ then the same holds if $D$ is replaced by $Q\cap D.$ Therefore $\mathrm{tr}(u_{reg})(Q\setminus\{\xi\})=\infty$ for any such $Q.$ Consequently $\xi\in\mathcal{S}_q(u_{reg}),$ which proves that $\mathcal{S}_{q,0}(u)\setminus b_q(\mathcal{S}_q(u))\subset\mathcal{S}_q(u_{reg}).$

On the other hand, if $\xi\in\mathcal{S}_q(u_{reg})\setminus b_q(\mathcal{S}_q(u))$ then there exists a $\frak T_q$-open neighborhood $D$ such that (\ref{5.57}) holds and $\mathrm{tr}(u_{reg})(D)=\infty.$ Since $u_{reg}$ is $\xs$-moderate, $\mathrm{tr}(u_{reg})$ is $\frak T_q$-perfect so that $\mathrm{tr}(u_{reg})(D)=\mathrm{tr}(u_{reg})(D\setminus\{\xi\})=\infty.$ Consequently, by (\ref{5.57}), $\mathrm{tr}(u)(D\setminus\mathcal{S}_q(u))=\infty.$ If $Q$ is any $\frak T_q$-open neighborhood of $\{\xi\}$ then $D$ can be replaced by $D\cap Q.$ Consequently $\mathrm{tr}u(Q\setminus\mathcal{S}_q(u))=\infty$ and we conclude that $\xi\in\mathcal{S}_{q,0}(u).$ This completes the proof of (\ref{5.56}).\hfill$\Box$
\begin{prop}
Let $F$ be a $\frak T_q$-closed set. Then $\mathcal{S}_q(U_F)=b_q(F).$\label{UF}
\end{prop}
\Proof
Let $\xi\in\mathbb{R}^N$ such that $F$ is $\qcap$-thin at $\xi.$ Let $Q$ be a $\frak T_q$-open neighborhood of $\xi$ such that $\widetilde{Q}\qsub F^c.$ Then $[U_F]_Q=U_{F\cap\widetilde{Q}}=0.$ Therefore $\xi\in{\mathcal{R}_q}(U_F).$

Conversely, assume that $\xi\in F\cap{\mathcal{R}_q}(U_F),$ thus there exists a $\frak T_q$-open neighborhood $Q$ of $\xi$ such that $[U_F]_Q$ is moderate. But $[U_F]_Q=U_{F\cap\widetilde{Q}}$ which implies $\qcap(F\cap\widetilde{Q})=0$ and $Q\subset\mathcal{R}(u).$ Now, note that $\qcap(F)\leq\qcap(F\cap Q)+\qcap(Q^c).$ Thus $F$ is $\frak T_q$-thin at $\xi.$\hfill$\Box$

\subsection{The initial value problem}
The following notations will be used in the sequel.

\begin{notation}
a: We denote by $\mathfrak{M}_+(\mathbb{R}^N)$ the space of positive outer regular Borel measure on $\mathbb{R}^N.$\smallskip

\noindent b: We denote by $\mathcal{C}_q(\mathbb{R}^N)$ the space of couples $(\tau,F)$ such that $F$ is $\frak T_q$-closed, $\tau\in\mathfrak{M}_+(\mathbb{R}^N),$ $\frak T_q$-{supp}$\,(\tau)\subset \widetilde{F^c}$ and $\tau\chi_{F^c}$ is $\frak T_q$-locally finite.\smallskip

\noindent 
c: We denote by $\mathbb{T}:\mathcal{C}_q(\mathbb{R}^N)\rightarrow\mathfrak{M}_+(\mathbb{R}^N)$ the mapping given by $\xn=\mathbb{T}(\tau,F)$ where $\xn$ is defined as in (\ref{5.42}) with ${\mathcal{R}_q}(u),\;\mathcal{S}_q(u)$ replaced by $F,\;F^c$ respectively. In this setting $\xn$ is the measure representation of the couple $(\tau,F)$.\smallskip

\noindent 
d: If $(\tau,F)\in\mathcal{C}_q(\mathbb{R}^N)$ the set
\be
F_\tau=\{\xi\in\mathbb{R}^N:\;\tau(Q\setminus F)=\infty\quad\forall Q\;\;\frak T_q\text{-open neighborhood of }\;\xi \}
\ee
is called the set of explosion points of $\tau.$
\end{notation}
\emph{Remark.} Note that $F_\tau\subset F$ (because $\tau\chi_{F^c}$ is $\frak T_q$-locally finite) and $F_\tau\subset\widetilde{F^c}$ (because $\tau$ vanishes outside this set). Thus
\be
F_\tau\subset b_q(F^c)\cap F.\label{5.59}
\ee
\begin{prop} Let $\xn$ be a positive Borel measure on $\mathbb{R}^N.$\\
(i) The initial value problem
\be
\prt_t  u-\xD u+|u|^{q-1}u=0,\;\;u>0\;\mathrm{in}\;Q_\infty=\mathbb{R}^N\times(0,T),\;\;\mathrm{tr}(u)=\xn\;\mathrm{in}\;\mathbb{R}^N\times\{0\}.\label{5.60}
\ee
possesses a solution if and only if $\xn\in\mathfrak{M}_q(\mathbb{R}^N).$\\
(ii) Let $(\tau,F)\in \mathcal{C}_q(\mathbb{R}^N)$ and put $\xn:=\mathbb{T}(\tau,F).$ Then $\xn\in\mathfrak{M}_q(\mathbb{R}^N)$ if and only if
\be
\tau\in\mathfrak{M}_q(\mathbb{R}^N),\;\;\;\;F=b_q(F)\bigcup F_\tau.\label{5.61}
\ee
(iii) Let $\xn\in\mathfrak{M}_q(\mathbb{R}^N)$ and denote
\bea
\nonumber
\mathcal{E}_\xn&:=&\{E:\;E\;\,\frak T_q\text{-quasi-closed},\;\xn(E\cap K)<\infty,\;\forall\;\mathrm{compact}\;K\subset\mathbb{R}^N\}\\
\mathcal{D}_\xn&:=&\{D:\;D\;\,\frak T_q\text{-open},\widetilde{D}\qeq E\;\mathrm{for\;some\;}E\in\mathcal{E}_\xn\}.\label{5.62}
\eea
Then a solution of (\ref{5.60}) is given by $u=v\oplus U_F$ where
\be
G:=\bigcup_{\mathcal{D}_\xn}D,\;\;\;\;F:=G^c,\;\;\;\;v:=\sup\{u_{\xn}\chi_E:\;E\in\mathcal{E}_\xn\}.\label{5.63}
\ee

Note that if $E\in\mathcal{E}_\xn$ then $\xn\chi_{E}$ is locally bounded Borel measure which does not charge sets of $\qcap$-capacity zero. Recall that if $\xm$ is a positive measure possessing these properties, then $u_\xm$ denotes the moderate solution with initial trace $\xm.$\\
(iv) The solution $u=v\oplus U_F$ is $\xs$-moderate and it is the unique solution of problem (\ref{5.60}) in the class of $\xs$-moderate solutions. Furthermore, $u$ is the largest solution of the problem.\label{35}
\end{prop}
\Proof    The proof is similar to the one in \cite{MV-CONT}.\\
\textbf{(A)} \emph{If} $u\in\mathcal{U}_+(Q_T)$
\be
\mathrm{tr}(u)=\xn\Rightarrow\xn\in\mathfrak{M}_q(\mathbb{R}^N).\label{5.64}
\ee
By Proposition \ref{24}, $u_{reg}$ is $\xs$-moderate and $u\approx_{{\mathcal{R}_q}(u)}u_{reg}.$ Therefore
$$\mathrm{tr}(u)\chi_{{\mathcal{R}_q}(u)}=\mathrm{tr}(u_{reg})\chi_{{\mathcal{R}_q}(u)}.$$
By Proposition \ref{27}, $\overline{\xm}_{\mathcal{R}_q}:=\mathrm{tr}(u_{reg})\in \mathfrak{M}_q(\mathbb{R}^N).$ If $v$ is defined as in (\ref{5.63}) then
\be
v=\sup\{[u]_F:\;F\;\text{ $\frak T_q$-closed}\; F\qsub \mathcal{R}_q(u)\}=u_{reg},\label{5.65}
\ee
where the second equality holds by definition. Indeed, by Theorem \ref{29}, for every $\frak T_q$-open set $Q,$ $[u]_Q$ is moderate if and only if $\xn(K\cap\widetilde{Q}\setminus A)<\infty$ for some set $A$ with capacity zero and for any compact $K$ subset of $\mathbb{R}^N.$ This means that $[u]_Q$ is moderate if and only if there exists $E\in\mathcal{E}_\xn$ such that $\widetilde{Q}\qeq E.$ When this is the case,
$$\mathrm{tr}([u]_Q)=\xm_{\mathcal{R}_q}(u)\chi_{\widetilde{Q}}=\xm_{\mathcal{R}_q}(u)\chi_{E}=\xn\chi_E.$$
Thus $v\geq u_{reg}.$ On the other hand, if $E\in\mathcal{E}_\xn,$ then $\widetilde{E}\qsub\mathcal{R}_q(u)$ and $\xm_{\mathcal{R}_q}(u)(K\cap\widetilde{E})=\xm_{\mathcal{R}_q}(u)(K\cap E)<\infty$ for any compact $K$ subset of $\mathbb{R}^N.$ Therefore by Proposition \ref{24}-(ix), $\widetilde{E}$ is regular, i.e, there exist a $\frak T_q$-open regular set $Q$ such that $E\qsub Q.$ Hence $u_{\xn\chi_E}\leq [u]_Q$ and we conclude that $v\leq u_{reg}.$ This proves (\ref{5.65}). In addition, if $E\cap\mathcal{S}_q(u)\neq\emptyset$ then $\xn(E)=\infty$, by Definition \ref{25}. Therefore $\xn$ is outer regular with respect to $\frak T_q$-topology.

Next we must prove that $\xn$ is essentially absolutely continuous. Let $Q$ be a $\frak T_q$-open set and $A$ a non-empty $\frak T_q$-closed subset of $Q$ such that $\qcap(A)=0.$ Either $\xn(Q\setminus A)=\infty$, in which case $\xn(Q\setminus A)=\xn(Q)$, or $\xn(Q\setminus A)<\infty$, in which case $Q\setminus A\subset\mathcal{R}_q(u)$ and
$$\xn(Q\setminus A)=\overline{\xm}_{\mathcal{R}_q}(Q\setminus A)=\overline{\xm}_{\mathcal{R}_q}(Q)<\infty.$$
Let $\xi\in A$ let $D$ be a $\frak T_q$-open subset of $Q$ such that $\xi\in D\subset\widetilde{D}\qsub Q.$ Let $B_n$ be a $\frak T_q$-open neighborhood of $A\cap\widetilde{D}$ such that $\qcap(B_n)<2^{-n}$ and $B_n\qsub D.$ Then
$$[u]_D\leq[u]_{E_n}+[u]_{B_n},\qquad E_n=\widetilde{D}\setminus B_n.$$
Since $\lim_{n\to\infty}[u]_{B_n}=0$ it follows that $[u]_D=[u]_{E_n}.$ Now since $E_n\subset \mathcal{R}_q(u),$ $\xn(E_n)\leq\xn(Q\setminus A)<\infty,$ we have by definition of $\xn$ and Proposition \ref{24}-(ix) that $[u]_{E_n}$ is moderate. Also in view of Lemma 2.8 and Lemma 2.7(ii)-\cite{MV-CPDE}, we have for some positive constant $C$
$$\int_0^T\int_{K}[u]_{E_n}^qdxdt\leq C\xn(E_n)\leq C\xn(Q\setminus A)<\infty,$$
for any compact $K\subset\mathbb{R}^N.$
Therefore
$$\int_0^T\int_{K}[u]_{D}^qdxdt<\infty\qquad\forall\;K\subset\mathbb{R}^N, K\text{ compact}.$$
which implies that $[u]_D$ is moderate and thus $D\subset \mathcal{R}_q(u).$ Since every point $A$ has a neighborhood $D$ as above we conclude that $A\subset\mathcal{R}_q(u)$ and hence $\xn(A)=\mathrm{tr}_{\mathrm{R}}(u)(A)=0.$
If $A$ is any a non-empty Borel subset of $Q$ such that $\qcap(A)=0,$ by inequality $\qcap(\widetilde{A})\subset c\qcap(A),$ we have that $\xn$ is
absolutely continuous and $\xn\in\mathfrak{M}_q(\mathbb{R}^N).$

Secondly we prove:\\
\textbf{(B)}\emph{ Suppose that $(\tau,F)\in\mathcal{C}_q(\mathbb{R}^N)$ satisfies (\ref{5.61}) and put $\xn=\mathbb{T}(\tau,F)$. Then the solution $u=v\oplus U_F$, with $v$ as in (\ref{5.63}), satisfies $\mathrm{tr}(u)=\xn.$}

\emph{By (\ref{5.64}), this also implies that $\xn\in\mathfrak{M}_q(\mathbb{R}^N).$}

Clearly $v$ is a $\xs$-moderate solution. The fact that $\tau$ is $\frak T_q$-locally finite in $F^c$ and essentially absolutely continuous relative to $\qcap$ implies that
\be
G:=F^c\subset\mathcal{R}_q(v),\;\;\;\;\mathrm{tr}(v)\chi_G=\tau_G.\label{5.66}
\ee
It follows from the definition of $v$ that $F_\tau\subset\mathcal{S}_q(v).$ Hence, by Proposition \ref{UF} and (\ref{5.54}) we have
\be
F=b_q(F)\bigcup F_\tau\subset\mathcal{S}_q(v)\bigcup\mathcal{S}_q(U_F)\subset\mathcal{S}_q(u)\subset F.\label{5.67}
\ee
Thus, $\mathcal{S}_q(u)=F,$ $v=u_{reg}$ and $\tau=\mathrm{tr}(u_{reg}).$ Thus $\mathrm{tr}^c(u)=(\tau,F)$ which is equivalent to $\mathrm{tr}(u)=\xn.$

Next we show:
\textbf{(C)}\emph{ Suppose that $(\tau,F)\in\mathcal{C}_q(\mathbb{R}^N)$ and that there exists a solution $u$ such that $\mathrm{tr}^c(u)=(\tau,F).$ Then
\be
\tau=\mathrm{tr}_{\mathcal{R}_q}(u)=\mathrm{tr}(u_{reg}),\;\;\;\;F=\mathcal{S}_q(u).\label{5.68}
\ee
If $U:=u_{reg}\oplus U_F$ then $\mathrm{tr}(U)=\mathrm{tr}(u)$ and $u\leq U$.  U is the unique $\xs$-moderate solution of (\ref{5.60}) and $(\tau,F)$ satisfies condition (\ref{5.61}).}
Assertion (\ref{5.68}) follows by Proposition \ref{24}-(i) and Definition \ref{25}. Since $u_{reg}$ is $\xs$-moderate, it follows, by Theorem \ref{27}, that $\tau\in\mathfrak{M}_q(\mathbb{R}^N).$

By Proposition \ref{24} (vi), $u\approx_{\mathcal{R}_q(u)}u_{reg}.$ Therefore $w:=u\ominus u_{reg}$ vanishes on $\mathcal{R}_q(u)$ so that $w\leq U_F.$ Note that $u-u_{reg}\leq w$ and therefore
\be
u\leq u_{reg}\oplus w\leq U.\label{5.69}
\ee
By their definitions $\mathcal{S}_{q,0}(u)=F_\tau$ and by Theorem \ref{29} (vi) and Proposition \ref{UF},
\bea
\nonumber
\mathcal{S}_q(U)&=&\mathcal{S}_q(u_{reg})\bigcup\mathcal{S}_q(U_F)=\mathcal{S}_q(u_{reg})\bigcup b_q(U_F)\\
              &=&\mathcal{S}_{q,0}(u)\bigcup b_q(U_F)= F_\tau\bigcup b_q(U_F).\label{5.70}
\eea
On the other hand, $\mathcal{R}_q(U)\supset\mathcal{R}_q(u_{\mathcal{R}_q})=\mathcal{R}_q(u)$ and, as $u\leq U,$ $\mathcal{R}_q(U)\subset\mathcal{R}_q(u).$ Hence $\mathcal{R}_q(U)=\mathcal{R}_q(u)$ and $\mathcal{S}_q(U)=\mathcal{S}_q(u).$ Therefore, by (\ref{5.68}) and (\ref{5.70}), $F=\mathcal{S}_q(U)=F_\tau\cup b_q(U_F).$ Thus $(\tau,F)$ satisfies (\ref{5.61}) and $\mathrm{tr}^c(U)=(\tau,F).$ The fact that $U$ is the maximal solution with this trace follows from (\ref{5.69}).

The solution $U$ is $\xs$-moderate because both $u_{reg}$ and $U_F$ are $\xs$-moderate solutions (concerning $U_F,$  see Proposition \ref{31}).

The uniqueness of the solution in the class of $\xs$-moderate solutions follows from Proposition \ref{27}-(iv).

Finally we prove:\\
\textbf{(D)} \emph{If $\xn\in\mathfrak{M}_q(\mathbb{R}^N)$ then the couple $(\tau,F)$ defined by
\be
v:=\sup\{u_{\xn}\chi_E:\;E\in\mathcal{E}_\xn\},\;\;\;\tau:=\mathrm{tr}(v),\;\;F=\mathcal{R}_q^c(v),\label{5.71}
\ee
satisfies (\ref{5.61}). This is the unique couple in $\mathcal{C}_q(\mathbb{R}^N)$ satisfying $\xn=\mathbb{T}(\tau,F).$
}
The solution $v$ is
 $\xs$-moderate so that $\tau\in\mathfrak{M}_q(\mathbb{R}^N).$

We claim that $u:=v\oplus U_F$ is a solution with initial trace $\mathrm{tr}^c(u)=(\tau,F).$ Indeed $u\geq v$ so that $\mathcal{R}_q(u)\subset\mathcal{R}_q(v).$ On the other hand since $\tau$ is $\frak T_q$-locally finite in $\mathcal{R}_q(v)=F^c,$ it follows that $\mathcal{S}_q(u)\subset F.$ Thus $\mathcal{R}_q(v)\subset\mathcal{R}_q(u)$ and we conclude that $\mathcal{R}_q(u)=\mathcal{R}_q(v)$ and $F=\mathcal{S}_q(u).$ This also implies that $v=u_{reg}.$

Finally
$$\mathcal{S}_q(u)=\mathcal{S}_q(v)\bigcup b_q(\mathcal{S}_q(U_F))=b_q(F)\bigcup F_\tau,$$
so that $F$ satisfies (\ref{5.61}).

The fact that, for $\xn\in\mathfrak{M}_q(\mathbb{R}^N),$ the couple $(\tau,F)$ defined by (\ref{5.71}) is the only one in $\mathcal{C}_q(\mathbb{R}^N)$ satisfying $\xn=\mathbb{T}(\tau,F)$ follows immediately from the definition of these spaces.\

At end, statements \textbf{A-D} imply (i)-(iv).\hfill$\Box$\\

\noindent \Remark If $\xn\in\mathfrak{M}_q(\mathbb{R}^N)$ then $G$ and $v$ as defined by (\ref{5.63}) have the following alternative representation:
\be
G:=\bigcup_{\mathcal{E}_\xn}E=\bigcup_{\mathcal{F}_\xn}Q,\;\;\;\;v:=\sup\{u_{\xn}\chi_Q:\;Q\in\mathcal{F}_\xn\},\label{5.72}
\ee
\be
\mathcal{E}_\xn:=\{Q:\;E\;\frak T_q\text{-open},\;\xn(Q\cap K)<\infty,\;\forall\;\mathrm{compact}\;K\subset\mathbb{R}^N\}.\label{5.73}
\ee

To verify this remark we first observe that Lemma \ref{6} implies that if $A$ is a $\frak T_q$-open set then there exists an increasing sequence of $\frak T_q$-quasi closed sets $\{E_n\}$ such that $A=\cup_{n=1}^\infty E_n.$ In fact, in the notation of Lemma \ref{6} (II)(i)-(ii), we may choose $E_n=F_n\setminus L$ where $L=A'\setminus A,$ is a set of capacity zero.

Therefore
$$\bigcup_{\mathcal{D}_\xn}D\subset\bigcup_{\mathcal{F}_\xn}Q\subset\bigcup_{\mathcal{E}_\xn}E:=H.$$
On the other hand, if $E\in\mathcal{E}_\xn$ then $\xm_{\mathcal{R}_q}(u)(K\cap\widetilde{E})=\xm_{\mathcal{R}_q}(u)(K\cap E)=\xn(K\cap E)<\infty,$ for any compact $K\subset\mathbb{R}^N.$ By Proposition \ref{24}-(ix), $\widetilde{E}$ is regular, i.e., there exists a $\frak T_q$-open regular set $Q$ such that $E\qsub Q.$ Thus $H=\bigcup_{\mathcal{D}_\xn}D.$

If $D$ is a $\frak T_q$-open regular set then $D=\cup_{n=1}^\infty E_n,$ where $\{E_n\}$ is an increasing sequence of $\frak T_q$-quasi closed sets. We infer
$$u_{\xn\chi_D}=\lim_{n\to\infty} u_{\xn\chi_{E_n}}.$$
Therefore
$$\sup\{u_{\xn\chi_{Q}}:\;Q\in\mathcal{D}_\xn\}\leq\sup\{u_{\xn\chi_{Q}}:\;Q\in\mathcal{F}_\xn\}\leq\sup\{u_{\xn\chi_{Q}}:\;Q\in\mathcal{E}_\xn\}.$$
On the other hand, if $E\in\mathcal{E}_\xn$, there exists a $\frak T_q$-open regular set $Q$ such that $E\qsub Q.$ Consequently the equality follows.
\setcounter{equation}{0}

\section{The equation $\prt_t  u-\xD u+Vu=0$}

Let $0<T\leq\infty$, $Q_T:=\mathbb{R}^N\ti (0,T)$, $C>0$ and $V:Q_T\rightarrow [0,\infty)$ be a Borel function satisfying
\begin{equation}\label{sch}0\leq V(x,t)\leq\frac{C}{t}\qquad\forall (x,t)\in Q_T.\end{equation}
In this section we prove a general representation theorem for positive solutions of
\begin{equation}\label{sch1}
\prt_t  u-\xD u+Vu=0\qquad\text{in }Q_T.
\end{equation}
\subsection{Preliminaries}
We recall that $\mathfrak{M}(\mathbb{R}^N)$ is the set of Radon measures on $\mathbb{R}^N$ and $\mathfrak{M}_+(\mathbb{R}^N)$ its positive cone.
\begin{defin}
Let $\xm\in\mathfrak{M}(\mathbb{R}^N)$. We say that $u$ is a weak solution of problem
\begin{equation}\label{prob*1}
\BA{lll}
\prt_t  u-\xD u+Vu=0\qquad&\text{ in }\; Q_T\\ 
\phantom{\xD +,V\chi_{\xO}}
u(.,0)=\xm\qquad&\text{ in }\; \BBR^N,
\EA
\end{equation}
if $u\in L^1_{loc}(\overline Q_T),$ $Vu\in L^1_{loc}(\overline Q_T)$ and there holds
\be
\int\int_{Q_T} u(-\xf_t-\xD\xf)dxdt+\int\int_{Q_T} Vu\xf dxdt=\int_{\mathbb{R}^N}\xf(x,0)d\xm,\label{weakun}
\ee
for all $\xf\in X(Q_T)$, where
$$X(Q_T)=\{\xf\in C_c(\BBR^N\ti[0,T)),\;\xf_t+\xD\xf\in L^\infty_{loc}(\overline{Q}_T)\}.$$
\end{defin}

\noindent \Remark The definition implies that for any $\gz\in C_c^2(\BBR^N)$, the function $t\mapsto \int \gz(x)u(x,t)dx$ can be extended by continuity on $[0,T]$ as a continuous function and
\begin{equation}\label{sch2}
\BA{lll}\displaystyle
\lim_{t\to 0}\myint{\BBR^N}{}\gz(x)u(x,t)dx=\myint{\mathbb{R}^N}{}\gz d\xm.
\EA
\end{equation}
Therefore $\norm {u(.,t)}_{L1(\Gw)}$ remains uniformly bounded on $(0,T)$ for any bounded open set $\Gw\subset\BBR^N$.
\begin{lemma}\label{lem1}
Let $\xm\in\mathfrak{M}_+(\mathbb{R}^N)$ and assume that there exists  a positive weak solution $u$ of problem (\ref{prob*1}) where $V$ satisfies (\ref{sch}). Then for any smooth bounded domain $\xO$ there exists a unique positive weak solution $v$ of problem
\begin{equation}\label{prob*}
\BA{lll}
\prt_t  v-\xD v+Vv=0\qquad&\text{ in }\; Q_T^\xO=\xO\times(0,T)\\ 
\phantom{\prt_t  v-\xD v+V}
v=0\qquad&\text{ on }\; \partial_lQ_T^\xO=\partial\xO\times(0,T)\\
\phantom{,,\xD v+V}
\!v(.,0)=\chi_{\xO}\xm\qquad&\text{ in }\; \xO,
\EA
\end{equation}
where $\chi_{\xO}$ is the characteristic function on $\xO$ and there holds $v\leq u$ in $\xO\times(0,T)$.
\end{lemma}
\Proof
Let $\{t_j\}$ be a decreasing sequence converging to $0$, such that $t_j<T$ for all
$j\in\mathbb{N}.$ We consider the following problem
$$
\BA{lll}
\prt_t  v-\xD v+Vv=0\qquad&\text{ in }\; \xO\times(t_j,T)&\\ 
\phantom{\prt_t  v-\xD v+V}
v=0\qquad&\text{ on }\; \!\partial\xO\times(t_j,T)\\ 
\phantom{,\xD v+V}
v(.,t_j)=u(.,t_j)\qquad&\text{ in }\;\xO.
\EA
$$
Since $u(.,t_j)\in L^1(\Gw)$ and $0\leq V\in L^\infty(\xO\times(t_j,T]),$ there exists a unique positive weak solution $v_j$ of the above problem, smaller than the solution $\BBH^\Gw[u(.,t_j)\chi_{\xO}]$, where $\BBH^\Gw$ is the heat operator in $Q^\Gw:=\Gw\ti(0,\infty)$ with zero boundary condition furthermore $v_j\leq u$ in $\xO\times(t_j,T)$ for all $j\in\mathbb{N}$. By standard parabolic estimates we may assume that the sequence $\{v_j\}$ converges locally uniformly in $\xO\times(0,T]$ to a nonnegative function $v$ smaller than $u$.
If $\xf\in C^{1,1;1}(\overline{Q_T^\xO})$ vanishes on $\partial_lQ_T^\xO$ and satisfies $\xf(x,T)=0,$ the following identity holds
\[
\int_{t_j}^T\int_\xO v_j(-\xf_t-\xD\xf)dxdt+\int_{t_j}^T\int_\xO Vv_j\xf dxdt
+\int_\xO\xf(x,T-t_j)u(x,T-t_j)dx
=\int_\xO\xf(x,0)u(x,t_j)dx,
\]
where, in the above equality, we have take as test function $\xf(.,.-t_j).$ It follows by the dominated convergence theorem, that $v$ is a weak solution of problem (\ref{prob*}).\hfill$\Box$
\begin{lemma}
Assume (\ref{sch}) holds and let $u$ be a positive weak solution of problem (\ref{weakun}) with $\xm\in\mathfrak{M}_+(\mathbb{R}^N).$ Then for any  $(x,t)\in\mathbb{R}^N\times(0,T],$ we have $$\lim_{R\rightarrow\infty}u_R= u,$$
where $\{u_R\}$ is the increasing sequence of the weak solutions of the problem (\ref{prob*}) with $\xO=B_R(0).$ Moreover, the convergence is uniform in any compact subset of $\mathbb{R}^N\times(0,T].$\label{lem2}
\end{lemma}
\Proof    By the maximum principle (see \cite[Remark 2.5]{MV-CPDE}), 
$$u_R\leq u_{R'}\leq u$$
for all $0<R\leq R'$. Thus $u_R\rightarrow w\leq u.$ Also by standard parabolic estimates, this convergence is locally uniformly. Now by dominated convergence theorem, we have that $w$ is a weak solution of problem (\ref{prob*1}) with initial data $\xm.$
We set $\widetilde{w}=u-w\geq0.$ Then $\widetilde{w}$ satisfies in the weak sense
\bea
\nonumber
\widetilde{w}_t-\xD \widetilde{w}+V\widetilde{w}&=&0
\qquad \text{in}\qquad \mathbb{R}^N\times(0,T),\\ \nonumber
w(x,t)&\geq&0\qquad \text{in}\qquad \mathbb{R}^N\times(0,T)\\ \nonumber
\widetilde{w}(x,0)&=&0\qquad \text{in}\;\;\mathbb{R}^N.
\eea
Since $\widetilde{w}$ satisfies in the weak sense
\bea
\nonumber
\widetilde{w}_t-\xD \widetilde{w}&\leq&0\qquad \text{in } \mathbb{R}^N\times(0,T),\\ \nonumber
w(x,t)&\geq&0\qquad \text{in } \mathbb{R}^N\times(0,T),
\\ \nonumber
\widetilde{w}(x,0)&=&0\qquad \text{in }\mathbb{R}^N,
\eea
We extend $\widetilde{w}$ by $0$ for $t\leq 0$, with the same notation and set $\widetilde{w}_n:=\widetilde{w}\ast J_{\ge_n}$ where $\{J_{\ge_n}\}$ is a sequence of mollifiers  in $\BBR^{N+1}$. Then $\widetilde{w}_n\leq 0$, therefore $\widetilde{w}=0.$ \hfill$\Box$
\begin{lemma}
Let $u\in C^{2,1}(\mathbb{R}^N\times(0,T])$ be a positive solution of
$$\prt_t  u-\xD u+Vu=0\qquad \mathrm{in}\qquad \mathbb{R}^N\times(0,T).$$
Assume that, for any $x\in\mathbb{R}^N,$ there exists an open bounded neighborhood $U$ of $x$ such that
$$\int_0^T\int_{U}u(y,t)V(y,t)dxdt<\infty$$ Then $u\in L^1(U\times(0,T))$ and there exists a unique positive Radon measure $\xm$ such that
$$\lim_{t\rightarrow0}\int_{\mathbb{R}^N}u(y,t)\xf(x)dx=\int_{\mathbb{R}^N}\xf(x)d\xm\qquad\forall\xf\in C_0^\infty(\mathbb{R}^N).$$\label{lem3}
\end{lemma}
\Proof
Since $Vu\in L^1(U\times(0,T))$ the following problem has a weak solution $v$ (see \cite{MV-CPDE}).
\bea
\nonumber
\prt_t  v-\xD v&=&Vu,\qquad \mathrm{in}\qquad U\times(0,T],\\ \nonumber
v(x,t)&=&0\qquad\mathrm{on}\;\; \partial U\times(0,T]\\ \nonumber
v(x,0)&=&0\qquad \mathrm{in}\;\;U.
\eea
Thus the function $w=u+v$ is a positive solution of the heat equation, thus there exists a unique Radon measure $\xm$ such that
$$\lim_{t\rightarrow0}\int_{U}w(y,t)\xf(x)dx=\int_{U}\xf(x)d\xm,\;\;\;\;\forall\xf\in C_0^\infty(U).$$
But the initial data of $v$ is zero, thus the result follows by a partition of unity and Lemma \ref{lem2}.\hfill$\Box$

\subsection{Representation formula for positive solutions}
Assume $V$ satisfies (\ref{sch}) in $Q_T$
 and let $u$ be a positive solution of (\ref{sch1}). If $\psi\in C^{2,1}( \mathbb{R}^N\times(0,T])$, we set
 $u(x,t)=e^{\psi(x,t)} v(x,t).$
Then $v$ satisfies
\be
\prt_t  v-\xD v-2\nabla\psi\nabla v-|\nabla\psi|^2v-2\xD\psi v+\left(\psi_t+\xD\psi+V\right) v=0\;\;\mathrm{in}\qquad \mathbb{R}^N\times(0,T].\label{eq1}
\ee
We choose $\psi$ to be the solution of the problem
\bea
\nonumber
-\psi_t-\xD\psi&=&V\qquad\text{in}\qquad \mathbb{R}^N\times(0,T]\\ \nonumber
\psi(x,T)&=&0\qquad\text{ in}\qquad \mathbb{R}^N.
\eea
Then 
\begin{equation}\label{rep}
\psi(t,x)=\int_t^T\int_{\mathbb{R}^N}\frac{1}{\left(4\pi(s-t)\right)^\frac{n}{2}}e^{-\frac{|x-y|^2}{4(s-t)}}V(x,s)dxds.
\end{equation}
By a straightforward calculation we verify that\smallskip

\noindent 1. $0\leq \psi\leq C\ln\frac{T}{t},$\smallskip

\noindent 2. $|\nabla\psi|\leq C_1(T)+C_2(\ln\frac{T}{t}).$\smallskip

\noindent Thus (\ref{eq1}) becomes
$$\prt_t  v-\xD v-\sum_{i=1}^n\left(2\psi_{x_i}v\right)_{x_i}-|\nabla\psi|^2v=0.$$
Since $\int_0^1|\ln t|^pdt<\infty$ for all $p\geq1$, we verify by $1$ and $2$ that
$$\int_0^T\sup_{x\in\mathbb{R}^N}|\psi|^qds<M_1<\infty\qquad\forall q\geq1$$
and
$$\int_0^T\sup_{x\in\mathbb{R}^N}|\nabla\psi|^qds<M_2<\infty\qquad\forall q\geq1.$$
For $A_{i,j}=\xd_{ij},$  $A_{i}=2\psi_{x_i}$ $B_i=0$ and $C=|\nabla\psi|^2$ we see that the above operator satisfies the condition $H$ in \cite{arons} for $R_0=\infty$ and $p=\infty.$ Thus there exists a kernel $\xG(x,t;y,s)$, defined in $Q_T\ti Q_T$ satisfying the estimates
\be
C_1(T,n,M_2)\frac{1}{\left(4\pi(t-s)\right)^\frac{n}{2}}e^{-A_1\frac{|x-y|^2}{4(t-s)}}\leq\xG(x,t;y,s)\leq
C_2(T,n,M_2)\frac{1}{\left(4\pi(t-s)\right)^\frac{n}{2}}e^{-A_2\frac{|x-y|^2}{4(t-s)}},\label{heatestimates}
\ee
where $A_1,\;A_2>0$ depend on $T,\;n,\;M_2$ with the property that $v$ admits the following representation formula:
\begin{equation}\label{rep1}
v(x,t)=\int_{\mathbb{R}^N}\xG(x,t;y,0)d\xm(y),
\end{equation}
where $\xm$ is a uniquely defined positive Radon measure on $\mathbb{R}^N$, and there holds
$$\lim_{t\rightarrow0}\int_{\mathbb{R}^N}\int_{\mathbb{R}^N}\xG(x,t;y,0)\xf(x)d\xm(y)dx=\int_{\mathbb{R}^N}\xf d\xm\qquad\forall\xf\in C_0^\infty(\mathbb{R}^N).$$
\\
Furthermore, if $e^{-\xg|x|^2}u_0\in L^2(\mathbb{R}^N)$ for some $\xg\geq0$, and if $u_0$ is continuous at $y$, then
\be
\lim_{t\rightarrow0}\int_{\mathbb{R}^N}\xG(x,t;y,0)u_0(x)dx=u_0(y).\label{limit}
\ee
Finally we have
\be
u(x,t)=e^{\psi(x,t)}\int_{\mathbb{R}^N}\xG(x,t;y,0)d\xm(y).\label{repres}
\ee

\setcounter{equation}{0}
\section{$\xs$-moderate solutions}
\subsection{Preliminaries}
\begin{prop}
Let $u\in\mathcal{U}_+(Q_T).$ Then
\be
\max(u_{\mathcal{R}_q},[u]_{\mathcal{S}_q(u)})\leq u\leq u_{reg}+[u]_{\mathcal{S}_q(u)}.\label{s4.6}
\ee
\end{prop}
\Proof    The principle of the proof is similar as the one in \cite{M-JAM}.\\
By Proposition \ref{24}-(vii), the function $v=u\ominus u_{reg}$ vanishes on $\mathcal{R}_q(u)$ i.e., $\frak T_q\text{-supp}\,(v)\subset\mathcal{S}_q(u).$ Thus $v$ is a solution dominated by $u$ and supported in $\mathcal{S}_q(u),$ which implies that $v\leq[u]_{\mathcal{S}_q(u)}$ by Definition \ref{32}. Since $u-u_{reg}\leq v$ this implies the inequality on the right hand side of (\ref{s4.6}). The inequality on the left hand side is obvious.\hfill$\Box$
\begin{prop}
Let $u\in\mathcal{U}_+(Q_T)$ and let $A,\;B$ be two disjoint $\frak T_q$-closed subsets of $\mathbb{R}^N.$ If $\frak T_q\text{-supp}\,(u)\subset A\cup B$ and $[u]_A,\;[u]_B$ are $\xs$-moderate then $u$ is $\xs$-moderate. Furthermore
\be
u=[u]_A\oplus[u]_B=[u]_A\vee[u]_B.\label{s4.7}
\ee\label{37}
\end{prop}
\Proof  The proof is same as in \cite{M-JAM}.\\
By Proposition \ref{27}-(iii) there exist two increasing sequences $\{\tau_n\},\;\{\tau_n'\}\subset W^{-\frac{2}{q},q}(\mathbb{R}^N)\cap\mathfrak{M}^b_+(\mathbb{R}^N)$ such that
$$u_{\tau_n}\uparrow[u]_A,\qquad u_{\tau_n'}\uparrow[u]_B.$$
By proposition \ref{33}, $\frak T_q$-{supp}$\,(\tau_n)\qsub A$ and $\frak T_q$-{supp}$\,(\tau_n')\qsub B.$ Thus $\qcap\left(\frak T_q\text{-supp}(\tau_n)\cap \frak T_q\text{-supp}(\tau_n')\right)=0,$ and
$$u_{\tau_n}\vee u_{\tau_n'}=u_{t_n}\oplus u_{t_n'}=u_{\tau_n+\tau_n'}.$$
By (\ref{13}) and Definition \ref{32},
\be
\max([u]_A,[u]_B)\leq u\leq[u]_A+[u]_B.\label{s4.8}
\ee
Therefore,
$$\max(u_{\tau_n},u_{\tau_n}')\leq u\Rightarrow u_{\tau_n+\tau_n'}\leq u.$$
On the other hand
$$u-u_{\tau_n+\tau_n'}\leq[u]_A-u_{\tau_n}+[u]_B-u_{\tau_n'}\downarrow0.$$
Thus
\be
\lim_{n\to\infty} u_{\tau_n+\tau_n'}=u,\label{s4.9}
\ee
so that $u$ is $\xs$- moderate.\\
The assertion (\ref{s4.7}) is equivalent to the statements: (a) $u$ is the largest solution dominated by $[u]_A+[u]_B$ and (b) $u$ is the smallest solution dominating $\max([u]_A,[u]_B).$
Let
$$u\leq w:=[u]_A\oplus[u]_B\leq[u]_A+[u]_B.$$
Thus we have $[u]_A\leq[w]_A.$
But $[w]_A\leq w\leq[u]_A+[u]_B\;\Rightarrow\;[w]_A-[u]_A\leq[u]_B.$
By Notation \ref{34} we have
$$v=[([w]_A-[u]_A)_+]_\dag\leq[u]_B,\qquad v\leq[w]_A,$$
that is
$$\frak T_q\text{-supp}\,(v)\subset A\qquad\mathrm{and}\qquad\frak T_q\text{-supp}\,(v)\subset B.$$
But $A\cap B=\emptyset,$ which implies $v=0$ and $[w]_A\leq[u]_A.$ Similarly, we have $[w]_B\leq[u]_B.$ Thus
$$[w]_A=[u]_A,\qquad[w]_B\leq[u]_B.$$
By (\ref{s4.8}) and the fact that for any Borel $E$ $[u]_E\leq[u]_{\widetilde{E}\cap A}+[u]_{\widetilde{E}\cap B},$ there holds
$$\mathcal{S}_q(u)=\mathcal{S}_q(w).$$
Let $Q$ be a $\frak T_q$-open regular set in $\mathcal{R}_q(w),$ then $Q\in\mathcal{R}_q(u).$ Using (\ref{13}), (\ref{14}) and the fact that $\frak T_q\text{-supp}\,(w)\subset A\cap B,$ we derive
$$[w]_{Q}\leq[w]_{\widetilde{Q}\cap A}+[w]_{\widetilde{E}\cap B}=[[w]_A]_{\widetilde{Q}}+[[w]_B]_{\widetilde{Q}}=
[u]_{\widetilde{Q}\cap A}+[u]_{\widetilde{Q}\cap B}.$$
Since $[w]_Q,\;[u]_Q$ are moderate solutions and $A\cap B=\emptyset,$ we have $[u]_{\widetilde{Q}\cap A}\oplus[u]_{\widetilde{Q}\cap B}\leq[u]_Q,$ which implies
$[w]_Q=[u]_Q.$ Thus by Proposition \ref{24}-(ii) $w_{\mathcal{R}_q}=u_{\mathcal{R}_q},$ and since $u$ is $\xs$-moderate by Proposition \ref{35} and the remark below we get
$$u\leq w\leq u_{\mathcal{R}_q}\oplus U_F.$$
By the uniqueness of $\xs$- moderate solutions (Theorem \ref{27}-(iv)), $w$ and $u$ coincide. This proves (a).\\
For the statement (b), we note that
$$u_{\tau_n+\tau_n'}=u_{\tau_n}\vee u_{\tau_n'}\leq[u]_A\vee[u]_B,$$
since $u_{\tau_n}\leq[u]_A$ and $u_{\tau_n'}\leq[u]_B.$ Thus the result follows by (\ref{s4.9}) and (\ref{s4.8}), by letting  $n$ tend to infinity.\hfill$\Box$
\subsection{Characterization of positive solutions of $\prt_t  u-\xD u+u^q=0$}
The following notation is used throughout the subsection.\\
Let $u\in\mathcal{U}_+(Q_T).$ Set
$$V=u^{q-1},$$
then $$V\leq\left(\frac{1}{q-1}\right)^{q-1}t^{-1}.$$
Thus $u\in C^{2,1}(\mathbb{R}^N\times(0,T]$ and satisfies
\be
\prt_t  u-\xD u+Vu=0,\qquad \mathrm{in}\qquad \mathbb{R}^N\times(0,1].\label{LV1}
\ee
Hence, by the representation formula (\ref{repres}), $u$ satisfies
\begin{equation}\label{repres2}
u(x,t)=e^\psi\int_{\mathbb{R}^N}\xG(x,t;y,0)d\xm(y),\qquad\forall\;t\leq T,
\end{equation}
where $\xm$ is Radon measure (see subsection 7.2). The measure $\gm$ is called the {\it extended initial trace} of $u$.
\\
For any Borel set $E$ set
$$\xm_E=\xm\chi_E\qquad\mathrm{and}\qquad(u)_E=e^\psi\int_{\mathbb{R}^N}\xG(x,t;y,0)d\xm_E,\qquad\forall\;t\leq T.$$
\begin{lemma}
Let $F$ be a compact subset of $\mathbb{R}^N.$ Then
$$
(u)_F\leq[u]_F,\qquad\forall t\leq T.
$$\label{lemma5.3}
\end{lemma}
\Proof    We follow the ideas of \cite{M-JAM}, adapted to the parabolic framework.\\
Let $A$ be a Borel subset of $\mathbb{R}^N.$ Let $0<\xb\leq\frac{T}{2}$ and let $v_\xb^A$ be the positive solution of
\begin{equation}\label{LV}
\begin{array} {lll}
\prt_t  v-\xD v+Vv=0\qquad&\text{in }\mathbb{R}^N\times(\xb,T] \\ 
\phantom{--.,,-}
v(.,\xb)=u(.,\xb)\chi_A(.)\qquad&\text{in }\mathbb{R}^N.
\end{array}
\end{equation}
Also let $w_\xb^A$ be the positive solution of
$$\begin{array} {lll}
\phantom{,,...,-}
\prt_t  w-\xD w+|w|^{q-1}w=0\qquad\phantom{,,...,-}&\text{in }\mathbb{R}^N\times(\xb,T]\\ 
\phantom{--,.-----..,-}
\!w(.,\xb)= \chi_A(.)u(.,\xb)\qquad&\text{in }\mathbb{R}^N.
\end{array}
$$
Then by the maximum principle $w_\xb^A\leq u$, which implies
$$0=\frac{dw_\xb^A}{dt}-\xD w_\xb^A+(w_\xb^A)^q\leq\frac{dw_\xb^A}{dt}-\xD w_\xb^A+Vw_\xb^A.$$
Thus $w_\xb^A$ is a supersolution of (\ref{LV}), and by the maximum principle (see \cite{arons} or Lemma \ref{lem2}), we have
$$v_\xb^A\leq w_\xb^A\leq u.$$
For any sequence $\{\xb_k\}$ decreasing to zero one can extract a subsequence $\{\xb_{k_n}\}$ such that $\{w_{\xb_{k_n}}^A\}$ and $\{v_{\xb_{k_n}}^A\}$
converge locally uniformly; we denote the limits $w^A$ and $v^A$ respectively (the limits may depend on the sequence). Then $w^A\in\mathcal{U}_+(Q_T)$ while $v^A$ is a solution of (\ref{LV1}), and
\be
v^A\leq w^A\leq[u]_{\widetilde{Q}},\qquad\forall Q\;\mathrm{open,}\;A\subset Q.\label{s5.3}
\ee
The second inequality follows from the fact that $\frak T_q\text{-supp}\,(w_\xb^A)\subset\widetilde{Q}$ for any $\xb.$

Now we set $v_{\xb_{k_n}}^A=e^\psi\widetilde{v}_n,$ where $\psi$ is the function in subsection 7.2. Then $\widetilde{v}_n$ is the solution of
\bea
\nonumber
\prt_t  v-\xD v-2\nabla\psi\nabla v-|\nabla\psi|^2v-2\xD\psi v+\left(\psi_t+\xD\psi+V\right) v=0\;\;&\mathrm{in}&\quad \mathbb{R}^N\times(\xb_{k_n},T].\\ \nonumber
v(.,\xb_{k_n})=\chi_A(.)\int_{\mathbb{R}^N}\xG(.,{\xb_{k_n}};y,0)d\xm(y)\;\;&\mathrm{ in}&\quad\mathbb{R}^N.
\eea
Now using the representation formula in \cite{arons}, we derive that for any open $Q\supset A$, there holds
\bea
\nonumber
\widetilde{v}_n(x,t)&=&\int_{\mathbb{R}^N}\chi_A(x)\xG(x,t-{\xb_{k_n}};y,0)\left(\int_{\mathbb{R}^N}\xG(x,{\xb_{k_n}};y,0)d\xm(y)\right)dx\\ \nonumber
&=&\int_{\mathbb{R}^N}\left(\int_{\mathbb{R}^N}\chi_A(x)\xG(x,t-{\xb_{k_n}};y,0)\xG(x,{\xb_{k_n}};y,0)dx\right)dxd\xm(y)\\ \nonumber
&\leq&\int_{\mathbb{R}^N}\left(\int_{\mathbb{R}^N}\chi_Q(x)\xG(x,t-{\xb_{k_n}};y,0)\xG(x,{\xb_{k_n}};y,0)dx\right)dxd\xm(y).
\eea
Therefore, by (\ref{limit}), estimate (\ref{heatestimates}) and using the fact that $\xG(x,t-s;y,0)$ is a continuous function for any $s<t$ (see \cite{arons}), we can let $k_n\to\infty$ in the above inequality and get
$$\lim_{n\to\infty}\widetilde{v}_n\leq\int_{\mathbb{R}^N}\xG(x,t;y,0)d\xm_{\widetilde{Q}}.$$
Hence
$$v^A\leq(u)_{\widetilde{Q}}.$$

We apply the same procedure to the set $A^c$ extracting a further subsequence of $\{\xb_{k_n}\}$ in order to obtain the limits $v^{A^c}$ and $w^{A^c}.$ Thus
$$v^{A^c}\leq w^{A^c}\leq[u]_{\widetilde{Q'}},\qquad\forall Q'\;\mathrm{open,}\;A^c\subset Q'.$$
Note that
$$v^{A}+v^{A^c}=u,\;\;\;v^{A}\leq(u)_{\widetilde{Q}},\;\;\;v^{A^c}\leq(u)_{\widetilde{Q'}}.$$
Therefore
\be
v^{A}=u-v^{A^c}\geq(u)_{(\widetilde{Q'})^c}.\label{s5.4}
\ee
Now, given $F$ compact, let $A$ be a closed set and $O$ an open set such that $F\subset O\subset A.$ Note that $A^c\cap F=\emptyset.$ By (\ref{s5.4}) with $Q'=A^c$
$$v^A\geq(u)_O.$$
By (\ref{s5.3})
$$
v^A\leq w^A\leq[u]_{\widetilde{Q}}\qquad\forall Q\;\mathrm{open,}\;A\subset Q,
$$
and consequently
\be
(u)_F\leq(u)_O\leq[u]_{\widetilde{Q}}.\label{s5.5}
\ee
By Lemma \ref{6}, we can choose a sequence of open sets $\{Q_n\}$ such that $\cap\widetilde{ Q}_j=E'\qeq F,$ thus by Proposition \ref{17}-(iii) $[u]_{Q_j}\downarrow[u]_F.$ The result follows by (\ref{s5.5}).\hfill$\Box$\medskip

In the next lemma we prove that the extended initial trace of a positive solution of (\ref{maineq}) is absolutely continuous with respect to the $\qcap$-capacity.

\begin{lemma}
Let $u\in \CU_+(Q_T)$,  $\gm$ its extended initial trace. If $E$ is a Borel set and $\qcap(E)=0$ then $\xm(E)=0.$\label{lemma5.4}
\end{lemma}
\Proof    The proof is similar as the one in \cite{M-JAM}.
If $F$ is a compact subset of $E$, then $\qcap(F)=0$ and therefore by Proposition \ref{sygklisi1}, $U_F=0.$ But $[u]_F=u\wedge U_F=0.$ Therefore, by Lemma \ref{lemma5.3} $(u)_F=0.$ Consequently $\xm(F)=0.$ As this holds for every compact subset of $E$ we conclude that $\xm(E)=0.$\hfill$\Box$\medskip

We recall that, if $\xn\in W^{-\frac{2}{q},q}(\mathbb{R}^N)\cap\mathfrak{M}^b_+(\mathbb{R}^N),$ then for any $T>0$, there exists a constant $C>0$ independent on $\xn$ (see Lemma 2.11-\cite{MV-CVPDE}) such that
\be
C^{-1}||\xn||_{ W^{-\frac{2}{q},q}(\mathbb{R}^N)}\leq||\BBH[\xn]||_{L^q(Q_T)}\leq C||\xn||_{ W^{-\frac{2}{q},q}(\mathbb{R}^N)},\label{caloric}
\ee
where $\BBH[\xn]$ is the solution of the heat equation in $Q_\infty$ with $\xn$ as initial data.
\begin{lemma}
Let $u\in \CU_+(Q_T)$,  $\gm$ its extended initial trace and $\xn\in W^{-\frac{2}{q},q}(\mathbb{R}^N)\cap\mathfrak{M}^b_+(\mathbb{R}^N).$ Suppose that there exists no positive solution of (\ref{maineq}) dominated by the supersolution $v=\inf\{u,\BBH[\xn]\}.$ Then $\xm\perp\xn.$\label{lemma5.5}
\end{lemma}
\Proof 
Set $V'=v^{q-1}$, then $v$ is a supersolution of
\be
\prt_t  w-\xD w+V'w=0\qquad \mathrm{in}\qquad \mathbb{R}^N\times(0,T].\label{s1}
\ee
We first claim that there exists no positive solution of the above problem dominated by $v.$ We proceed by contradiction in assuming that there exists a positive solution $w\leq v$ of (\ref{s1}). Then $w$ satisfies
$$\prt_t  w-\xD w+w^q\leq \prt_t  w-\xD w+V'w=0.$$
Since 
$$\norm v_{L^q(Q_T)}\leq \norm {\BBH[\xn]}_{L^q(Q_T)}\approx \norm {\xn}_{W^{-\frac{2}{q},q}(\mathbb{R}^N)}$$
this implies that $w$ is a positive moderate solution of (\ref{maineq}) dominated by $v,$ contrary to assumption.
Now for any $t\leq T,$ we have by representation formula (\ref{repres}),
\bea
\nonumber
\inf\{u,\BBH[\xn]\}&=&\inf\left\{e^\psi\int_{\mathbb{R}^N}\xG(x,t;y,0)d\xm(y),\BBH[\xn]\right\}\\ \nonumber
&\geq&\inf\left\{\int_{\mathbb{R}^N}\xG(x,t;y,0)d\xm(y),\BBH[\xn]\right\}\\ \nonumber
&\geq& C\inf\left\{\BBH[\xm](\frac{t}{A_2},x),\BBH[\xn](t,x)\right\}\\ \nonumber
&\geq& C\inf\left\{\BBH[\xm](\frac{t}{\max(A_2,1)},x),\BBH[\xn](\frac{t}{\max(A_2,1)},x)\right\},
\eea
where, in the above inequalities, we have used estimates (\ref{heatestimates}) and the constants $C>0, \;A_2>0$ therein.\\
Now since $\inf\left\{\BBH[\xm](\frac{t}{\max(A_2,1)},x),\BBH[\xn](\frac{t}{\max(A_2,1)},x)\right\}$ is a supersolution of $\prt_t  w-\frac{1}{\max(A_2,1)}\xD w=0,$ there exists a positive Radon measure $\widetilde{\xn}$ such that
$$\lim_{t\rightarrow0}\int_{\mathbb{R}^N}\xf(x)\inf\left\{\BBH[\xm](\frac{t}{\max(A_2,1)},x),\BBH[\xn](\frac{t}{\max(A_2,1)},x)\right\}dx=
\int_{\mathbb{R}^N}\xf(x)d\widetilde{\xn}\quad\forall \xf\in C_0^\infty(\mathbb{R}^N).$$
Thus in view of Lemmas \ref{lem2} and \ref{lem3}, there exists a positive weak solution $\widetilde{v}\leq v$ of the problem
\bea
\nonumber
\prt_t  w-\xD w+V'w&=&0\qquad \mathrm{in}\qquad \mathbb{R}^N\times(0,T].\\ \nonumber
w(.,0)&=&\widetilde{\xn}\qquad \mathrm{in}\qquad \mathbb{R}^N,
\eea
and by the first claim it yields $\widetilde{\xn}=0.$

By the Lebesgue-Radon-Nikodym Theorem we can write $d\xn=\xf d\xm+d\xs,$ where $0\leq\xf\in L^1_{loc}(\mathbb{R}^N,\xm)$ and $\xs\perp\xm.$
Thus we have
\bea
\nonumber
0&=&\lim_{t\rightarrow0}\int_{\mathbb{R}^N}\xf(x)\inf\left\{\BBH[\xm](\frac{t}{\max(A_2,1)},x),\BBH[\xn](\frac{t}{\max(A_2,1)},x)\right\}dx\\ \nonumber
&\geq&\lim_{t\rightarrow0}\int_{\mathbb{R}^N}\xf(x)h(\frac{t}{\max(A_2,1)},x,y)\min\{f,1\}(y)d\xm(y)dx\\ \nonumber
&=&\lim_{t\rightarrow0}\int_{\mathbb{R}^N}\xf(y)\min\{f,1\}(y)d\xm(y)=0,
\eea
where, we recall it, $h(t,x,y)$ in the heat kernel in $Q_\infty$. Hence $f=0$ and $\xn\perp\xm.$\hfill$\Box$
\begin{lemma}
Let $u\in \CU_+(Q_T)$,  $\gm$ its extended initial trace and suppose that for every $\xn\in \mathfrak{M}^b_+(\mathbb{R}^N)\cap W^{-\frac{2}{q},q}(\mathbb{R}^N)$ there exists no positive solution of (\ref{maineq}) dominated by $v=\inf(u,\BBH[\xn]).$ Then $u=0.$\label{lemma5.6}
\end{lemma}
\Proof    The proof is similar  as the one in \cite{M-JAM}.
By Lemma \ref{lemma5.5},
$$\xm\perp\xn\qquad\forall\xn\in W^{-\frac{2}{q},q}(\mathbb{R}^N)\cap\mathfrak{M}^b_+(\mathbb{R}^N).
$$
Suppose that $\xm\neq0.$ By Lemma \ref{lemma5.4}, $\xm$ vanishes on sets of $\qcap$-capacity zero. Thus, there exists an increasing sequence $\{\xn_k\}\subset W^{-\frac{2}{q},q}(\mathbb{R}^N)\cap\mathfrak{M}^b_+(\mathbb{R}^N)$ which converges to $\xm.$ Thus $\xm\perp\xn_k$ and for every $k\in\BBN$ there exists a Borel set $A_k\subset\mathbb{R}^N$ such that
$$\xm(A_k)=0\;\text{ and }\;\xn_k(A_k^c)=0.$$
Therefore, if $A=\cup_k A_k$ then
$$\xm(A)=0\;\text{ and }\;\xn_k(A^c)=0\quad\forall k.$$
Since $\xn_k\leq\xm$ we have $\xn_k(A)=0$ and therefore $\xn_k=0.$ Contradiction.\hfill$\Box$
\begin{lemma}
Let $u\in\CU_+(Q_T)$. Then $[u]_{\mathcal{S}_q(u)}$ is $\xs$-moderate.\label{lemma5.7}
\end{lemma}
\Proof
To simplify the notations we put $u_\CS=[u]_{\mathcal{S}_q(u)}$ and denote $F:=\frak T_q\text{-supp}\,(u_\CS).$ Incidentally, $F\subset\mathcal{S}_q(u);$ since if $\mathcal{S}_q(u)$ is thin at $\xi$, then $\mathcal{S}_q(u)^c\cup\{\xi\}$ is $\frak T_q$-open and $\mathcal{S}_q(u)^c\cup\{\xi\}\qeq\mathcal{S}_q(u)^c.$ Thus by definition of $F,$ we see that $F$ consists precisely of the $\qcap$-thick points of $\mathcal{S}_q(u).$ The set $\mathcal{S}_q(u)\setminus F$ is contained in the singular set of $u_{\mathcal{R}_q}.$

For $\xn\in W^{-\frac{2}{q},q}(\mathbb{R}^N)\cap\mathfrak{M}^b_+(\mathbb{R}^N)$ we denote by $u_\xn$ the solution of (\ref{maineq}) with initial trace $\xn.$ Put
\be
u^*:=\sup\{u_\xn:\;\xn\in W^{-\frac{2}{q},q}(\mathbb{R}^N)\cap\mathfrak{M}^b_+(\mathbb{R}^N),\;u_\xn\leq u_\CS\}.\label{s5.11}
\ee
By Lemma \ref{lemma5.6} the family over which the supremum is taken is not empty. Therefore $u^*$ is a positive solution of (\ref{maineq}) and, by Proposition \ref{27}-(iii), it is $\xs$-moderate, thus it is the largest $\xs$-moderate solution dominated by $u_\CS$. We denote by $\{\gn_n\}\subset W^{-\frac{2}{q},q}(\mathbb{R}^N)\cap\mathfrak{M}^b_+(\mathbb{R}^N)$ an increasing sequence such that $u^*=\lim_n{n\to\infty}u_{\gn_n}$.

Let $F^*=\frak T_q\text{-supp}\,( u^*).$ Then $F^*$ is $\frak T_q$-closed and $F^*\subset F.$ Suppose that
$$\qcap(F\setminus F^*)>0.$$
Then there exists a compact set $E\subset F\setminus F^*$ such that $\qcap(E)>0$ and $(F^*)^c=:Q^*$ is a $\frak T_q$-open set containing $E.$ Furthermore by Lemma \ref{q-open} there exists a $\frak T_q$-open set $Q'$ such that $E\qsub Q\subset\widetilde{Q'}\qsub Q^*.$ Since $Q'\qsub\frak T_q\text{-supp}\,(u_\CS),$ $[u_\CS]_{Q'}>0$ and therefore by Lemma \ref{lemma5.6}, there exists a positive measure $\tau\in W^{-\frac{2}{q},q}(\mathbb{R}^N)\cap\mathfrak{M}^b_+(\mathbb{R}^N)$ supported in $\widetilde{Q'}$ such that $u_\tau\leq u_\CS.$ As the $\frak T_q$-{supp}$\,(\tau)$ is a $\frak T_q$-closed set disjoint from $F^*$, it follows that the inequality $u^*\geq u_\tau$ does not hold. On the other hand, since $\tau\in W^{-\frac{2}{q},q}(\mathbb{R}^N)\cap\mathfrak{M}^b_+(\mathbb{R}^N)$ and $u_\tau\leq u_\CS,$ it follows that $u_\tau\leq u^*.$ This contradiction shows that
\be
\qcap(F\setminus F^*)=0.\label{s5.12}
\ee

Since $u_{\gn_n}\uparrow u^*$, $\frak T_q\text{-supp}\,(u_{\gn_n})\subset \frak T_q\text{-supp}\,(u^*):=F^*$. Therefore
there exists a $\frak T_q$-closed set $F_0^*\subset F^*$ such that $\mathcal{S}_q(u^*)=F_0^*$  and $\mathcal{R}_q(u^*)=(F^*_0)^c.$ Suppose that
$$\qcap(F\setminus F_0^*)>0.$$
Let $Q'$ be a $\frak T_q$-open subset of $\mathcal{R}_q(u^*)$ such that $[u_\CS]_{Q'}$ is a moderate solution, then $\widetilde{Q'}\qsub\mathcal{R}_q(u^*)$ and $[u^*]_{\widetilde{Q'}}$ is a moderate solution of (\ref{maineq}), i.e.,
$$
\int_0^T\int_{\mathbb{R}^N}[u^*]_{\widetilde{Q'}}^q\xf(x)dxdt<\infty\qquad\forall\xf\in C_0(\mathbb{R}^N),\,\xf\geq 0.
$$
On the other hand $Q'$ is a $\frak T_q$-open subset of $F=\frak T_q\text{-supp}\,(u_\CS);$ therefore the initial trace of $[u^*]_{\widetilde{Q'}}$ has no regular part, i.e.,
$$
\CR_q([u^*]_{\widetilde{Q'}})=0\text{ and }\mathcal{S}_q([u^*]_{\widetilde{Q'}})=\frak T_q\text{-supp}\,([u^*]_{\widetilde{Q'})};
$$
we say that $[u^*]_{\widetilde{Q'}}$ is {\it a purely singular solution of (\ref{maineq})}. It follows that $v:=\left[[u_\CS]_{\widetilde{Q'}}-[u^*]_{\widetilde{Q'}}\right]_\dag$ is a purely singular solution of (\ref{maineq}).

Let $v^*$ be defined as in (\ref{s5.11}) with $u$ replaced by $v.$ Then $v^*$ is a singular $\xs$-moderate solution of (\ref{maineq}). Since $v^*$ is smaller than $u$ and since it is $\xs$-moderate it is dominated by $u^*$. On the other hand, since $v^*$ is singular and $\frak T_q\text{-supp}\,(v^*)\qsub\widetilde{Q'}\qsub\mathcal{R}_q(u^*)$ it follows that $u^*$ is not larger or equal to $v^*,$ i.e. $(v^*-u^*)_+$ is not identically zero.
Since both $u^*$ and $v^*$ are $\xs$-moderate, it implies that there exists $\tau\in W^{-\frac{2}{q},q}(\mathbb{R}^N)\cap\mathfrak{M}^b_+(\mathbb{R}^N)$ such that $u_\tau\leq v^*$, and $(u_\tau-u^*)$ is not identically zero. Therefore $u^*\lneq\max(u^*,u_\tau).$ The function $\max(u^*,u_\tau)$ is a subsolution of (\ref{maineq}) and the smallest solution above it, denoted by $Z$, is strictly larger than $u^*.$ However $u_\tau\leq v^*\leq u^*$ and consequently $Z=u^*.$

This contradiction proves that $\qcap(Q')=0,$ for any set $Q'\subset\mathcal{R}_q(u^*)$ such that $[u]_{Q'}$ is moderate solution, that is $\qcap(\mathcal{R}_q(u^*))=0$ which implies
\be
\qcap(F\setminus F_0^*)=0.\label{s5.13}
\ee
In conclusion, $u^*$ is $\xs$-moderate, $\frak T_q\text{-supp}\,(u^*)\subset F$ and $F_0^*=\mathcal{S}_q(u^*)\qeq F.$ Therefore, by Proposition \ref{35} and the remark below, $u^*=U_F.$ Since by definition $u^*\leq u\leq U_F,$ it follows $u^*=u.$ \hfill$\Box$
\begin{theorem}
Every positive solution of (\ref{maineq}) is $\xs$-moderate.
\end{theorem}
\Proof    We borrow the ideas of the proof to \cite{M-JAM}.
By Proposition \ref{24}-(i), $\mathcal{R}_q(u)$ has regular decomposition $\{Q_n\}$. Furthermore
$$v_n=[u]_{Q_n}\uparrow u_{\mathcal{R}_q}.$$
Thus the solution $u_{\mathcal{R}_q}$ is $\xs$-moderate and
$$u\ominus u_{\mathcal{R}_q}\leq [u]_{\mathcal{S}_q(u)}.$$
Put
$$u_n=v_n\oplus [u]_{\mathcal{S}_q(u)}.$$
By Lemma \ref{lemma5.7} we have that $[u]_{\mathcal{S}_q(u)}$ is $\xs$-moderate solution, thus by Proposition \ref{37}, as $\widetilde{Q}_n\cap\mathcal{S}_q(u)=\emptyset,$ it follows that $u_n$ is $\xs$-moderate. As $\{u_n\}$ is increasing it follows that $\overline{u}=\lim_{n\to\infty} u_n$ is a $\xs$-moderate solution of (\ref{maineq}). In addition
$$v_n\vee[u]_{\mathcal{S}_q(u)}=u_n=v_n\oplus [u]_{\mathcal{S}_q}(u)\Rightarrow\;\;\max(u_{\mathcal{R}_q},[u]_{\mathcal{S}_q(u)})\leq\overline{u}\leq u_{\mathcal{R}_q}+[u]_{\mathcal{S}_q(u)}.$$
This further implies that $\mathcal{S}_q(u)=\mathcal{S}_q(\overline{u}).$ By construction we have $$[u]_{Q_n}=v_n\leq[\overline{u}]_{Q_n}$$
Letting $n\rightarrow\infty$ we have by Proposition \ref{24}
$$u_{\mathcal{R}_q}\leq \overline{u}_{\mathcal{R}_q}\Rightarrow u_{\mathcal{R}_q}=\overline{u}_{\mathcal{R}_q},$$
thus $\mathrm{tr}(u)=\mathrm{tr}(\overline{u})$ and since $\overline{u}\leq u,$ we have by Proposition \ref{35} and the uniqueness of $\xs$-moderate solutions that $u=\overline{u}.$\hfill$\Box$

\end{document}